\documentclass{acta_2024}
\usepackage{amsmath}

\usepackage[longnamesfirst]{natbib} \usepackage{har2nat}
\setcitestyle{comma,aysep={}} 
\shortcites{} 

\usepackage{newtxtext} 
\usepackage{newtxmath} 
\DeclareSymbolFont{CMlargesymbols}{OMX}{cmex}{m}{n} 
\DeclareMathDelimiter{(}{\mathopen} {operators}{"28}{CMlargesymbols}{"00}
\DeclareMathDelimiter{)}{\mathclose}{operators}{"29}{CMlargesymbols}{"01}
\DeclareMathAlphabet\mathcal{OMS}{cmsy}{m}{n} 
\SetMathAlphabet\mathcal{bold}{OMS}{cmsy}{b}{n} 
\usepackage{microtype} 
\usepackage{type1cm} 

\pagerange{\pageref{firstpage}--\pageref{lastpage}}
\doi{XXXXXXXX}  
\usepackage[twoside,
	paperheight=247mm,
	paperwidth=174mm,
	marginratio={3:5,2:3},
	left=18mm,
	top=20mm]{geometry}
    
\usepackage{mathtools}

\numberwithin{figure}{section}
\numberwithin{table}{section}

\usepackage[cmyk]{xcolor}
\usepackage{graphics,graphicx}
\usepackage{caption}
\usepackage{subcaption}
\usepackage{sidecap}
\usepackage{tikz}
\usetikzlibrary{calc,positioning,shapes}
\usetikzlibrary{patterns,decorations.pathmorphing,decorations.markings,decorations.pathreplacing}
\usepackage{pgfplots}
\pgfplotsset{compat=newest}

\usepackage[UKenglish]{babel}
\usepackage{enumerate} 

\newcommand{\ignore}[1]{}
\allowdisplaybreaks

\usepackage{xspace}
\usepackage{bm}

\usepackage{cleveref}

\graphicspath{{figures/}}

\newtheorem{remark}{Remark}[section]
\crefname{remark}{Remark}{Remarks}
\Crefname{remark}{Remark}{Remarks}

\newcommand{\N}{\ensuremath\mathbb{N}}

\newcommand{\RR}{\ensuremath\mathbb{R}}

\newcommand{\T}{\ensuremath\mathsf{T}}

\newcommand{\integrate}[1]{\,\mathrm{d}#1}

\newcommand{\ds}{\integrate{s}}
\newcommand{\dt}{\integrate{t}}

\DeclareMathOperator{\diag}{diag}

\newcommand{\calD}{\mathcal{D}}

\newcommand{\calH}{\mathcal{H}}

\newcommand{\calM}{\mathcal{M}}
\newcommand{\calN}{\mathcal{N}}

\newcommand{\calP}{\mathcal{P}}
\newcommand{\calQ}{\mathcal{Q}}
\newcommand{\calR}{\mathcal{R}}
\newcommand{\calS}{\mathcal{S}}
\newcommand{\calT}{\mathcal{T}}

\newcommand{\calV}{\mathcal{V}}

\newcommand{\calZ}{\mathcal{Z}}

\newcommand{\ff}{\ensuremath{\bm{f}}}

\newcommand{\Ltwo}{L^2}

\newcommand{\Linfty}{L^\infty}

\definecolor{kit-green}{RGB}{0, 150, 130}
\definecolor{kit-blue}{RGB}{70, 100, 170}
\definecolor{kit-royalblue}{RGB}{0, 45, 76}
\definecolor{kit-iceblue}{RGB}{30, 53, 69}
\definecolor{kit-red}{RGB}{162, 34, 35}
\definecolor{kit-yellow}{RGB}{252, 229, 0}
\definecolor{kit-orange}{RGB}{223, 155, 27}
\definecolor{kit-lightgreen}{RGB}{140, 182, 60}
\definecolor{kit-purple}{RGB}{163, 16, 124}
\definecolor{kit-brown}{RGB}{167, 130, 46}
\definecolor{kit-cyan}{RGB}{35, 161, 224}
\definecolor{kit-gray}{RGB}{0, 0, 0}

\colorlet{optSubspace}{kit-green!80}
\colorlet{optSubspaceDark}{kit-green}
\colorlet{suboptSubspace}{kit-blue!80}
\colorlet{suboptSubspaceDark}{kit-blue}

\newcommand{\PDEspace}{\mathscr{Q}}
\newcommand{\pivotSpace}{\mathscr{P}}

\newcommand{\timeInt}{\calT}
\newcommand{\statePDE}{q}
\newcommand{\statePDEQ}{\statePDE_{\!{}_\calQ}}
\newcommand{\statePDEQz}{\statePDE_{\!{}_\calQ,0}}
\newcommand{\state}{\mathbf{\statePDE}}
\newcommand{\stateDim}{N}
\newcommand{\stateSpace}{\calZ}
\newcommand{\stateSpaceElem}{z}

\newcommand{\reduce}[1]{\hat{#1}}
\newcommand{\statePDEred}{\reduce{\statePDE}}
\newcommand{\stateRed}{\reduce{\state}}
\newcommand{\stateDimRed}{n}

\newcommand{\param}{\boldsymbol{\mu}}
\newcommand{\paramSet}{\calD}
\newcommand{\paramDim}{d'}

\newcommand{\inpVar}{\boldsymbol{u}}
\newcommand{\inpVarDim}{N_{\mathrm{in}}}

\newcommand{\outVar}{\boldsymbol{y}}
\newcommand{\outVarDim}{N_{\mathrm{out}}}
\newcommand{\outVarRed}{\reduce{\outVar}}

\newcommand{\HankelOperator}{\calH}
\newcommand{\ctrlG}{\calP_{\mathrm{c}}}
\newcommand{\obsG}{\calP_{\mathrm{o}}}

\newcommand{\basisFunction}{\phi}
\newcommand{\reducedBasis}{\varphi}

\newcommand{\bfalpha}{\boldsymbol{\alpha}}
\newcommand{\bfx}{\boldsymbol{x}}
\newcommand{\bfh}{\boldsymbol{h}}

\newcommand{\bff}{\boldsymbol{f}}
\newcommand{\bfq}{\boldsymbol{q}}
\newcommand{\bfmu}{\boldsymbol{\mu}}
\newcommand{\bfZ}{\boldsymbol{Z}}
\newcommand{\bfz}{\boldsymbol{z}}
\newcommand{\bfbeta}{\boldsymbol{\beta}}
\newcommand{\bfP}{\boldsymbol{P}}

\newcommand{\bfe}{\boldsymbol{e}}

\newcommand{\bfQ}{\boldsymbol{Q}}
\newcommand{\bfC}{\boldsymbol{C}}
\newcommand{\bfS}{\boldsymbol{S}}
\newcommand{\bfR}{\boldsymbol{R}}
\newcommand{\bfv}{\boldsymbol{v}}
\newcommand{\bfV}{\boldsymbol{V}}
\newcommand{\bfA}{\boldsymbol{A}}
\newcommand{\bfB}{\boldsymbol{B}}

\newcommand{\Ared}{\reduce{\bfA}}
\newcommand{\Bred}{\reduce{\bfB}}
\newcommand{\Cred}{\reduce{\bfC}}

\newcommand{\bftheta}{\boldsymbol{\theta}}

\newcommand{\bfJ}{\boldsymbol{J}}
\newcommand{\bfeta}{\boldsymbol{\eta}}
\newcommand{\bfF}{\boldsymbol{F}}
\newcommand{\bfGamma}{\boldsymbol{\Gamma}}
\newcommand{\bfE}{\boldsymbol{E}}

\newcommand{\Vcal}{\mathcal{V}}
\newcommand{\bfy}{\boldsymbol{y}}
\newcommand{\bfW}{\boldsymbol{W}}
\newcommand{\bfb}{\boldsymbol{b}}

\newcommand{\Dcal}{\mathcal{D}}

\newcommand{\bfU}{\boldsymbol{U}}
\newcommand{\bfSigma}{\boldsymbol{\Sigma}}
\newcommand{\bfr}{\boldsymbol{r}}

\newcommand{\bfpsi}{\boldsymbol{\psi}}

\newcommand{\trafo}{\mathrm{T}}
\newcommand{\shiftOper}{\mathscr{T}}
\newcommand{\pathVar}{\gamma}
\newcommand{\pathVarRed}{\hat{\boldsymbol{\gamma}}}
\newcommand{\coeff}{\hat{q}}
\newcommand{\stateRedSPOD}{\hat{\bfq}}
\newcommand{\massMatrix}{\boldsymbol{M}}
\newcommand{\massMatrixShift}{\boldsymbol{M}_s}

\newcommand{\fred}{\hat{\bff}}

\newcommand{\featuremap}{g}

\newcommand{\woff}{\bftheta_{\text{off}}}
\newcommand{\bfD}{\boldsymbol{D}}
\newcommand{\bfI}{\boldsymbol{I}}

\newcommand{\gdim}{S}
\newcommand{\spts}{m}
\newcommand{\ntmu}{M}
\newcommand{\offdim}{P}

\newcommand{\PhiDECoder}{\Phi}

\newcommand{\timeStep}{\tau} 

\newcommand{\bfVa}{\bfV_1}
\newcommand{\bfVb}{\bfV_2}
\newcommand{\bfG}{\boldsymbol{G}}
\newcommand{\acronym}[1]{#1\xspace}

\newcommand{\POD}{\acronym{POD}}		

\title[Nonlinear model reduction]{Nonlinear model reduction for transport-dominated problems}

\author[J.~S.~Hesthaven, B.~Peherstorfer, B.~Unger]{Jan S.~Hesthaven\\Karlsruhe Institute of Technology \and Benjamin Peherstorfer\\Courant Institute of Mathematical Sciences, New York University
\and Benjamin Unger\\Karlsruhe Institute of Technology}

\begin{document}

\label{firstpage}
\maketitle

\begin{abstract}
This article surveys nonlinear model reduction methods that remain effective in regimes where linear reduced-space approximations are intrinsically inefficient, such as transport-dominated problems with wave-like phenomena and moving coherent structures, which are commonly associated with the Kolmogorov barrier. The article organizes nonlinear model reduction techniques around three key elements---nonlinear parametrizations, reduced dynamics, and online solvers---and categorizes existing approaches into transformation-based methods, online adaptive techniques, and formulations that combine generic nonlinear parametrizations with instantaneous residual minimization.
\end{abstract}

\tableofcontents 

\section{Introduction}
We provide a brief, non-technical overview of model reduction and discuss the challenge of reducing transport-dominated problems. 

\subsection{Model reduction}
In many applications in computational science and engineering, simulations with a numerical model appear within an outer loop, so that multiple numerical simulations are required, at many different parameter configurations, inputs, and initial conditions. Examples include design, control, optimization, uncertainty quantification, and inverse problems. 
The high cost of repeatedly performing numerical simulations is a major barrier to solving outer-loop and other many-query applications.\footnote{See \citet{doi:10.1137/16M1082469} for a discussion about the distinction of outer-loop and many-query applications in the context of multi-fidelity methods. In this survey, we use the terms interchangeably.}  
Model reduction seeks reduced models to perform numerical simulations at greatly reduced costs and has therefore emerged as a key technique for making outer-loop applications tractable; see \citet{Ant05,rozza_reduced_2008,benner2015survey,AntBG20,P22AMS,hesthaven_pagliantini_rozza_2022,annurev:/content/journals/10.1146/annurev-fluid-121021-025220} for surveys. 

Model reduction exploits two fundamental properties of outer-loop applications \cite{rozza_reduced_2008}.
First, across the range of inputs, parameters, and initial conditions relevant to an outer-loop application, the solutions of the numerical simulations typically vary in a structured, rather than arbitrary, manner across the generic solution space. The class of model reduction methods considered in this survey aims to identify this structured variability and leverages it to construct reduced models with reduced solutions that depend on only a small number of degrees of freedom compared to the solutions of the original high-fidelity numerical model; henceforth the full model, and its corresponding full solutions. 
Second, the performance of outer-loop applications is governed primarily by the aggregate computational cost of all numerical simulations. This focus on aggregate costs allows model reduction to incur one-time high pre-processing costs in an offline (or training) phase to construct reduced models, provided that subsequent reduced-model simulations in an online (or test) phase are sufficiently inexpensive for the pre-processing costs to be amortized over the outer loop. 
Leveraging these two opportunities has enabled model reduction to play a key role in making outer-loop applications tractable. 

\begin{remark}
Model reduction is an umbrella term that includes a wide range of approaches beyond the methods emphasized in this survey, including reduced-physics models and purely data-driven surrogate models. In this survey, we focus on so-called projection-based reduced models that are built from full-model information and that represent a reduced solution with a low-dimensional state that evolves under reduced dynamics; even though for nonlinear parametrizations the resulting constructions are not necessarily projections in the classical vector-space sense.  
\end{remark}

\subsection{Linear model reduction and its limitation}\label{sec:Intro:LinMOR}
The classical approach to (projection-based) model reduction is approximating full-model solutions in a reduced space of low dimension, relative to the dimension of the generic full solution space over which the full model is formulated. 
A reduced space is a problem-dependent, low-dimensional approximation space, typically a subspace of the full-solution Hilbert space. It is meant to approximate the set of all full-model solutions over the inputs, parameters, and initial conditions of interest in an outer-loop application; a set that is often referred to as the solution manifold. 
The reduced solutions are linear approximations of the full solutions in the reduced space, which motivates referring to these methods as linear model reduction in the following. We stress that in linear model reduction a single, global reduced space is used. Examples of linear model reduction methods include balanced truncation \cite{1084254,1102568}, proper orthogonal decomposition \cite{sirovich1,holmes},  interpolation methods \cite{Ant05,doi:10.1137/090776925}, and the reduced basis method \cite{Maday2002,rozza_reduced_2008}.

There are linear model reduction methods that are applicable to full models with nonlinear dynamics, meaning numerical models whose governing operators---and thus the time evolution---depend nonlinearly on the solution function. The term linear in linear model reduction refers solely to the ansatz for the reduced solution: the reduced solution is sought in a low-dimensional space and can therefore be represented as a linear combination of basis functions. Consequently, although the approximation ansatz is linear, the resulting reduced model can still exhibit nonlinear dynamics---analogous to finite-element discretizations of nonlinear partial differential equations (PDEs), where the solution field is also represented by a linear combination of basis functions.  
 
Tremendous progress has been made on linear model reduction to turn it into a mature and reliable technology, including effective ways to construct expressive reduced spaces, a posteriori error estimation, handling problems with nonlinear solution field dependence, and structure preservation. We refer to surveys \cite{antoulas-modredsurvey-2001,rozza_reduced_2008,benner2015survey,hesthaven_pagliantini_rozza_2022,annurev:/content/journals/10.1146/annurev-fluid-121021-025220} and monographs \cite{Ant05,QuaR14,HesRS16,QuaMN16,BenCOW17,AntBG20,BenGQRSS21, BenGQRSS21a} that provide overviews of the progress in linear model reduction.

However, linear model reduction has a major, fundamental limitation: Solutions of the full model have not only to lie approximately on a smooth, low-dimensional manifold in the full-model solution space but have to be well approximated by a low-dimensional space with a linear vector space structure. 
The requirement that reduced solutions are linear approximations in a vector space imposes a strong structural restriction. Even if the solution manifold is smooth and of low dimension, it may be strongly curved or otherwise far from linear and therefore cannot be captured well with a low-dimensional vector space. 
For a class of problems related to well-behaved diffusion (elliptic/parabolic) problems, it has been shown that solution manifolds can be well approximated with reduced spaces, i.e., there exist sequences of reduced spaces with increasing dimension $\stateDimRed$ that achieve approximation errors that decay exponentially $\mathcal{O}(\mathrm e^{-\stateDimRed})$ in the reduced dimension $\stateDimRed$ in a suitable  metric \cite{MADAY2002289,buffa_maday_patera_prudhomme_turinici_2012}. For this class of problems, linear model reduction is well suited---and this class has been the primary focus of model reduction for a long time. 
In contrast, transport-dominated problems---where coherent structures such as wave fronts or phase transitions travel through the spatial domain---can lead to solution manifolds that are rough in the sense that linearly approximating these solutions in a reduced space leads to a slow error decay with the reduced dimension \cite{RowM00,OhlR13,P22AMS}. Linear model reduction is ineffective for such transport-dominated problems. 
For example, \citet{OhlR16,GreU19,ArbGU25} derive lower bounds on the error that behave as $\mathcal{O}(1/\sqrt{\stateDimRed})$ when solution manifolds corresponding to instances of the linear advection and wave equation are approximated with linear approximations in reduced spaces of dimension $\stateDimRed$---a substantially slower error decay compared to the exponential decay obtained for some diffusion-dominated problems. This limitation is often referred to as the Kolmogorov barrier. 

Besides these approximation-theoretic error bounds, the limitation of linear model reduction for transport- and wave-like problems has been empirically noted on a wide range of problems of interest in practice, including strongly advection and chemically reacting flows \cite{doi:10.2514/6.2018-1183,HUANG2023112356}, wildfire dynamics \cite{BlaSU21b}, pattern formation \cite{10.1098/rsta.2021.0206}, solar weather predictions \cite{Issan2023}, pulsed detonation combustors \cite{10.1007/978-3-319-98177-2_17}, kinetic plasma problems \cite{NicT21}, and rotating detonation waves \cite{Mendible2020,Uy2024}.

\subsection{Nonlinear model reduction}
Nonlinear model reduction aims to overcome the limited expressivity of linear model reduction by seeking nonlinear approximations of full-model solutions. 
Thus, the reduced solution is allowed to depend nonlinearly on its degrees of freedom, which can provide a more efficient representation of full-model solutions than linear approximations in a reduced space. 
In other words, the trial manifold induced by the nonlinear approximations can represent the solution manifold accurately with fewer degrees of freedom than linear reduced-space approximations. 
From an approximation-theoretic perspective, allowing nonlinearity in the approximations can overcome the fundamental expressivity limitation of linear reduced-space approximations and so mitigate the Kolmogorov barrier. 
However, nonlinear approximations are just one ingredient of nonlinear model reduction.  
In contrast to linear model reduction, the reduced solutions given by the nonlinear approximations are no longer within a linear vector space. In particular, there is in general no global linear projection operator onto the trial manifold and no associated orthogonality principle, so standard Galerkin constructions and many standard stability arguments cannot be carried over verbatim from the full model. 
Nonlinear model reduction therefore requires revisiting the formulation of reduced models, their analysis, and numerical solution methods to account for the nonlinearity of the approximations, which is the topic of this survey.

\begin{remark}
Nonlinear approximations have a long history in numerical analysis, well predating recent nonlinear model reduction developments. For example, \citet{DeVore_1998} investigates the approximation-theoretic aspects of nonlinear approximations over approximations in linear spaces as well as procedures to construct nonlinear approximations. Neural networks are investigated from an approximation-theoretic perspective in \citet{Pinkus_1999,DeVore_Hanin_Petrova_2021,pmlr-v49-telgarsky16}. In finite-element and related numerical methods for PDEs, adaptive mesh refinements lead to approximation spaces that are chosen in a solution-dependent way, so the resulting approximation is a nonlinear function of the underlying degrees of freedom \cite{doi:10.1137/0715049,BERGER198964,BERGER1984484}. Likewise, approximations with rational functions are fundamentally nonlinear in their degrees of freedom, which has been shown to be more effective than polynomial and fixed linear-space approximations for functions with sharp features or near-singular behavior; see, e.g., \citet{BraessBook,TrefethenRatFun}. 
While nonlinear model reduction can build on many of these advances, a key unique feature of nonlinear model reduction is that reduced models are not used to better approximate a single solution function but to build a nonlinear reduced model that can be repeatedly used for a range of parameters, inputs, and initial conditions to accurately approximate full-model solutions. 
\end{remark}

\subsection{Outline of the paper}
After this introduction, we briefly discuss linear model reduction methods in \Cref{sec:LinMOR} to set the stage. The limitation of linear model reduction to reduce transport- and wave-like problems is discussed in \Cref{sec:Kolmogorov}. \Cref{sec:nonlinMOR} provides an overview of nonlinear model reduction and introduces nonlinear parametrizations, reduced dynamics, and online solvers as the three elements of nonlinear model reduction. We then discuss three categories of nonlinear model reduction methods. First, methods based on nonlinear transformations in \Cref{sec:Template}. Second, methods that adapt the reduced space online in \Cref{sec:ADEIM}. Third, methods that build on generic nonlinear parametrizations and employ instantaneous residual minimization in the online phase in \Cref{sec:InstantaneousResMin}. Conclusions and an outlook are given in \Cref{sec:Outlook}.

\section{Linear model reduction}\label{sec:LinMOR}
This section briefly reviews linear model reduction techniques, which seek linear approximations of full-model solutions. Linear model reduction is typically separated into an offline phase, where a reduced model is constructed, and an online phase, where simulations are performed with the reduced model. We discuss the three steps of the offline phase of linear model reduction, which are (i) collecting training data, (ii) training a reduced space, and (iii) constructing a reduced model. We then discuss the lifting bottleneck, which hinders an efficient online phase, and present  empirical interpolation as a way to circumvent the lifting bottleneck and so enable efficient online reduced-model solutions. 
We conclude with a numerical experiment that demonstrates the limited effectiveness of linear model reduction for transport-dominated problems. 

\subsection{Full model}\label{sec:FullModel}
We consider time- and parameter-dependent PDEs of the form 
\begin{equation}
    \label{eqn:PDE}
    \partial_t \statePDEQ(t, \bfx; \param) = f(\bfx, \statePDEQ(t,\cdot;\param); \param)
\end{equation}
where $t\in\timeInt=[0, T]$ is the temporal variable, $\bfx\in\Omega\subseteq\RR^d$ denotes the spatial variable, and $\param\in\paramSet\subseteq\RR^{\paramDim}$ a parameter. In this work, we typically consider situations where the dimension $d$ of the spatial domain is low, $d \in \{1, 2, 3\}$, and the dimension $\paramDim$ of the parameter is at most moderate. Problem \eqref{eqn:PDE} is equipped with an appropriate initial condition $\statePDEQ(0, \cdot; \param)\colon \Omega \to \RR$ for all $\param \in \paramSet$. The solution function $\statePDEQ(t, \cdot; \param)\colon \Omega \to \RR$ is an element of a suitable Hilbert space $\calQ$ that imposes boundary conditions so that the PDE problem is well posed. Note that several extensions to the PDE problem~\eqref{eqn:PDE} are possible: the  problem may also depend on a control input, the right-hand side function $f$ may depend on time $t$, and the solution may be vector-valued, as in systems of PDEs. While model reduction methods exist for these settings, we restrict our attention to the formulation in \eqref{eqn:PDE} to ease exposition.  

Let the solution function~$\statePDEQ(t, \cdot; \bfmu)\colon \Omega \to \RR$ be approximated in a finite-dimensional subspace $\calQ_{\stateDim} \subseteq \calQ$ at all times $t \in \mathcal{T}$ and parameters $\param \in \paramSet$. The space $\calQ_{\stateDim}$ has dimension $\stateDim \in \N$ and is spanned by the basis functions $\phi_1, \dots, \phi_\stateDim\colon \Omega \to \RR$. For example, the finite-dimensional space $\calQ_N$ could be obtained with a finite-elements approach \cite{FEMTheory}. An approximation of $\statePDEQ(t, \cdot; \param)$ in~$\calQ_{\stateDim}$ is a linear combination
\begin{align}
    \label{eqn:PDEapprox:FOM}
    \statePDE(t, \bfx; \param) =  \sum_{i=1}^{\stateDim} \statePDE_i(t, \param) \basisFunction_i(\bfx)\,,
\end{align}
with $N$ coefficients $\statePDE_1(t, \param), \dots, \statePDE_{\stateDim}(t, \param) \in \RR$ that depend on time $t$ and parameter $\param$.  
The ansatz \eqref{eqn:PDEapprox:FOM} can then be inserted into a suitable weak (variational) formulation of the PDE problem \eqref{eqn:PDE} and (Petrov--)Galerkin conditions can be imposed so that the residual vanishes when tested against a test space. It is important to note that, in model reduction, it is typically assumed that the error of the solution~\eqref{eqn:PDEapprox:FOM} in~$\calQ_{\stateDim}$ compared to the corresponding solution $q_{\mathcal{Q}}$ in the space $\mathcal{Q}$ can be controlled with standard numerical analysis tools and is negligible for the purpose of model reduction. Thus, whenever we construct a reduced model via model reduction, it is sufficient to derive a reduced model that is accurate with respect to the solution $q$ given in \eqref{eqn:PDEapprox:FOM} in $\mathcal{Q}_N$. Correspondingly, we refer to the solution $q$ in $\calQ_{\stateDim}$ as the full-model solution. 

Collecting the coefficients of the full-model solution \eqref{eqn:PDEapprox:FOM} into  
\begin{equation}
    \label{eq:PDEFullSolutionSemiDisc}
    \state(t, \param) = \begin{bmatrix} \statePDE_1(t, \param), \dots, \statePDE_{\stateDim}(t, \param)\end{bmatrix}^{\T} \in \RR^{\stateDim},
\end{equation}
gives a vector that describes a function $q(t, \cdot; \bfmu)$ in the space $\mathcal{Q}_N$. Such a coordinate-based representation $\bfq(t, \bfmu)$ of the full-model solution $q(t, \cdot; \bfmu)$ leads to the semi-discrete full model, which is a system of ordinary differential equations (ODEs), 
\begin{equation}
    \label{eqn:FOM}
    \begin{aligned}
        \frac{\mathrm d}{\mathrm dt} \state(t, \param) &= \ff(\state(t, \param); \param), \\
        \state(0; \param) &= \state_0(\param)\,,
    \end{aligned}
\end{equation}
with initial condition $\state_0(\param) \in \mathbb{R}^N$. Correspondingly, the vector $\bfq(t, \bfmu)$ is often referred to as state of the semi-discrete full model. 
There could be a mass matrix resulting from the spatial discretization, which we formally include in the vector field~$\ff$. The mass matrix needs to be appropriately treated in actual implementations. 

\begin{remark} Throughout this survey, we consider parametric full models without explicitly formulating an input--output map. Nevertheless, an important class of model reduction methods originates from systems and control theory, where the explicit definition of inputs and outputs is central. In that setting, the full model is viewed as an input--output dynamical system and reduction can be posed in terms of approximating the induced input--output map. 
Representative approaches include balanced truncation, moment matching and interpolation, and related projection-based techniques, which provide reduced models that approximate the full system's input--output response while often admitting stability and error guarantees, at least in case of linear system dynamics. We refer to \citet{Ant05,benner2015survey,AntBG20} for overviews of system- and control-theoretic model reduction. 
\end{remark}

\subsection{Overview of linear model reduction}
Linear model reduction seeks approximations of the full-model solutions in a low-dimensional subspace $\calV \subseteq \calQ_{\stateDim}$ of dimension $\stateDimRed \ll \stateDim$. 

The reduced space $\calV$ has to be sufficiently expressive to well approximate the full-model solutions over the time interval $\timeInt$ and the parameter domain $\paramSet$. 
Formally, the full-model solutions induce a so-called solution manifold,  
\begin{equation}\label{eq:Prelim:SolManifold}
\mathcal{M} = \{ q(t, \cdot; \param) \,|\, t \in \timeInt, \param \in \paramSet\}\,,
\end{equation}
which is the set of all full-model solutions corresponding to the time and parameter range of interest. We follow the convention of referring to $\mathcal{M}$ as the solution manifold, even though the set $\mathcal{M}$ does not necessarily have a manifold structure; see for instance \citet[Ex.~2.9]{Haasdonk2017}. The solution manifold could be defined as the set of solutions in the space $\mathcal{Q}$ instead of the $\stateDim$-dimensional full-model space~$\mathcal{Q}_{\stateDim}$.

A reduced solution is denoted as a function $\statePDEred\colon \timeInt \times \Omega \times \mathcal{D} \to \RR$. For a given $t$ and $\bfmu$, the function $\hat{q}(t, \cdot; \param)\colon \Omega \to \mathbb{R}$ is an element of $\calV$ and approximates the full-model solution $q(t, \cdot; \param)$. Because $\hat{q}(t, \cdot; \bfmu)$ is in the reduced space $\calV$, it can be represented as a linear combination
\begin{equation}
    \label{eq:Prelim:RedQLinComb}
    \hat{q}(t, \bfx; \param) = \sum_{i = 1}^{\stateDimRed} \hat{q}_i(t, \param) \varphi_i(\bfx)
\end{equation}
of basis functions $\varphi_1, \dots, \varphi_{\stateDimRed}\colon \Omega \to \RR$ of the reduced space $\Vcal$. The coefficients $\hat{q}_1(t, \param), \dots, \hat{q}_\stateDimRed(t, \param) \in \RR$ of the linear combination \eqref{eq:Prelim:RedQLinComb} depend on time $t$ and parameter $\bfmu$. In contrast, the basis functions that span the space $\Vcal$ are fixed over all times $t$ and parameters $\param$. Thus, at all times $t$ and parameters $\param$, the reduced solution~\eqref{eq:Prelim:RedQLinComb} is a linear approximation in $\Vcal$ in the sense that the coefficients (degrees of freedom) $\hat{q}_1(t, \param), \dots \hat{q}_{\stateDimRed}(t, \param)$ enter linearly in $\hat{q}$ when $\hat{q}$ is interpreted as a function of the coefficients. The coefficients of the reduced solution \eqref{eq:Prelim:RedQLinComb} can be obtained via a variational formulation of the PDE problem \eqref{eqn:PDE} and Galerkin or Petrov--Galerkin projection onto $\Vcal$, analogous to how the full-model solution~\eqref{eqn:PDEapprox:FOM} in~$\mathcal{Q}_{\stateDim}$ is obtained from the PDE problem. 
In particular, the reduced solution~$\hat{q}(t, \cdot; \param)$ at all times $t$ and parameters $\param$ is computed in $\Vcal$ so that it depends only on $\stateDimRed$ degrees of freedom, which are the $\stateDimRed$ coefficients in the representation. This is in contrast to the full-model solution \eqref{eqn:PDEapprox:FOM} that depends on $\stateDim$ coefficients with typically $\stateDim \gg \stateDimRed$.

\subsection{The three steps of the offline phase of linear model reduction}\label{sec:LinMOR:OfflineStep}
We identify three main steps of the offline phase of linear model reduction: Step 1 is collecting information about the full model, which is used to construct a suitable reduced space $\Vcal$ that can well approximate the full solutions lying on the solution manifold $\mathcal{M}$. Within the scope of this survey, the information is typically obtained in the form of full-model solutions at selected times and parameters, which constitute the training data. Step 2 is training a reduced space~$\Vcal$ on the training data to well approximate $\mathcal{M}$. The third step is constructing a reduced model that describes a system of equations for the coefficients of reduced solutions \eqref{eq:Prelim:RedQLinComb} in the reduced space~$\Vcal$. 
Once a reduced model has been constructed with these three steps in the offline phase, it then can be used in the online phase to compute reduced solutions within an outer-loop application. 
We will now discuss the three steps of the offline phase in more detail, and then the online phase in \Cref{sec:Prelim:EIM}.

\subsubsection{Offline Step 1: Collecting training data}\label{sec:Prelim:LMOR:Step1}
A standard approach to collecting training data is to solve the full model \eqref{eqn:FOM} for the parameters in a set of training parameters $\paramSet_{\text{train}} = \{\param_1, \dots, \param_{\ntmu}\}\subseteq\paramSet$, which results in a training data set 
\begin{equation}\label{eq:QTrainSet}
    \mathcal{Q}_{\text{train}} = \{\bfq(t, \param_i) \,|\, i = 1, \dots, \ntmu,\, t \in [0, T]\}\,,
\end{equation}
where we use the semi-discrete representation given in \eqref{eq:PDEFullSolutionSemiDisc} of full-model solutions. 

Selecting suitable training parameters $\param_1, \dots, \param_{\ntmu}$ so that the training data set~$\mathcal{Q}_{\text{train}}$ is informative for constructing a reduced space is challenging. A large body of work addresses this problem using greedy procedures that construct the training set incrementally: starting from an initial parameter sample, the reduced model is repeatedly evaluated over a candidate parameter set, an error indicator identifies the parameter where the reduced model performs worst, and this parameter is then added to the training set. This procedure is repeated until the error indicator is below a user-specified tolerance \cite{10.1115/1.1448332,FLD:FLD867,M2AN_2008__42_2_277_0,doi:10.1137/100795772,M2AN_2014__48_1_259_0,RandGreedy}. Such greedy strategies often critically build on advances in a posteriori error estimation \cite{M2AN_2005__39_1_157_0,FLD:FLD867,Haasdonk2011,CRMATH_2012__350_3-4_203_0}. Greedy methods can also explicitly target quantities of interest by constructing reduced spaces so that reduced solutions accurately approximate these quantities of interest; see, e.g.,  \citet{BUITHANH2007880}. 

\subsubsection{Offline Step 2: Training a reduced space}\label{sec:Prelim:LMOR:Step2}
The training data $\mathcal{Q}_{\text{train}}$ from Offline Step 1 can be used to train a low-dimensional space~$\mathcal{V}$ of dimension $n \ll N$, which is typically achieved by constructing basis functions $\varphi_1, \dots, \varphi_n$ that span $\Vcal\subseteq \mathcal{Q}_{\stateDim}$. Consequently, the basis functions $\varphi_1, \dots, \varphi_{\stateDimRed}$ of $\mathcal{V}$ can be represented as linear combinations of the basis functions $\phi_1, \dots, \phi_{\stateDim}$ that span the full-model solution space $\mathcal{Q}_{\stateDim}$,
\begin{equation}
    \label{eq:Prelim:RedBasisAsLinComb}
    \varphi_i = \sum_{j=1}^\stateDim v_i^{(j)}\phi_j\,,\qquad i = 1, \dots, \stateDimRed\,.
\end{equation}
We collect the coefficients into the vectors $\bfv_i = [v_i^{(1)}, \dots, v_i^{(\stateDim)}]^{\top} \in \mathbb{R}^N$ for $i = 1, \dots, \stateDimRed$.  
The coefficient vectors $\bfv_1, \dots, \bfv_{\stateDimRed}$ span an $\stateDimRed$-dimensional subspace of $\RR^{\stateDim}$. 
Thus, we can represent the reduced space by the basis functions $\varphi_1, \dots, \varphi_n$ that span a subspace of $\mathcal{Q}_N$ or, equivalently, by the basis vectors $\bfv_1, \dots, \bfv_n$ that span a subspace of $\mathbb{R}^N$.  

One approach to construct a basis $\bfv_1, \dots, \bfv_n$, and thus implicitly basis functions $\varphi_1, \dots, \varphi_n$ of the form \eqref{eq:Prelim:RedBasisAsLinComb} that span $\Vcal$, is given by the singular value decomposition~(SVD). For this, consider $K$ discrete time points $0 = t_0 < t_1 < \dots < t_K = T$ and the corresponding snapshot matrix of the training data set $\mathcal{Q}_{\text{train}}$,
\begin{equation}
    \label{eq:QtrainSnapshotMatrix}
    \bfQ_{\text{train}} = [\bfq(t_0; \bfmu_1), \dots, \bfq(t_K; \bfmu_1), \bfq(t_0; \bfmu_2), \dots, \bfq(t_K; \bfmu_{\ntmu})] \in \mathbb{R}^{N \times (K + 1)\ntmu}\,.
\end{equation}
Computing the SVD of $\bfQ_{\text{train}}$ and using the first $\stateDimRed$ left-singular vectors corresponding to the $\stateDimRed$ largest singular values gives a basis $\bfv_1, \dots, \bfv_n$ that spans a subspace of~$\RR^{\stateDim}$. 
In model reduction, this process of computing a basis matrix based on the SVD of the snapshot matrix is known as proper orthogonal decomposition (POD); see, e.g., \citet{sirovich1,holmes} for early usage in the context of fluid mechanics and \citet{GubV17} for a numerical analysis point of view. There is a wide range of other techniques than POD to construct reduced spaces, including greedy methods \cite{10.1115/1.1448332,FLD:FLD867,BUITHANH2007880,doi:10.1137/100795772,RandGreedy}, interpolation-based methods \cite{Ant05,doi:10.1137/090776925,AntBG20}, and balancing \cite{1084254,1102568}. 

The outcome of Offline Step 2 is a basis of an $\stateDimRed$-dimensional reduced space~$\Vcal$, either represented in terms of the basis functions $\varphi_1, \dots, \varphi_{\stateDimRed}\in \mathcal{Q}_{\stateDim}$ or the coefficients vectors $\bfv_1, \dots, \bfv_{\stateDimRed}\in\RR^{\stateDim}$. 

\subsubsection{Offline step 3: Constructing a reduced model}\label{sec:Prelim:LMOR:Step3}
The third offline step of linear model reduction is constructing a reduced model, 
which means deriving a system of equations that determines the coefficients $\hat{q}_1(t, \bfmu), \dots, \hat{q}_n(t, \bfmu)$ of a reduced solution $\hat{q}$ as represented in \eqref{eq:Prelim:RedQLinComb}. 

Analogous to how the full model was derived, the reduced model can be obtained with a variational formulation and a Galerkin or Petrov--Galerkin projection, but now onto the reduced space $\Vcal$ rather than the full-model space $\mathcal{Q}_N$. 
In the semi-discrete setting, the reduced model becomes a system of ODEs: The representation \eqref{eq:Prelim:RedBasisAsLinComb} allows us to identify a function $\hat{q}(t, \cdot; \bfmu) \in \Vcal$ as given in  \eqref{eq:Prelim:RedQLinComb} with the coefficient vector
\begin{equation}\label{eq:Prelim:RedCoeffVector}
\hat{\bfq}(t, \bfmu) = [\hat{q}_1(t, \bfmu), \dots, \hat{q}_n(t, \bfmu)]^{\top} \in \mathbb{R}^n\,,
\end{equation}
which is an element of the $n$-dimensional subspace of $\mathbb{R}^n$ that is spanned by the columns of the basis matrix $\bfV = [\bfv_1, \dots, \bfv_n] \in \mathbb{R}^{N \times n}$. We note that in many cases, the reduced basis vectors $\bfv_1, \dots, \bfv_n$ (and correspondingly the reduced basis functions $\varphi_1, \dots, \varphi_n)$ are constructed so that they are orthonormal. In particular, in these cases, the basis matrix $\bfV$ has orthonormal columns. If not otherwise noted, we assume that $\bfV$ has orthonormal columns in the following.  

Building on the basis matrix $\bfV$ and the semi-discrete full model \eqref{eqn:FOM},  the semi-discrete reduced model is
\begin{equation}
    \label{eqn:ROM}
    \frac{\mathrm d}{\mathrm dt} \hat{\bfq}(t, \bfmu) = \hat{\bff}(\hat{\bfq}(t, \bfmu); \bfmu)\,,
\end{equation}
where the right-hand side function $\hat{\bff}$ is obtained via Galerkin projection as
\begin{equation}\label{eq:Prelim:HatFGalerkin}
\hat{\bff}(\hat{\bfq}(t, \bfmu); \bfmu) = \bfV^{\top}\bff(\bfV\hat{\bfq}(t, \bfmu); \bfmu)\,.
\end{equation}
Analogously to the state $\bfq(t, \bfmu)$ of the semi-discrete full model, the vector $\hat{\bfq}(t, \bfmu)$ is often referred to as the reduced state. The initial condition $\hat{\bfq}(0; \bfmu) \in \RR^{\stateDimRed}$ can be obtained, for example, via the projection $\hat{\bfq}(0; \bfmu) = \bfV^{\top}\bfq(0; \bfmu)$. 

While we define the right-hand side function via Galerkin projection in \eqref{eq:Prelim:HatFGalerkin}, we note that Petrov--Galerkin projections, with test and trial spaces that are different, is essential for stability when the right-hand side function $\bff$ stems from non-symmetric PDE operators \cite{QuaR14,HesRS16,MehU23}. Model reduction methods based on interpolation also typically derive test and trial spaces that are different from each other so that the reduced model is obtained with Petrov--Galerkin projection \cite{Ant05,doi:10.1137/090776925,AntBG20}. Using oblique projections in a Petrov--Galerkin setting can also improve the predictive capabilities of reduced models as demonstrated by \citet{OttPR22,Otto2023}. 
Improving robustness and stability is also one of the motivations of least-squares Petrov--Galerkin methods that were introduced by \citet{carlberg2011efficient} and \citet{CARLBERG2017693}. The works \citet{dahmen_plesken_welper_2014,doi:10.1137/18M1176269,doi:10.1137/23M1613402} explicitly target transport-dominated problems and develop reduced models that remain stable for these problems; however, their emphasis is on robustness and stability of linear model reduction rather than overcoming the Kolmogorov barrier of linear approximations. A different approach is to combine linear model reduction formulations with stabilizing terms (e.g., \citet{PACCIARINI20141}) or to develop online stabilization techniques (e.g., \citet{CRMATH_2016__354_12_1188_0,BALAJEWICZ2016224}).  

It is important to notice that the right-hand side function~$\hat{\bff}$ of the reduced model~\eqref{eqn:ROM} is not necessarily linear in the reduced state $\hat{\bfq}(t, \bfmu)$, even though $\hat{\bfq}(t, \bfmu)$ is a linear approximation of the full-model state $\bfq(t, \bfmu)$ in the $n$-dimensional subspace spanned by the columns of $\bfV$. This stresses once more that linear model reduction methods are applicable to nonlinear full models and the term linear in linear model reduction refers solely to the linear ansatz for the reduced solution. 

\begin{remark}
One generally distinguishes between intrusive and non-intrusive model reduction. Non-intrusive methods learn both the reduced representation such as the reduced space as well as the reduced dynamics primarily from snapshot data; see, e.g.,  \citet{Ghattas_Willcox_2021,annurev:/content/journals/10.1146/annurev-fluid-121021-025220}. In contrast, in this survey, we focus on intrusive methods that derive the reduced dynamics from the governing equations via, e.g., projection or more general residual minimization techniques. We briefly revisit non-intrusive model reduction in \Cref{sec:Outlook}. 
\end{remark}

\subsection{Online computations, the lifting bottleneck, and empirical interpolation}\label{sec:Prelim:EIM}
We now discuss the online phase of linear model reduction, which concerns performing numerical simulations efficiently with the reduced model for a given new parameter and initial condition. In particular, we discuss the lifting bottleneck, which describes the issue that applying a projection directly to derive the reduced right-hand side function can still require online computations that scale unfavorably with the full space dimension $\stateDim$, which can result in little to no speedup. We then briefly review empirical interpolation as a strategy to overcome the lifting bottleneck. 

\subsubsection{Online efficiency} In the online phase, the reduced model is used to perform simulations---also referred to as model evaluations---for new parameters $\bfmu$ and initial conditions by integrating the reduced dynamical system \eqref{eq:Prelim:HatFGalerkin} in time. This process requires repeated evaluations of the reduced right-hand side function $\hat{\bff}$ and, for implicit time-integration schemes, typically also evaluations of derivatives of $\hat{\bff}$ within, e.g.,  Newton iterations. Consequently, the computational costs of evaluating $\hat{\bff}$ are typically dominating the online costs. 

We speak of online-efficient reduced models if the computational costs of the online phase scale independently of the dimension of the full model $\stateDim$ \cite{rozza_reduced_2008}. A classical setting in which online efficiency can be often achieved is when the full-model right-hand side depends linearly on the state and affinely on the parameter, i.e.,
\begin{align*}
    \bff(\bfq;\param) = \sum_{i=1}^\ell \rho_i(\param)\boldsymbol{A}_i\bfq,\qquad \rho_i\colon\paramSet\to\RR, \quad\boldsymbol{A}_i\in\RR^{\stateDim\times\stateDim}.
\end{align*}
In this case, the parameter-independent reduced operators $\bfV^{\top}\boldsymbol{A}_i\bfV$ can be pre-computed offline, while the assembly and evaluation of $\hat{\bff}$ for a given parameter incur only costs independent of the full-model dimension $\stateDim$ \cite{rozza_reduced_2008,QuaR14,HesRS16}. If one needs to evaluate the reduced model only at a certain time instant or a small time window, then online efficiency can be further improved using contour integral methods; see for instance work by \citet{GugLM21} and \citet{GugM23}. Outside such special structures, however, evaluating the reduced right-hand side $\hat{\bff}$ generally still involves operations in the high-dimensional full-model space, which fundamentally limits online efficiency.

\subsubsection{Lifting bottleneck}\label{sec:Prelim:LiftingBottleneckSub} To see why online efficiency is difficult to achieve in general, recall that the reduced right-hand side $\hat{\bff}$ in~\eqref{eq:Prelim:HatFGalerkin} is obtained from the full-model right-hand side~$\bff$ via (Petrov--)Galerkin projection. 
Evaluating $\hat{\bff}$ at a reduced state $\hat{\bfq}(t, \param) \in \RR^{\stateDimRed}$ at time $t$ and parameter $\param$ therefore typically proceeds as follows: the reduced state is first lifted to the high-dimensional space as $\bfV\hat{\bfq}(t, \param) \in \RR^{\stateDim}$, then the full-model right-hand side function $\bff$ is evaluated at this lifted state $\bfV\hat{\bfq}(t, \param)$, which requires evaluating all $N$ component functions, and finally, the output of $\bff$ is projected onto the reduced space spanned by the columns of $\bfV$.

In special cases, such as the affine and linear setting discussed above,  the terms involved in the evaluation of $\hat{\bff}$ can be pre-computed in the offline phase and then reused online \cite{rozza_reduced_2008}. However, in more general situations, the lifting of the reduced state, the evaluation of $\bff$ in $\RR^N$, and the subsequent projection imply that computing the reduced right-hand side $\hat{\bff}$ incurs computational costs that scale with the full-model dimension $N$. The need to repeatedly lift quantities to the high-dimensional space and perform computations there during online evaluation is commonly referred to as the lifting bottleneck. 

\subsubsection{Empirical interpolation}\label{sec:EmpiricalRegression} The empirical interpolation method has been introduced by \citet{Barrault2004} to circumvent the lifting bottleneck. Empirical interpolation replaces the projection in \eqref{eq:Prelim:HatFGalerkin} with interpolation so that only a few of the $N$ component functions of the full-model right-hand side function need to be evaluated. 
Several other methods have been proposed to circumvent the lifting bottleneck with similar interpolation and regression strategies, such as missing point estimation \cite{AstWWB04,4668528}, empirical operator learning \cite{MR2537223,MR2914310}, Gauss--Newton approximated tensors \cite{carlberg2011efficient}, and the energy-conserving sampling and weighting method that focus on preserving structure and stability of underlying full-model discretizations \cite{https://doi.org/10.1002/nme.4668,https://doi.org/10.1002/nme.4820}. 

\paragraph{Steps of empirical interpolation} We focus here on the empirical interpolation method by \citet{Barrault2004}, and in particular on its discrete counterpart, the discrete empirical interpolation method introduced by \citet{Chaturantabut2010}. In the offline phase, pairwise distinct interpolation points are selected, 
\begin{equation}\label{eq:EIM:InterpolationPoints}
p_1, \dots, p_n \in \{1, \dots, N\}\,,
\end{equation}
which correspond in the discrete setting to indices of the $N$ component functions of the full-model right-hand side function $\bff$. 
The interpolation points give rise to the interpolation points matrix
\[
\bfP = \begin{bmatrix}\bfe_{p_1} & \dots & \bfe_{p_n}\end{bmatrix}^{\top} \in \mathbb{R}^{N \times n}
\]
so that when applying a vector $\bfz = [z_1, \dots, z_N]^{\top} \in \mathbb{R}^N$ to $\bfP^{\top}$, the components corresponding to indices $p_1, \dots, p_n$ are selected:
\[
\bfP^{\top}\bfz = \begin{bmatrix} z_{p_1}&  \dots & z_{p_n}\end{bmatrix}^{\top}\,.
\]
In implementations, the interpolation points matrix $\bfP$ is not assembled and instead its application to a vector is mimicked via slicing operations \cite{Chaturantabut2010}. 
The reduced right-hand side function is then obtained by interpolating the full-model right-hand side $\bff$ at the interpolation points $p_1, \dots, p_n$ in the basis $\bfv_1, \dots, \bfv_n$, i.e.,
\begin{equation}\label{eq:Prelim:EIMHatF}
\tilde{\bff}(\hat{\bfq}; \bfmu) = (\bfP^{\top}\bfV)^{-1}\bfP^{\top}\bff(\bfV\hat{\bfq}; \bfmu)\,.
\end{equation}
A reduced model based on empirical interpolation is analogous to the model given in \eqref{eqn:ROM} except that the reduced right-hand side $\tilde{\bff}$ given in \eqref{eq:Prelim:EIMHatF} is used, rather than the right-hand side $\hat{\bff}$ defined in \eqref{eq:Prelim:HatFGalerkin} via projection. We remark that empirical interpolation can also be applied to individual terms of the right-hand side function~$\bff$. For example, \citet{Barrault2004,Chaturantabut2010} apply empirical interpolation only to nonlinear terms of $\bff$ while projecting the linear terms. Moreover, one may employ a separate reduced space for the approximation of the nonlinear terms, i.e., use distinct reduced spaces for the state and for the empirical-interpolation approximation of the nonlinear terms of $\bff$. 

In the online phase, using $\tilde{\bff}$ given via empirical interpolation \eqref{eq:Prelim:EIMHatF} as right-hand side function in the reduced model can lead to runtime speedups even when using~$\hat{\bff}$ obtained via the projection \eqref{eq:Prelim:HatFGalerkin} does not. If the full-model right-hand side function~$\bff$ satisfies some structural properties, e.g., each component function of $\bff$ can be evaluated independently of all other component functions, then computing the term $\bfP^{\top}\bff(\bfV\hat{\bfq}(t, \bfmu); \bfmu)$ of \eqref{eq:Prelim:EIMHatF} can require evaluating only $n$ instead of all $N$ component functions of $\bff$. Thus, in this case, empirical interpolation circumvents the lifting bottleneck. We refer to discussions by \citet{Barrault2004} and \citet{Chaturantabut2010} for more details.

\paragraph{Selection of interpolation points} 
The selection of the interpolation points  \eqref{eq:EIM:InterpolationPoints} is critical for the accuracy of empirical-interpolation approximations. A range of methods has been developed for selecting the interpolation points. \citet{Barrault2004} and \citet{Chaturantabut2010} propose greedy methods that sequentially add points based on the interpolation residual in each step. The QDEIM algorithm introduced by \citet{drmac-gugercin-DEIM-2016} selects interpolation points based on the pivoting elements of the QR factorization of the transposed basis matrix $\bfV^{\top}$, which can lead to sharper error bounds than the ones corresponding to the greedy selection of the interpolation points; this QR-based sampling approach is related to max-volume submatrix selection \cite{drmac-gugercin-DEIM-2016}.

\paragraph{Oversampling and empirical regression} 
While standard empirical interpolation enforces interpolation at exactly $n$ sampling points, a common generalization is to relax this requirement by using more sampling points than the reduced-space dimension. Selecting $n_s > n$ points $p_1, \dots, p_n, p_{n + 1}, \dots, p_{n_s}$ leads to empirical regression,
\[
\tilde{\bff}(\hat{\bfq}; \bfmu) = (\bfP^{\top}\bfV)^{+}\bfP^{\top}\bff(\bfV\hat{\bfq}; \bfmu)\,,
\]
where $(\bfP^{\top}\bfV)^{+}$ denotes the Moore--Penrose pseudo inverse of $\bfP^{\top}\bfV$. Taking $n_s > n$ points is referred to as oversampling and has been used in the context of model reduction first by methods motivated by the gappy POD method \cite{Everson1995} such as  \cite{AstWWB04,4668528,carlberg2011efficient}. One situation where oversampling is critical has been identified by \citet{ARGAUD2018354}, who show that adding noise or perturbations can lead to instabilities in empirical interpolation when using $n = n_s$ points. \citet{doi:10.1137/19M1307391} propose to randomly oversample, i.e., $n_s > n$, and prove that a randomized selection of points can mitigated this instability. 

It is critical to note that neither the greedy algorithm introduced by \citet{Barrault2004,Chaturantabut2010} nor the QDEIM algorithm by \citet{drmac-gugercin-DEIM-2016} can be directly applied to selecting $n_s > n$ points. Instead, a range of oversampling strategies are available \cite{4668528,carlberg2011efficient}. One line of work focuses on  minimizing the norm of the sampling operator $\bfP^{\top}\bfV$ and interprets adding a point as a low-rank update to the operator $\bfP^{\top}\bfV$ to derive selection criteria \cite{zimm-will-2016,doi:10.1137/19M1307391}.

\subsection{Linear model reduction for diffusion- and transport-dominated problems}\label{sec:Intro:LinNotEnough}
Let us now consider a numerical example to demonstrate when linear model reduction becomes inefficient for transport-dominated problems. The following presentation follows \citet{P22AMS}. We distinguish between diffusion- and transport-dominated problems. 

\begin{figure}[t]
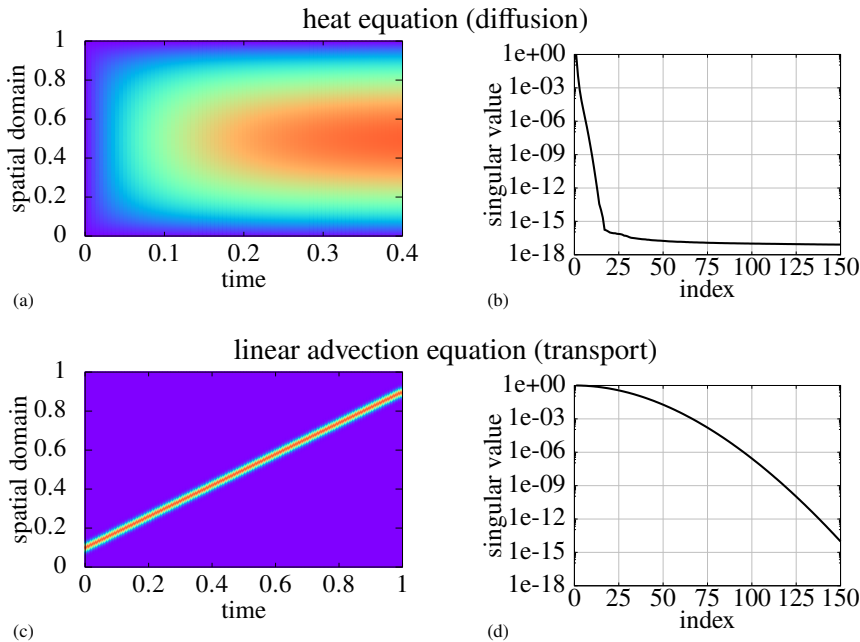

\centering
\begin{tabular}{ll}
\multicolumn{2}{c}{\text{heat equation (diffusion)}}\vspace*{-0.6cm}\\
{\huge\resizebox{0.46\columnwidth}{!}{\input{figures/HeatSolutionTimeSpace}}} &
{\Huge\resizebox{0.39\columnwidth}{!}{\input{figures/HeatSingVal}}}\vspace*{-0.4cm}\\
\hspace*{-0.15cm}\scriptsize (a) & \hspace*{-0.1cm}\scriptsize (b)\vspace*{0.25cm}\\
\multicolumn{2}{c}{\text{linear advection equation (transport)}}\vspace*{-0.6cm}\\
{\huge\resizebox{0.46\columnwidth}{!}{\input{figures/AdvecSolutionTimeSpace}}} &
{\Huge\resizebox{0.39\columnwidth}{!}{\input{figures/AdVecSingValGlobal}}}\vspace*{-0.4cm}\\
\hspace*{-0.15cm}\scriptsize (c) & \hspace*{-0.1cm}\scriptsize (d)
\end{tabular}
\caption{The snapshot matrix  corresponding to the heat-equation example has rapidly decaying singular values, indicating that a low-dimensional subspace can capture the snapshot set efficiently. By contrast, the advection equation leads to a substantially slower decay, highlighting the limitations of linear reduced-space approximations for transport-dominated dynamics. (First published in Notices of the American Mathematical Society in 69, Number 5 (2022), published by American Mathematical Society. \textcopyright ~2022 American Mathematical Society.)}
\label{fig:Need:NumExp}
\end{figure}

\subsubsection{Diffusion problem: heat equation}
As an example of a diffusion-dominated problem, consider the heat equation
\begin{equation}\label{eq:HeatEquation}
\partial_t q_{\mathcal{Q}}(t, x; \mu) - \mu\partial_x^2 q_{\mathcal{Q}}(t, x; \mu) = 1\,,\qquad x \in \Omega\,,
\end{equation}
over the one-dimensional spatial domain $\Omega = [0, 1]$ with homogeneous Dirichlet boundary conditions. The initial condition is the constant function that evaluates to zero in $\Omega$. The parameter $\mu$ corresponds to the heat-conductivity coefficient and is set to $\mu = 1$ in this example. 
Figure~\ref{fig:Need:NumExp}a shows a numerical solution of the heat equation \eqref{eq:HeatEquation} over time and space. The numerical solution is computed with linear finite elements on $N = 1024$ nodes (including boundary points) in $\Omega$ and implicit Euler in time with $K = 400$ time steps of time-step size $10^{-3}$ over the time interval $\mathcal{T} = [0, 0.4]$. Let us now collect snapshots \eqref{eq:QTrainSet} as in Offline Step 1 of linear model reduction outlined in \Cref{sec:Prelim:LMOR:Step1}. Recall that $\mu = 1$ is fixed. We take snapshots at the $K + 1$ discrete time points to form the snapshot matrix $\bfQ_{\text{train}}$ as defined in \eqref{eq:QtrainSnapshotMatrix} and plot the normalized singular values in \Cref{fig:Need:NumExp}b. The normalized singular values are the ratios $\bar{\sigma}_i = \sigma_i/\sigma_1$ for $i = 1, \dots, K+1$ so that the first normalized singular value is one. 

The normalized singular values plotted in \Cref{fig:Need:NumExp}b decay rapidly in the sense that the space spanned by the first $n = 20$ left-singular vectors approximates the snapshot data up to double precision in the relative Frobenius norm. The dimension of the snapshot vectors obtained from the full model is $N = 1024$ in this example (including the two boundary points), which is more than $50$ times higher than the reduced dimension $n$. The fast decay of the singular values shows that the snapshots collected in this heat equation example can be well approximated with a low-dimensional space $\Vcal$. While the singular values provide no guarantee about the approximation error of full solutions other than the ones in the snapshot matrix, the singular values are still often used as a heuristic for how much reduction can be expected.

\subsubsection{Transport problem: linear advection equation}
Let us now consider a canonical example of a transport problem, namely the linear advection equation,
\begin{equation}\label{eq:LinAdvProblem}
\partial_t q_{\mathcal{Q}}(t, x; \mu) + \mu \partial_x q_{\mathcal{Q}}(t, x; \mu) = 0\,,
\end{equation}
over the real line, where $\mu$ corresponds to the advection speed and is fixed to $\mu = 4/5$ in this example. 
The initial condition $q_{\mathcal{Q},0}: \mathbb{R} \to \mathbb{R}$ is a Gaussian probability density function with mean $10^{-1}$ and standard deviation $1.5 \times 10^{-2}$. 
The closed-form solution of the problem \eqref{eq:LinAdvProblem} is given by $(t, x) \mapsto q_{\mathcal{Q},0}(x - 4t/5)$ because we fixed $\mu = 4/5$. We evaluate the analytical solution at the $N = 1024$ equidistant points in $\Omega = [0, 1]$ and take $K = 400$ equidistant time steps in $[0, 1]$, leading to $K + 1$ snapshots in time. The snapshots are plotted in the time-space representation in \Cref{fig:Need:NumExp}c. The normalized singular values of the corresponding snapshot matrix are plotted in \Cref{fig:Need:NumExp}d. 

The normalized singular values decay is slower than in the case of the diffusion-dominated problem. In particular, the results appear to be consistent with an algebraic decay, in contrast to the rapid decay observed for the diffusion case. 
A space of dimension $n\geq 150$ is needed to represent the snapshot data up to double precision in the relative Frobenius norm, which yields a compression ratio $ N/n$ about one order of magnitude lower than for the diffusion problem. 
The slow decay is consistent with the geometric structure of the solution manifold for this transport problem: the dynamics only translate a localized feature (the Gaussian bump) through the domain. This means that the snapshots are shifted copies of the same function, and a vector space must use many basis functions to represent the same feature at different spatial locations.
Thus, the issue is not that the solution manifold is high-dimensional---in fact, in this example, it is one-dimensional because the parameter $\mu$ is fixed---but rather that it is inefficient to approximate this solution manifold by a low-dimensional vector space, as in linear model reduction. 

In summary, this numerical illustration indicates that the efficiency of linear model reduction, as measured by the dimension of the reduced space required to reach a given accuracy, can be much lower for transport- than for diffusion-dominated problems.

\section{Kolmogorov barrier}
\label{sec:Kolmogorov}
A central tool for understanding the performance limits of linear model reduction is the Kolmogorov $\stateDimRed$-width, which measures the best possible worst-case error achievable by approximations restricted to $\stateDimRed$-dimensional subspaces. In this section, we briefly review this concept and discuss lower bounds on the $\stateDimRed$-width, which show the fundamental limitations of linear model reduction methods for transport-dominated problems and motivate the notion of a  Kolmogorov barrier. This section is intended to provide only a glimpse at the relevant approximation-theoretic concepts. For comprehensive treatments of the approximation-theoretic foundations of model reduction, we refer to  \citet{CohD15,10.1093/imanum/dru066,BenCOW17,Cohen2022}.

\subsection{Best linear approximations in vector spaces: the Kolmogorov \texorpdfstring{$\stateDimRed$}{n}-width}

Let $\mathcal{Z}$ be a Hilbert space and let $\calN \subseteq \mathcal{Z}$ be a subset. We might think of $\mathcal{Z}$ as the space over which the PDE problem \eqref{eqn:PDE} is defined or the finite-dimensional approximation space corresponding to the full model \eqref{eqn:FOM}. Correspondingly, the subset $\calN$ can be informally thought of as the set of all solutions to the PDE problem~\eqref{eqn:PDE} or~the full model \eqref{eqn:FOM}, respectively. For a given $\stateDimRed$-dimensional subspace $\stateSpace_{\stateDimRed}\subseteq\mathcal{Z}$, we define the deviation of $\calN$ from~$\stateSpace_{\stateDimRed}$ as
\begin{equation}\label{eq:K:DFirst}
    d(\calN,\stateSpace_{\stateDimRed}) = \sup_{\stateSpaceElem\in\calN} \inf_{\stateSpaceElem_{\stateDimRed}\in\stateSpace_{\stateDimRed}} \|\stateSpaceElem - \stateSpaceElem_{\stateDimRed}\|,
\end{equation}
which can be interpreted as the worst-case situation for the best approximation. Because $\mathcal{Z}$ is a Hilbert space, the infimum is obtained by the orthogonal projection of $\stateSpaceElem$ onto $\stateSpace_{\stateDimRed}$, which always exists as $\stateSpace_{\stateDimRed}$ is of finite dimension $n$. Furthermore, the projection map is linear, which means restricting best approximations to lie in a subspace of a Hilbert space lead to linear (subspace-based) approximations. 

Let us now take the infimum of \eqref{eq:K:DFirst} over all $\stateDimRed$-dimensional subspaces of $\mathcal{Z}$ and so define the Kolmogorov $n$-width
\begin{align}
    \label{eqn:KolmogorovNWidth}
    d_n(\calN,\stateSpace) = \inf_{\substack{\stateSpace_{\stateDimRed}\leq \stateSpace\\\dim(\stateSpace_{\stateDimRed})\leq \stateDimRed}} d(\calN,\stateSpace_{\stateDimRed}) = \inf_{\substack{\stateSpace_{\stateDimRed}\leq \stateSpace\\\dim(\stateSpace_{\stateDimRed})\leq \stateDimRed}} \sup_{\stateSpaceElem\in\calN} \inf_{\stateSpaceElem_{\stateDimRed}\in\stateSpace_{\stateDimRed}} \|\stateSpaceElem - \stateSpaceElem_{\stateDimRed}\|,
\end{align}
an idea first proposed by \citet{Kol36}. Following \citet{Pin85}, an $\stateDimRed$-dimensional subspace $\stateSpace_{\stateDimRed}^\star$ of $\stateSpace$ is called an optimal subspace for $d_n(\calN,\stateSpace)$ if 
\begin{align*}
    d_n(\calN,\stateSpace) = d(\calN,\stateSpace_{\stateDimRed}^\star).
\end{align*}

\subsection{Fast and slow decay of Kolmogorov $n$-width and Kolmogorov barrier}\label{sec:KolmogorovNonLinExamples}
Let us now apply the concept of the Kolmogorov $n$-width to the setting of model reduction. As mentioned above, we might think of $\mathcal{Z}$ as the space $\mathcal{Q}$ over which the PDE problem \eqref{eqn:PDE} is defined and of $\calN$ as the set of all solutions. Analogously, we can think of $\mathcal{Z}$ as the finite-dimensional approximation space $\mathcal{Q}_N$ corresponding to the full model \eqref{eqn:FOM} and of the subset $\calN$ as the set of all solutions \eqref{eq:Prelim:SolManifold} to the full model. 

Under this setting, the Kolmogorov $n$-width \eqref{eqn:KolmogorovNWidth} is the smallest worst-case error that can be achieved when approximating the solutions in an $n$-dimensional reduced space. Consequently, how fast or slow the Kolmogorov $n$-width decays with the dimension $n$ of the reduced space has major implications on the efficiency of linear model reduction. If it decays quickly, then a reduced low-dimensional space is often sufficient to meet a given error tolerance. Thus, in such a situation, linear model reduction can be efficient, at least from an approximation-theoretic perspective. In contrast, if the Kolmogorov $n$-width decays slowly, then the dimension of the reduced space to meet the error tolerance might be so high that linear model reduction becomes inefficient. In this case, we speak of the Kolmogorov barrier, because it presents a hard limitation of linear model reduction

Let us consider examples from the model reduction literature to clarify what fast and slow decay means in this context: \citet{MADAY2002289} consider  elliptic coercive PDEs over a one-dimensional parameter $\mu$ and the corresponding set of solutions. \citet{MADAY2002289} show that in this situation there is an upper bound on the Kolmogorov $n$-width that decays exponentially fast in the dimension $n$ of the reduced space. Thus, increasing the dimension of the reduced space in linear model reduction can achieve an exponential error reduction. 
Analog results of a rapid, albeit not necessarily exponential, decay for more general elliptic problems are shown in  \citet{doi:10.1142/S0219530511001728,CohD15,10.1093/imanum/dru066,Bachmayr2017Kolmogorov}. 
As a remark, we add that beyond existence results,  \citet{doi:10.1137/100795772,buffa_maday_patera_prudhomme_turinici_2012,RandGreedy} show that sequences of reduced spaces can be constructed with greedy methods so that the same decay up to constants as that of the $n$-width is attained.

Broadly speaking, often diffusion-dominated problems, such as the one discussed in \Cref{sec:Intro:LinNotEnough}, fall in this category for which the Kolmogorov $n$-width decays quickly and thus for which linear model reduction is well suited.

\subsection{Slow $n$-width decay for linear advection problem}
\label{subsec:Kolmogorov:lowerBound}
Let us now consider examples of PDE problems for which lower bounds on the Kolmogorov $n$-width of solution manifolds have been shown.

Let us consider the linear advection problem over the one-dimensional spatial domain $\Omega = [0, 1]$ and time domain $[0, 1]$,
\begin{equation}\label{eq:LinAdvProblemEq}
\partial_t q_{\mathcal{Q}}(t, x; \mu) + \mu\partial_x q_{\mathcal{Q}}(t, x; \mu) = 0\,,
\end{equation}
with the conditions
\[
q_{\mathcal{Q}}(0, x; \mu) = 0\,,\qquad q_{\mathcal{Q}}(t, 0; \mu) = 1\,,
\]
for all $\mu \in [0, 1]$. 
Let us consider the corresponding solution manifold for just $\mu = 1$,
\begin{equation}\label{eq:LinAdvNWidthProblem:SolManifoldMCal}
\mathcal{M} = \{q_{\mathcal{Q}}(t, \cdot; \mu)\, | \, t \in [0, 1], \mu = 1\}\,,
\end{equation}
which consists of the step functions 
\begin{equation}\label{eq:NWidthExample:StepFunSol}
q_{\mathcal{Q}}(t, x; \mu) = \begin{cases} 1\,, &0 \leq x \leq t\mu\,,\\
0\,, &t\mu < x \leq 1\,,\end{cases}
\end{equation}
that have a discontinuity somewhere in the spatial domain $[0, 1]$. Notice that we use the set of solutions $\mathcal{M}$ not over a finite-dimensional (full-model) space $\mathcal{Q}_N$ here as in \eqref{eq:Prelim:SolManifold} but over  $\mathcal{Q} = L^2([0, 1])$.
\citet{OhlR16} show that the Kolmogorov $n$-width with respect to the $L^2$ norm is lower bounded as
\[
d_n(\mathcal{M}, \mathcal{Q}) \geq c \frac{1}{\sqrt{n}}\,,
\]
where $n$ is the dimension of the reduced space and $c > 0$ is a constant. 

We can interpret the bound as showing that there cannot exist a sequence of subspaces of $L^2$ that achieve a faster approximation error in the sense of the Kolmogorov $n$-width \eqref{eqn:KolmogorovNWidth} as $1/\sqrt{n}$---which is a substantially slower error decay than the exponential error decay that can be achieved for the one-parameter diffusion-dominated problem discussed in \Cref{sec:KolmogorovNonLinExamples}. 
This result is also in agreement with the numerical experiments shown in \Cref{sec:Intro:LinNotEnough}, where the example based on the linear advection equation empirically exhibits a slower decay  of the projection error given by the singular values, whereas the diffusion-dominated problem given by the instance of the heat equation leads to an exponential decay of the singular values. We stress once more that the decay of the singular values does not directly correspond to the decay of the Kolmogorov $n$-width; however, the singular values can provide an empirical indication under certain conditions; see \Cref{sec:Intro:LinNotEnough}. 

\citet{ArbGU25} show that if the initial condition in the linear advection example is chosen smoother, then this smoothness is also reflected in the decay rate of the Kolmogorov $n$-width. Analogous results for the wave equation are shown in \citet{GreU19}.

Overall, these approximation-theoretic results, as well as the empirical evidence discussed in \Cref{sec:Intro:LinNotEnough} and the cited literature in \Cref{sec:Intro:LinMOR}, suggest that transport-dominated problems with wave-like phenomena and moving coherent, localized features lead to slowly decaying $n$-widths, which fundamentally limit linear model reduction techniques. 

\begin{remark} In this survey, we focus on time-dependent problems; however, also for stationary problems, linear model reduction can be come inefficient due to the Kolmogorov barrier. For example, \citet{doi:10.1137/16M1059904,Tad20,ALIREZAMIRHOSEINI2023111739} consider problems in which the parameter $\param$ induces the transport of a coherent structure in the spatial domain. 
\end{remark}

\subsection{Kolmogorov n-width and linear control systems}
Throughout this survey, we focus on a finite-dimensional parameter space. To complement the picture, we present here a particular example with an infinite-dimensional parameter space, given by linear, finite-dimensional, time-invariant control systems of the form 
\begin{equation}
	\label{eqn:LTI}
    \quad\begin{aligned}
        \dot{\state}(t) &= \bfA\state(t) + \bfB\inpVar(t),\\
        \outVar(t) &= \bfC\state(t),
    \end{aligned}
\end{equation}
with matrices $\bfA\in\RR^{\stateDim\times\stateDim}$, $\bfB\in\RR^{\stateDim\times\inpVarDim}$, and $\bfC\in\RR^{\outVarDim\times\stateDim}$. The symbols $\inpVar$ and $\outVar$ denote the input and output of the system, respectively, and the input $\inpVar\in\Ltwo(-\infty,0;\RR^{\inpVarDim})$ takes the role of the parameter. We assume that the control system~\eqref{eqn:LTI} is asymptotically stable, i.e., the eigenvalues of the matrix $\bfA$ are contained in the open left-half complex plane. In terms of model reduction, we focus on the approximation of the input-to-output operator
\begin{align}
	\label{eqn:HankelOperator}
	\HankelOperator\colon \Ltwo(-\infty,0;\RR^{\inpVarDim}) \to \Ltwo(0,\infty;\RR^{\outVarDim}), \qquad \inpVar \mapsto \outVar,
\end{align}
which maps past inputs to future outputs. In the literature, the operator~$\HankelOperator$ is known as the Hankel operator. Introducing the impulse response $\bfh(t) = \bfC\exp(t\bfA)\bfB$, the Hankel operator is given via the convolution integral
\begin{align*}
	(\HankelOperator \inpVar)(t) = \int_{-\infty}^0 \bfh(t-s) \inpVar(s)\ds.
\end{align*}
It is well-known that $\HankelOperator$ is a finite-rank operator of rank at most $\stateDim$; cf.~\citeasnoun[Sec.~5.4]{Ant05} or \citeasnoun[Sec.~5.1]{Fra87}. Consequently, $\HankelOperator$ is a compact operator and its singular values, commonly known as the Hankel singular values, can be computed as the square roots of the eigenvalues of the matrix product $\ctrlG\obsG$, where $\ctrlG,\obsG\in\RR^{\stateDim\times\stateDim}$ are the unique solutions of the Lyapunov equations
\begin{equation}
	\label{eqn:Lyap}
	\bfA\ctrlG + \ctrlG \bfA^\top + \bfB\bfB^\top = 0\qquad\text{and}\qquad \bfA^\top \obsG + \obsG \bfA + \bfC^\top C = 0.
\end{equation}
We sort the Hankel singular values in decreasing order and denote the $i$th Hankel singular value by $\sigma_i(\HankelOperator)$ for $i=1,\ldots,\stateDim$. For $\stateSpace = \Ltwo(-\infty,0;\RR^{\inpVarDim})$,  \citet{UngG19} prove that the Kolmogorov $\stateDimRed$-widths coincide with the Hankel singular values in the sense that 
\begin{equation}
        \label{eqn:KolmogorovHankelSVD}
        d_{\stateDimRed}(\HankelOperator(\{\inpVar\in\calZ \mid \|\inpVar\|_{\calZ} \leq 1\}),\calZ) = \sigma_{\stateDimRed+1}(\HankelOperator).
    \end{equation}

We emphasize that in contrast to the previous results, \eqref{eqn:KolmogorovHankelSVD} is not a lower or upper bound, but an exact identity that can be explicitly computed. While the computation of the $\stateDimRed$-widths, i.e., the computation of the Hankel singular values, is computationally well explored and feasible in the large-scale case, see \citeasnoun{Sim16} for a survey, the computation of an explicit reduced model of the form
\begin{align*}
    \quad\begin{aligned}
        \dot{\stateRed}(t) &= \Ared\stateRed(t) + \Bred\inpVar(t),\\
        \outVarRed(t) &= \Cred\stateRed(t),
    \end{aligned}
\end{align*}
that realizes this approximation quality is a nontrivial task. It is known in the literature as optimal Hankel norm approximation, and was first solved by \citeasnoun{Glo84}, relying on results from \citeasnoun{AdaAK71}. Extensions to descriptor systems are presented in \citeasnoun{CaoSW15} and \citeasnoun{BenW20}.

While the spaces for the solution operator in~\eqref{eqn:HankelOperator} allow for a precise analytical treatment, in practical applications one is often interested in the solution operators
\begin{align*}
	\calS_1&\colon \Ltwo(0,\infty;\RR^{\inpVarDim})\to\Ltwo(0,\infty;\RR^{\outVarDim}), & \inpVar\mapsto \outVar,\\
	\calS_2&\colon\Ltwo(0,\infty;\RR^{\inpVarDim})\to\Linfty(0,\infty;\RR^{\outVarDim}), & \inpVar\mapsto \outVar,
\end{align*}
which correspond to approximations in the Hardy $\calH_\infty$ and $\calH_2$ norm. Nevertheless, studying the Hankel operator is of great interest, as the famous a-priori error bound of balanced truncation is a bound on the $\calH_\infty$ norm in terms of the truncated Hankel singular values.

We conclude this short excursion to linear control systems with four remarks. First, the analysis in \citeasnoun{UngG19} demonstrates that the Kolmogorov $\stateDimRed$-widths for control systems can be interpreted as the dual concept to the method of active subspaces developed by \citeasnoun{Con15}. Second, if the matrices $\bfA$, $\bfB$, and $\bfC$ have a parameter dependency with a parameter $\param$ in a compact parameter set~$\paramSet$, then \citet[Thm.~3]{UngG19} prove that the Kolmogorov $\stateDimRed$-widths are lower bounded by the maximum of the Hankel singular values over the parameter set. Consequently, one can improve upon the Kolmogorov $\stateDimRed$-width by allowing parameter-dependent projection matrices, as demonstrated in \citeasnoun{WitTKAS16}, \citeasnoun{GosGU21}, and \citeasnoun{HunMMS22}. Third, an extension to bilinear control systems is presented in \citet{ZuyFB24}. Fourth, the Kolmogorov $\stateDimRed$-widths critically depend on the space $\calZ$ and the associated norm. In view of transport-dominated problems and associated methods, the standard $L_p$-norms may paint an incomplete picture, and additional metrics such as a Wasserstein distance may provide different insights \cite{Ehrlacher2020}; see also \Cref{sec:TMOR:MetricLitReview}.

\section{The elements of nonlinear model reduction}
\label{sec:nonlinMOR}
Nonlinear model reduction aims to overcome the Kolmogorov barrier by seeking nonlinear instead of linear approximations of full-model solutions. 

\subsection{Overcoming the Kolmogorov barrier with nonlinear parametrizations}
We now introduce what we mean by a parametrization and then discuss that nonlinear parametrizations can lead to nonlinear approximations that overcome the Kolmogorov barrier. 

\subsubsection{Parametrizations}
A parametrization is a map with signature
\begin{equation}\label{eq:NMOR:Param}
    \hat{q}\colon \mathbb{R}^n \times \Omega \to \mathbb{R}\,,\quad (\bftheta, \bfx) \mapsto \hat{q}(\bftheta, \bfx),
\end{equation}
so that for an $n$-dimensional weight vector $\bftheta = [\theta_1, \dots, \theta_n]^{\top} \in \mathbb{R}^n$, the function $\hat{q}(\bftheta, \cdot): \Omega \to \mathbb{R}$ is an element $\hat{q}(\bftheta, \cdot) \in \mathcal{Z}$ of a function space of interest $\mathcal{Z}$. For example, $\mathcal{Z}$ could be the space $\mathcal{Q}$ over which we defined the PDE problem \eqref{eqn:PDE} or the space $\mathcal{Q}_N$ over which the full-model solutions \eqref{eqn:PDEapprox:FOM} are defined. 
More precisely, one could first define a map $\hat{Q}: \mathbb{R}^n \to \mathcal{Z}$ to obtain $\hat{Q}(\bftheta) \in \mathcal{Z}$, which can then be evaluated at points in the spatial domain $\Omega$ over which $\mathcal{Z}$ is defined as $\hat{Q}(\bftheta)(\bfx)$. It will be convenient for us to directly work with \eqref{eq:NMOR:Param} and write  $\bftheta \mapsto \hat{q}(\bftheta, \cdot)$ if we mean $\hat{Q}$. At this point, we do not specify how to obtain a weight vector $\bftheta$ for approximating an element of a function space of interest. A map that produces a weight vector for a given function is sometimes referred to as encoder map (and correspondingly $\bftheta \mapsto \hat{q}(\bftheta, \cdot)$ is then referred to as decoder map). In model reduction, however, we typically do not have access to the functions that we would like to encode in the online phase because these are the full-model solutions. Instead, the weight vectors have to be computed via reduced dynamics in the online phase; see the forthcoming \Cref{sec:NMOR:Ingredients}. 

We call a parametrization  \eqref{eq:NMOR:Param} linear (or affine) if it depends linearly (or affinely) on the weight vector $\bftheta$, i.e., if the map $\bftheta \mapsto \hat{q}(\bftheta, \cdot)$ is linear (or affine). It is important to note that linearity refers to the dependence on $\bftheta$ and not $\bfx$, i.e., the map $\hat{q}(\bftheta, \cdot): \Omega \to \mathbb{R}, \bfx \mapsto \hat{q}(\bftheta, \bfx)$  can still be nonlinear in the spatial coordinate $\bfx$ for all $\bftheta \in \mathbb{R}^n$. Furthermore, the weight vector $\bftheta(t, \bfmu)$ can depend nonlinearly on quantities such as time $t$ and parameter $\bfmu$.

\begin{remark}
The vector $\bftheta$ is often called the parameter vector; however, we follow the convention in model reduction and refer to $\bfmu$ as the parameter and thus opted to call $\bftheta$ a weight vector, in preparation for nonlinear parametrizations such as neural networks. In the context of model reduction, the weight vector $\bftheta$ is also referred to as the reduced or latent state because, as we will see in later sections, it can represent the reduced solution. 
\end{remark}

\subsubsection{Linear parametrizations and linear model reduction}
Recall the reduced ansatz in linear model reduction given by the linear combination \eqref{eq:Prelim:RedQLinComb}. This ansatz can equivalently be written as a function of the coefficients (or reduced state \eqref{eq:Prelim:RedCoeffVector}) as 
\begin{equation}\label{eq:NMOR:LinParam}
\hat{q}(\bftheta(t, \bfmu), \bfx) = \sum_{i = 1}^{n} \theta_i(t, \bfmu) \varphi_i(\bfx)\,,
\end{equation}
which makes the linear dependence on the coefficients explicit. Notice that we identify in \eqref{eq:NMOR:LinParam} the coefficients of the linear combination \eqref{eq:Prelim:RedQLinComb} with $n$ weights $\theta_1(t, \bfmu), \dots, \theta_n(t, \bfmu)$, which are collected into the weight vector $\bftheta(t, \bfmu)$ and 
depend on time $t$ and parameter $\bfmu$. From this perspective, we can interpret the reduced ansatz~\eqref{eq:Prelim:RedQLinComb} as a linear parametrization \eqref{eq:NMOR:LinParam} of the reduced space $\mathcal{V}$, which is spanned by the reduced basis functions $\varphi_1, \dots, \varphi_n$, with the $n$-dimensional weight vector $\bftheta(t, \bfmu)$. In other words, linear model reduction represents reduced solutions with linear parametrizations. 

The set of all reduced solutions that can be represented by the linear parametrization \eqref{eq:NMOR:LinParam} is the reduced space $\Vcal$, because for every $\bftheta \in \mathbb{R}^n$ we have $\hat{q}(\bftheta, \cdot) \in \Vcal$, and conversely every element in $\Vcal$ can be written as $\hat{q}(\bftheta, \cdot)$ for some $\bftheta \in \mathbb{R}^n$. The degrees of freedom of a reduced solution are the weights in the vector $\bftheta$, which can depend on time $t$ and parameter $\bfmu$, but the reduced space $\Vcal$ given by the basis $\varphi_1, \dots, \varphi_n$ remains fixed and independent of time and parameter.

\subsubsection{Nonlinear parametrizations}
Let us now consider nonlinear parametrizations, i.e.,  parametrizations that are nonlinear in the weight vector $\bftheta$. Given a nonlinear parametrization, we can consider the set of all functions that can be reached by varying the weight $\bftheta$ as
\begin{equation}\label{eq:IRM:HatM}
    \widehat{\mathcal{M}} = \{\hat{q}(\bftheta, \cdot) \mid \bftheta \in \mathbb{R}^n\}\,,
\end{equation}
which is often referred to as trial set or trial manifold. 
While linear parametrizations lead to sets $\widehat{\mathcal{M}}$ with a vector-space structure, nonlinear parametrizations correspond to more general trial sets that are not necessarily vector spaces and may therefore represent target solutions more efficiently with respect to the number of degrees of freedom $n$. Broadly speaking, the nonlinear dependence on the weight vector~$\bftheta$ means that nonlinear parametrizations may circumvent the Kolmogorov barrier because the induced set of functions $\widehat{\mathcal{M}}$ does not necessarily have a vector-space structure anymore.   

We can interpret the nonlinear dependence on the weight vector $\bftheta$ also as allowing the representation itself to change. To illustrate this point, consider the parametrization
\begin{equation}\label{eq:NMOR:NonLinParam}
    \hat{q}(\bftheta, \bfx) = \sum_{i = 1}^r \beta_i \psi_i(\alpha_i, \bfx)\,.
\end{equation}
Here, the weight vector $\bftheta = [\bfalpha^\top, \bfbeta^\top]^\top \in \mathbb{R}^n$ with $n = 2r$ consists of $r$ weights $\bfalpha = [\alpha_1, \dots, \alpha_r]^{\top} \in \mathbb{R}^r$ that enter nonlinearly through the representation functions $\psi_1, \dots, \psi_r\colon \mathbb{R} \times \Omega \to \mathbb{R}$ and $r$ weights $\bfbeta = [\beta_1, \dots, \beta_r]^{\top} \in \mathbb{R}^r$ that enter linearly. Unlike linear parametrizations, where the representation is fixed by the basis functions
$\varphi_1, \dots, \varphi_n$, the weights $\alpha_1, \dots, \alpha_r$ in \eqref{eq:NMOR:NonLinParam}
directly influence the functions $\psi_1(\alpha_1, \cdot), \dots, \psi_r(\alpha_r, \cdot)$. As a result, the representation adapts to the specific function being approximated.
For a full-model solution $q(t, \cdot; \bfmu) \in \mathcal{M}$, the weight vector $\bftheta(t, \bfmu)$ may therefore change not only the coefficients of a fixed expansion, but also the shape of the underlying representation itself. 
This is fundamentally different from linear model reduction, where the representation corresponds to the basis functions $\varphi_1, \dots, \varphi_n$ that span a vector space, which is fixed in advance and the weights in the linear parametrization only determine the coefficients of the associated linear combination.

\begin{remark}
There is an analog to linear versus nonlinear parametrizations in machine learning. Linear parametrizations typically correspond to kernel approximations \cite{scholkopf2002learning}. Once a kernel is picked, the learning problem reduces to finding a linear combination of the kernels, which means that the representation is fixed and the coefficients (weights) of approximations in the kernel space enter linearly. In contrast, nonlinear parametrizations are given, e.g., by neural networks with at least one hidden layer and nonlinear activation functions \cite{LeCun2015,Goodfellow-et-al-2016}. The inner layers correspond to the representation and are adapted with the weights. In other words, the network adapts the basis functions themselves, which is called representation learning \cite{Rumelhart1986,10.1109/TPAMI.2013.50}.
\end{remark}

\begin{remark}
In addition to parametrizations in which the weight vector $\bftheta(t, \bfmu)$ evolves in time, one can also consider global-in-time parametrizations, in which time is treated as another input and the weight vectors do not depend on time; this is analogous in spirit to space--time discretizations. We do not pursue global-in-time parametrizations here and instead focus on sequential- or local-in-time parametrizations, where the time dependence enters via the weight vector, which aligns with the widely used model reduction paradigm in which reduced states are advanced by reduced dynamical systems. We refer the reader to \citet{24OTDDTO} for a discussion of global-in-time versus local-in-time (nonlinear) parametrizations. 
\end{remark}

\subsubsection{Example: Nonlinear parametrizations and linear advection} \label{sec:NMOR:ExampleNonLinMOROnAdvEq}
Let us revisit the linear advection equation problem of \Cref{subsec:Kolmogorov:lowerBound}. The corresponding solution manifold consists of the step functions given in \eqref{eq:NWidthExample:StepFunSol}.  
We view the step functions \eqref{eq:NWidthExample:StepFunSol} as the restriction to $[0, 1]$ of a translated reference profile, 
\begin{equation}\label{eq:LinAdvSolutionExample}
q_{\mathcal{Q}}(t, x; \mu) = q_0(x - t\mu),
\end{equation}
where $q_0\colon \mathbb{R} \to \mathbb{R}$ is defined as
\[
q_0(x) = \begin{cases}
1, & x \leq 0\,,\\
0, & x > 0\,.\end{cases}
\]
Representing the solution to the linear advection problem as  \eqref{eq:LinAdvSolutionExample}, reveals a key difficulty: the solution can be written as a composition of the translation $(t, x, \mu) \mapsto x - t\mu$ and the initial function $q_0$. A linear parametrization  \eqref{eq:NMOR:LinParam} can efficiently represent superpositions of features because it is a linear combination of basis functions. However, representing a translation of a feature via superposition means that the feature must be deactivated (low weight in the superposition) at its old spatial location and activated (high weight in the superposition) at its new spatial location. A process that can require a prohibitively large number of basis functions to describe a smooth motion of the features. 

With a nonlinear parametrization such as the one defined in \eqref{eq:NMOR:NonLinParam}, translation can be described more efficiently by composing a translation with the superposition of features. To be concrete, let us derive a nonlinear parametrization of the form \eqref{eq:NMOR:NonLinParam} that can exactly represent the step function  \eqref{eq:LinAdvSolutionExample}. One instance is 
\begin{equation}\label{eq:LinAdvExampleWithNonlinParam}
\hat{q}(\bftheta(t, \mu), x) = \beta_1(t, \mu)\psi_1(\alpha_1(t, \mu), x),
\end{equation}
where $\beta_1(t, \mu) = 1$ and $\alpha_1(t, \mu) = t\mu$, yielding the weight vector $\bftheta(t, \mu) = [\alpha_1(t, \mu), \beta_1(t, \mu)]^{\top}$. The function $\psi_1$ is given as
\[
\psi_1(\alpha, x) = q_0(x - \alpha)\,.
\]
Thus, the step functions \eqref{eq:LinAdvSolutionExample}, which include all functions on the solution manifold~\eqref{eq:LinAdvNWidthProblem:SolManifoldMCal} of the linear advection problem in \Cref{subsec:Kolmogorov:lowerBound}, can exactly be represented with the nonlinear parametrization \eqref{eq:NMOR:NonLinParam} with only $n = 2$ weights. In contrast, the error of using a linear parametrization \eqref{eq:NMOR:LinParam} cannot decay faster than with a rate $1/\sqrt{n}$ in $L^2$, as discussed in  \Cref{sec:KolmogorovNonLinExamples}. 

Even though this example is synthetic, it already shows that with a suitable choice of a nonlinear parametrization, phenomena such as transport and superposition of features can be composed, instead of having to approximate all phenomena with only a superposition of features. This additional flexibility can lead to more efficient parametrizations in terms of the number of weights.

\subsection{Extensions of Kolmogorov $n$-width to nonlinear parametrizations}\label{sec:NonLinWidth}
There are several generalizations of the Kolmogorov $n$-width concepts to nonlinear parametrizations. 
One widely used concept builds on an autoencoder-like structure: an encoder maps an element of the set of functions of interest $\calN \subseteq \calZ$ to an $n$-dimensional vector, i.e., the weight $\bftheta$ in our terminology. A decoder maps the $n$-dimensional vector back onto a function in $\mathcal{Z}$. One is then interested in the worst-case approximation error of the composition of the encoder and decoder and its decay with the dimension $n$; see, e.g., \citet{DeVore1993,DeVore_1998}. We revisit autoencoders from a numerical perspective in \Cref{sec:InstantaneousResMin}.

The existence of an encoder-decoder pair that achieves a low approximation error is insufficient for meaningful numerical computations. Essential is the regularity of the encoder and decoder map because it excludes pathological parametrizations and so relates to numerical stability. This motivates the stable manifold widths considered in \citet{Cohen2022}, which ask for at least Lipschitz bounds on the encoder and decoder. 

\citet{Cohen2023NonlinearCB} further restrict the encoder to consist of $n$ linear functionals of the solution function that is to be encoded. The decoder can remain nonlinear.  The expressivity of such nonlinear approximations is quantified by the so-called sensing numbers. In particular, \citet{Cohen2023NonlinearCB} discuss examples of step functions that are transported over time for which the error of such approximations decays faster than the error achieved with linear approximations. Indeed, it is shown in \citet{Cohen2022,Cohen2023NonlinearCB} that the sensing numbers and the stable manifold width attain the same decay rates. Numerical experiments related to sensing numbers are presented in \citet{GLAS2025550}.

\subsection{Nonlinear parametrizations with offline and online weights}\label{sec:NPM:Examples}
There is a wide range of nonlinear parametrizations used in numerical analysis, machine learning, and computational science and engineering, which also motivates nonlinear parametrizations for model reduction. We provide only a brief overview in this section and revisit them in detail in later sections of this survey. 

\subsubsection{Offline and online weights}\label{sec:NonParam:OffOnWeights}
We now consider parametrizations that can depend on an offline weight vector $\woff \in \mathbb{R}^{\offdim}$,
\begin{equation}\label{eq:NonLin:ParamWithOfflineWeight}
    \hat{q}\colon \mathbb{R}^{n} \times \Omega \times \mathbb{R}^{{\offdim}} \to \mathbb{R}\,,\qquad (\bftheta, \bfx; \bftheta_{\text{off}}) \mapsto \hat{q}(\bftheta, \bfx; \bftheta_{\text{off}})\,, 
\end{equation}
so that for any fixed $\woff$, we obtain a parametrization $\hat{q}(\cdot, \cdot; \woff)\colon \mathbb{R}^n \times \Omega \to \mathbb{R}$ in  agreement with the definition in \eqref{eq:NMOR:Param}. When the dependence on the offline weight vector $\woff$ is clear from the context, we suppress it in the notation and simply write $\hat{q}(\bftheta, \bfx)$ as in \eqref{eq:NMOR:Param}.

In the context of model reduction, the offline weight vector $\woff$ is fixed in the offline phase and thus independent of time $t$ and parameter $\bfmu$. In contrast, the weight vector $\bftheta$ plays the role of a reduced (latent) state, analogous to the coefficient or state vector \eqref{eq:Prelim:RedCoeffVector} in linear model reduction, and thus can change in the online phase and depend on time $t$ and parameter $\bfmu$. We sometimes refer to $\bftheta$ as the online weight vector to clearly distinguish it from the offline weight vector $\woff$. 

We broadly distinguish between a priori versus trained choices of the offline weights. In the trained case, the vector $\woff$ (or parts of it) is fitted to snapshot data in the offline phase, analogous to how a reduced basis spanning the reduced space $\Vcal$ is constructed from snapshots in linear model reduction.
In contrast, in the a priori case, the offline weight vector $\woff$ is either empty or prescribed without training, for instance by choosing an analytic transformation with a given shift. The online weight vector $\bftheta$ still has to be determined in the online phase.

In model reduction, we are often particularly interested in trained nonlinear parametrizations, where the offline weight vector $\woff$ is learned offline to provide an expressive parametrization so that only a low-dimensional online weight vector~$\bftheta$ needs to be computed online.

\subsubsection{Examples of nonlinear parametrizations with offline--online weights}
A common class of nonlinear parametrizations consists of parametrizations that are polynomial in the online weights. A quadratic example is
\begin{equation}\label{eq:NonLinExp:Poly}
\hat{q}(\bftheta, \bfx; \woff) = \sum_{i = 1}^n \theta_i \varphi_i(\bfx; \woff) + \sum_{1 \leq i \leq j \leq n} \theta_i\theta_j \psi_{ij}(\bfx; \woff)
\end{equation}
with functions $\varphi_i(\cdot; \woff)\colon \Omega \to \mathbb{R}$ and $\psi_{ij}(\cdot; \woff)\colon \Omega \to \mathbb{R}$ for $i, j = 1, \dots, n$. These functions (or their coefficients in a discretization) are fitted to snapshot data by fitting the offline weight vector $\woff$; we will revisit such polynomial approximations in \Cref{sec:InstantaneousResMin}.

Composing a linear expansion with a nonlinear transformation yields another important class. Let $\trafo(\cdot; \woff)\colon \mathbb{R}^r \times \Omega \to \Omega$ be a transformation that depends on the offline weight vector. Additionally, let $\varphi_1(\cdot; \woff), \dots, \varphi_r(\cdot; \woff)$ be a set of basis functions. With a decomposition of the online weight vector as $\bftheta = [\bfalpha; \bfbeta]$, with $\bfalpha \in \mathbb{R}^r$ and $\bfbeta \in \mathbb{R}^r$, one obtains
\begin{equation}\label{eq:NonLinExp:Trafo}
\hat{q}(\bftheta, \bfx; \woff) = \sum_{i = 1}^r\beta_i \varphi_i(\trafo(\bfalpha, \bfx; \woff); \woff)\,,
\end{equation}
where we assume $n = 2r$ for ease of exposition. For example, $\trafo(\alpha, x) = x - \alpha$ could model translation given by a scalar shift $\alpha$ as in \Cref{sec:NMOR:ExampleNonLinMOROnAdvEq}. Here, the offline weight  vector $\woff$ can be trained to find the basis $\varphi_1, \dots, \varphi_r$ and fix parameters of the transformation $\trafo$, while the online weight vector $\bftheta$ selects the shift and coefficients. We will revisit such nonlinear parametrizations in \Cref{sec:Template}. 

A class of nonlinear parametrizations that we have briefly touched on in~\eqref{eq:NMOR:NonLinParam} are parametrizations that adapt the representation. Let the online weight vector $\bftheta = [\bfalpha; \bfbeta] \in \mathbb{R}^n$ be decomposed into the feature vector $\bfalpha = [\alpha_1, \dots, \alpha_{n_1}]^{\top} \in \mathbb{R}^{n_1}$ and the coefficient vector $\bfbeta = [\beta_1, \dots, \beta_{n_2}]^{\top} \in \mathbb{R}^{n_2}$ with $n = n_1 + n_2$. Furthermore, let $\psi_1(\cdot, \cdot; \woff), \dots, \psi_{n_2}(\cdot, \cdot; \woff): \mathbb{R}^{n_1} \times \Omega \to \mathbb{R}$ be representation functions that depend on the offline weight vector $\woff$ and are trained on snapshot data offline. Then the nonlinear parametrization
\begin{equation}\label{eq:NonlinParam:AdaptiveInTime}
\hat{q}(\bftheta, \bfx; \woff) = \sum_{i = 1}^{n_2} \beta_i \psi_i(\bfalpha, \bfx; \woff)\,,
\end{equation}
is linear in the coefficients $\bfbeta$ for a fixed $\bfalpha$ but nonlinear in the online weight vector~$\bftheta$ through the feature vector $\bfalpha$. In particular, if $\bftheta(t, \bfmu)$ varies with time (and with parameters), then the representation 
\begin{align*}
    \psi_1(\bfalpha(t, \bfmu), \cdot; \woff), \dots, \psi_{n_2}(\bfalpha(t, \bfmu), \cdot; \woff)\colon \Omega \to \mathbb{R}
\end{align*} 
can vary with time as well. We will discuss such nonlinear parametrizations in \Cref{sec:ADEIM}. 

Another common class is given by neural networks. As an example, let us consider a fully connected feedforward neural network with a linear output layer,
\begin{equation}\label{eq:ExampleNonLinParam-NeuralNetwork}
\hat{q}(\bftheta, \bfx; \woff) = \mathcal{C}_L(\sigma(\mathcal{C}_{L-1}(\dots \sigma(\mathcal{C}_1(\bfx)) \dots )))
\end{equation}
where $\sigma$ denotes an activation function that is applied component-wise and $\mathcal{C}_{\ell}(\bfy) = \bfW_{\ell}\bfy + \bfb_{\ell}$ denotes the $\ell$-th layer with matrix $\bfW_{\ell}$ and bias $\bfb_{\ell}$ of suitable dimensions. One can then collect some of the weight matrices  $\bfW_1, \dots, \bfW_L$ and biases $\bfb_1, \dots, \bfb_L$ into the offline weight vector $\woff$ and others into the online weight vector $\bftheta$. One approach is to include all weight matrices and biases of the first $\ell$ layers into the offline weight vector, and those of the last $L - \ell$ layers into the online weight vector $\bftheta$. But there are various other variants that we will discuss in more detail in \Cref{sec:InstantaneousResMin}. 

\subsection{The elements of nonlinear model reduction}\label{sec:NMOR:Ingredients}
We identify three elements that are essential for nonlinear model reduction---the choice of the parametrization, the construction of reduced dynamics for the online weights, and the numerical solvers for efficient online evaluations---which we will revisit throughout the remainder of the survey.

\subsubsection{Overview of elements of nonlinear model reduction}
Up to this point, our discussion has focused primarily on approximation-theoretic aspects of nonlinear model reduction, in particular on how nonlinear parametrizations can overcome the limitations inherent to linear model reduction and the associated Kolmogorov barrier. While expressive nonlinear parametrizations are a crucial ingredient, they represent only one component of a successful nonlinear model reduction approach.

Beyond the choice of a nonlinear parametrization, one must construct reduced dynamics that are compatible with this parametrization, that is, evolution equations for the online weights whose solutions yield accurate reduced solutions over time and across parameters and can be solved robustly, analogous to suitable formulations and discretizations for full models. Moreover, the ultimate purpose of model reduction is to achieve lower computational costs compared to full-model simulations at a prescribed accuracy, and not only to reduce the number of degrees of freedom from $N$ to $n$. Thus, achieving $n \ll N$ alone is insufficient unless the reduced model can be evaluated efficiently in the online phase. Efficient online computations typically require both suitable reduced dynamics and efficient numerical methods for time integration and nonlinear solvers.

Guided by these considerations, we identify and distinguish three elements of nonlinear model reduction:
\begin{enumerate}
  \item[(i)] nonlinear parametrizations that provide expressive approximations and circumvent the Kolmogorov barrier,
  \item[(ii)] reduced dynamics that are compatible with the chosen parametrization and govern the evolution of the online weights in an accurate and robust manner, and
  \item[(iii)] efficient online evaluation strategies and numerical solvers that render reduced-model simulations computationally advantageous compared to full-model simulations.
\end{enumerate}

These elements mirror standard concepts of numerically solving PDEs: one chooses an approximation class (e.g., a finite element space), a formulation that defines a problem (e.g., Galerkin projection, least-squares, mixed methods), and then a numerical method for the efficient solution. 

\begin{remark}
Error estimation and certification often form an additional, distinct element in model reduction. While rigorous and efficient a posteriori error estimation is mature for linear model reduction \cite{M2AN_2005__39_1_157_0,FLD:FLD867,rozza_reduced_2008,M2AN_2008__42_2_277_0,Haasdonk2011}, it is in its early stages for nonlinear model reduction. We exemplary refer to \citet{klein2025multifidelitylearningreducedorder} and \citet{BlaSU20}. However, we note that nonlinear reduced models can be combined with the full models in multi-fidelity methods to obtain accuracy guarantees for approximations of outer-loop results; see \citet{doi:10.1137/16M1082469} for a survey on multi-fidelity methods. 
\end{remark}

\subsubsection{Element 1: Nonlinear parametrizations}
The first element is the nonlinear parametrization that defines the ansatz of the reduced solutions used to approximate the full-model solutions. Typically, one has to specify a suitable parametrization architecture $\hat{q}$ such as a polynomial ansatz~\eqref{eq:NonLinExp:Poly} or a neural network \eqref{eq:ExampleNonLinParam-NeuralNetwork} and then one determines the offline weights $\woff$ in the offline phase by training on snapshot data.

Training on snapshot data to fit the offline weights $\woff$ aims to make the parametrization capture the dominant features of the full-model solutions so that the dimension $n$ of the online weight vector $\bftheta \in\mathbb{R}^n$ can be kept low. In addition, the parametrization should satisfy suitable regularity properties, e.g., small changes in the online weight vector $\bftheta$ lead to controlled changes in $\hat{q}(\bftheta, \cdot; \woff)$, which is important for robust online computations; see also \Cref{sec:NonLinWidth}. 

More generally, the parametrization has to be chosen with the intended reduced dynamics and online evaluations in mind, because its structure can strongly affect the costs of evaluating $\hat{q}$ and of assembling quantities needed for evolving the reduced dynamics online.

\subsubsection{Element 2: Reduced dynamics}
The second element of nonlinear model reduction is to impose reduced dynamics on the weights of the nonlinear parametrization, i.e., to derive equations that determine how the weight vector $\bftheta(t,\bfmu)$ evolves in time and depends on the parameters so that $\hat q(\bftheta(t,\bfmu),\cdot)$ approximates the full-model solution in some suitable  sense.

In the setting considered throughout this survey, $\bftheta(t,\bfmu)$ is propagated over a time interval $t\in[0,T]$ for parameters $\bfmu\in\calD$ by enforcing the governing equations on the nonlinear parametrization in a reduced sense. The analog in linear model reduction is the derivation of the reduced model for the reduced coordinates; however, for nonlinear parametrizations, standard techniques such as Galerkin projection cannot be applied verbatim because the classical projection framework critically relies on the vector space structure underlying linear parametrizations.

A key aspect of nonlinear model reduction is therefore to develop  suitable reduced-dynamics formulations that are compatible with nonlinear parametrizations and support reliable online computations. In particular, the reduced dynamics critically rely on the map $\bftheta \mapsto \hat{q}(\bftheta, \cdot; \woff)$ and so their robustness depends on the regularity and conditioning of this map. In practice, it is desirable that the resulting methods can handle regions where the map becomes poorly conditioned.

\subsubsection{Element 3: Efficient online evaluation strategies and solvers}
The third element concerns the efficient online solution of the reduced dynamics, i.e., the numerical methods to advance the online weight vector $\bftheta(t, \bfmu)$ in time (or to solve the corresponding reduced algebraic problem in the steady case). 

Depending on the time discretization, this involves time integrators and, for implicit schemes, solving nonlinear systems at each time step, typically via iterative methods that can require repeated solves with Jacobian (or Jacobian-like) systems induced by the parametrization $\hat{q}$. 
For many reduced dynamics that are compatible with nonlinear parametrizations (e.g., instantaneous residual minimization), each time step can be interpreted as an optimization problem to find an update to the weight vector $\bftheta(t, \bfmu)$. Consequently, online computations in nonlinear model reduction are often closer in spirit to optimization than to classical reduced linear system solves.
Moreover, in many variational formulations, the reduced dynamics are obtained by requiring the residual to be orthogonal to test directions, which leads to inner products and norms involving the nonlinear parametrization $\hat{q}$ and often also its derivatives with respect to $\bftheta$. In the online phase, these quantities must be computed by numerical quadrature or some other numerical means. Standard techniques that apply in linear (and affine) settings to pre-compute local matrices in the offline phase and then to rapidly assemble a matrix online by superposition are typically not directly applicable anymore, because the quantities depend nonlinearly on $\bftheta$ through $\hat{q}$ and therefore cannot, in general, be assembled from pre-computed terms alone. This is closely related in spirit to the lifting bottleneck in linear model reduction. Thus, even more so than in linear model reduction, a reduction in the number of degrees of freedom from $N$ to $n$ does not by itself imply a reduction in online costs.

\section{Transformation-based methods} \label{sec:Template}
The class of nonlinear model reduction methods surveyed in this section is based on transformations. 
The central idea is to identify a suitable coordinate system in which the solution manifold becomes amenable to linear reduction. 
Rather than directly approximating the solution in the original variables, one introduces one or multiple transformations that capture dominant nonlinear features such as transport, rotation, or deformation.
A common feature of transformation-based methods is that the transformation is either motivated by physical insight about the underlying problem or, when learned from data, admits a meaningful physical interpretation a posteriori. 

\subsection{Model reduction with transformation-based methods}
The motivational example of the one-dimensional advection equation discussed in \Cref{subsec:Kolmogorov:lowerBound}, and the representation of the solution as a translation of a template in~\Cref{sec:NMOR:ExampleNonLinMOROnAdvEq}, suggests the introduction of one or multiple explicit transformation maps $\trafo$, typically taking the form of coordinate or transport maps. 
Such transformations may be prescribed a priori, as in the case of a known spatial shift, or inferred from data. Once suitable transformations are identified, there are two principal ways in which they can be incorporated into a model reduction approach:
\begin{enumerate}
	\item[(i)] The transformation may be applied directly to the solution, so that linear model reduction techniques are applied to the transformed full-model solutions $\statePDE \circ \trafo$, rather than to $\statePDE$ itself.
	\item[(ii)] Alternatively, the transformation may act on the reduced basis functions $\varphi_i$, yielding an approximation based on transformed modes, as outlined in~\eqref{eq:NonLinExp:Trafo}.
\end{enumerate}
Depending on the properties of the transformation, the two approaches may be equivalent in some cases. 
In general, however, they lead to distinct reduced models. 
This distinction becomes particularly relevant when no diffeomorphism exists that transforms the governing equations into a form that admits an efficient linear reduced representation. 
An example of such a situation arises in the modeling of pulsed detonation combustors~\cite{10.1007/978-3-319-98177-2_17}.

Moreover, transformations need not be restricted to mappings of the spatial or temporal domain. 
More general transformation operators $\shiftOper$ acting directly on the reduced basis functions can be considered, leading to transformed modes of the form $\shiftOper_i \varphi_i$. 
In all these variants, the combination of linear model reduction with transformations yields reduced models with a clear physical interpretation, as the solution is decomposed into distinct features whose evolution is described through explicit transformations tailored to the application at hand.

\subsection{Literature review of transformation-based methods}
We survey literature on nonlinear model reduction based on transformations. 

\subsubsection{Transformations via group actions and shifts}
A classical and widely explored idea is the exploitation of known symmetries of the underlying system, such as translational or rotational invariance. 
In this setting, the dynamics associated with a group action are separated from the remaining evolution of the solution field, which is then approximated using linear model reduction techniques~\cite{Kirby1992999,RowM00,ClarenceWRowley_2003,doi:10.1137/17M1126576}.
Closely related is the method of freezing introduced by \citet{BeyT04}, which forms the conceptual basis of several reduction techniques. 
For instance, \citet{OhlR13} apply this idea in the context of model reduction to the Burgers' equation and show that once the dominant transport is factored out, the remaining solution evolution exhibits low-dimensional structure and can be efficiently reduced using linear techniques.

A limitation of many symmetry-based approaches is that they typically account for a single global group action, such as a uniform translation. 
As a consequence, capturing multiple simultaneously transported structures, for instance, several waves traveling at different speeds, poses significant challenges without further extensions. 
One such extension is the shifted POD introduced in~\citet{ReiSSM18}, which represents the solution as a superposition of low-rank components, each associated with its own shift. 
This allows features moving at different velocities to be approximated using a small number of modes, while the corresponding shifts are inferred directly from snapshot data rather than prescribed a priori. Extensions, variants, and generalizations of the shifted POD are presented for instance in \citet{BlaSU20,Rei21,PAPAPICCO2022114687,BurKR25,KraMZRS25}. In a similar vein, the work \citet{doi:10.1137/17M1113679} aims to reverse transport effects by shifting solution snapshots back to a reference template. 
An iterative procedure is used to extract multiple transported structures sequentially, and linear model reduction is then applied to this aligned data. The work of \citet{CagMS19} introduces a data-driven calibration approach that trains a transformation to make the transformed solution snapshots linearly reducible. In the online phase, linear model reduction is applied to the transformed solution, and both the reduced state and the transformation are advanced, driven by the residual. 
The work \citet{GreJU22} starts with linear reduced approximations and then adds ridge terms that are learned from data and that can account for transported features. A non-intrusive approach that also uses shifts is introduced in \citet{Issan2023}. Composing linear approximations with global transport dynamics that evolve along characteristic curves of conservation laws is proposed in \citet{doi:10.1137/20M1316998}. In particular, the transport dynamics based on the characteristics can be computed efficiently in the online phase. A recent generalization of this work is presented by \citet{RimW26}. 

\subsubsection{Optimal transport and formulations in metric spaces}\label{sec:TMOR:MetricLitReview}
Another line of research on transformation-based model reduction is inspired by optimal transport theory.  \citet{PhysRevE.89.022923} compute transport maps between snapshots and use these maps to advect the reduced basis functions in time. 
A different perspective is pursued by \citet{Ehrlacher2020} and \citet{BlickhanBattisti2023}, who formulate model reduction in metric spaces, in particular in the Wasserstein space. This means that solutions are represented as measures and reduced in Wasserstein geometry via, e.g., barycenters, which can be loosely interpreted as a transformation-induced nonlinear representation that is well suited for solutions with moving features.  
Extensions of these ideas to higher-dimensional settings via sparse Wasserstein barycenters have been proposed in~\citet{Do2025}. Other approaches introduce concepts of optimal transport into deep learning methods \cite{doi:10.1137/23M1604680} and build a bridge to registration-based methods \cite{doi:10.1137/23M1570715}. Closely related is work relying on the Radon cumulative distribution transform \cite{KolPR16}, studied in the context of model reduction by, e.g., \citet{RenWM21} and \citet{LonBJSI25}; see \citet{Rim18} for computational aspects regarding the Radon transform.

\subsubsection{Registration-based methods}
Registration-based approaches compute transformations directly from data, often in an equation-agnostic and purely data-driven way inspired from image registration, with the goal of aligning dominant features across parameter values or time instances to make them linearly reducible. Such registration concepts have been proposed for model reduction in various forms. In parametric settings where the parameter $\bfmu$ induces non-smooth behavior in the solution manifold,  \citet{doi:10.1137/16M1059904,Wel20} proposes to first learn suitable transformations that align moving features, and then perform interpolation in the transformed space, where the dependence on $\bfmu$ is smoother. We note that \citet{AmsF08} already use nonlinear mappings to perform Grassmann interpolation between local reduced bases in the context of model reduction. 

The work \citet{Tad20} introduces an equation-independent variational registration procedure, which has been combined with efficient online computations in \citet{FerTZ22}.   
The work \citet{SARNA2021114168} pursues a purely data-driven approach that learns transformations and combines it with a regression-based approach to compute reduced states online. Several works refine how transformations are constructed and deployed online:  \citet{Grundel2024} combine transformations with adaptive reduced meshes to efficiently compute reduced approximations. \citet{SARNA2024104065}  compute transformations via monotonic feature matching of discontinuities, while \citet {10.1007/978-3-031-86169-7_42} uses moving meshes. Registration can also be carried out in a space-time setting by learning a mapping on the combined space-time domain as in \citet{TaddeiZhang2021,kleikamp2022nonlinear}. 
Other work on registration-based methods include~\citet{MonicaNonino2023}, who explicitly target fluid-structure interaction problems, as well as \citet{iollo2026mathematicalaspectsregistrationmethods} on the foundations of registration on bounded domains.

\subsubsection{Lagrangian formulations and other transformations}
Transformation-based model reduction can also be formulated in a Lagrangian framework, where the dynamics are described in coordinates that move with the flow. 
The Lagrangian basis method introduced by~\citet{mojgani2017lagrangianbasismethoddimensionality} is an early example of this perspective. A driven-approach that extends dynamic mode decomposition to Lagrangian settings is presented by~\citet{LU2020109229,Lu_2021}. 
Finally, more problem-specific transformations have been proposed that explicitly identify shock curves or moving interfaces: \citet{Tommasao} explicitly track the shock curve to leverage piecewise-smooth components, while  \citet{RAZAVI2025113576} and \citet{doi:10.2514/6.2025-1571} compute registration maps that align shocks. \citet{ALIREZAMIRHOSEINI2023111739} also aligns dominated features but does so implicitly via an online optimization problem. 

\subsection{Shifted proper orthogonal decomposition (shifted POD)}
We now discuss in more detail the transformation-based method shifted POD, which was introduced in \citet{ReiSSM18} and further developed in, e.g., \citet{10.1007/978-3-319-98177-2_17,BlaSU20,Rei21,PAPAPICCO2022114687,Sch23,BurKR25,KraMZRS25}. Despite its name, the method is not tied to a shift transformation and thus we present it with general transformations. We specifically discuss the elements of nonlinear model reduction for shifted POD and conclude with a short discussion on the design of the transformation.

\subsubsection{Shifted POD: Motivational example}
To motivate shifted \POD, we discuss here again our prototypical transport-dominated problem, namely the one-dimensional advection equation; cf.~\Cref{subsec:Kolmogorov:lowerBound}. In particular, we consider the following instance of the advection equation PDE problem: set the time domain to $\timeInt = (0,1)$ and the spatial domain to $\Omega = (0,1)$, and consider 
\begin{align}
\label{eqn:advectionEquation:periodicBC}
	\begin{aligned}
	\partial_t \statePDEQ(t, \bfx) &= \partial_{\bfx} \statePDEQ(t,\bfx), & (t,\bfx) &\in \timeInt\times\Omega,\\
	\statePDEQ(t,0) &= \statePDEQ(t,1), & t&\in\timeInt,\\
	\statePDEQ(0,\bfx) &= \statePDEQz(\bfx), & \bfx&\in\Omega,
	\end{aligned}
\end{align}
with periodic boundary conditions but without a $\bfmu$ dependence. The solution can be represented as $\statePDEQ(t,\bfx) = \statePDEQz(\bfx-t)$ for $(t,\bfx)\in\timeInt\times\Omega$. 
Although we have already seen that the Kolmogorov $\stateDimRed$-widths decay at best only polynomially for this problem class, we explicitly compute the best linear approximation of the solution. The best linear basis $\{\reducedBasis_1, \ldots, \reducedBasis_{\stateDimRed}\}$ can be computed by solving the optimization problem
\begin{align}
	\label{eqn:POD:minimization}
	\begin{aligned}
		&\min \tfrac{1}{2} \int_\timeInt \|\statePDEQ(t,\cdot) - \sum_{j=1}^{\stateDimRed} \langle \statePDEQ(t,\cdot), \reducedBasis_j\rangle_{\pivotSpace} \reducedBasis_j\|^2_{\pivotSpace} \dt,\\
		&\text{s.t. } \reducedBasis_i\in \PDEspace \text{ and } \langle \reducedBasis_i, \reducedBasis_j\rangle_{\pivotSpace} = \delta_{ij} \text{ for } i,j=1,\ldots,\stateDimRed,
	\end{aligned}
\end{align}
where $\delta_{ij}$ denotes the Kronecker delta and $\pivotSpace$ a Hilbert space, e.g., the pivot space. For the specific example of the advection equation~\eqref{eqn:advectionEquation:periodicBC} we choose the Hilbert spaces $\PDEspace = H^1_{\mathrm{per}}(\Omega)$ and $\pivotSpace = \Ltwo(\Omega)$. The optimization problem~\eqref{eqn:POD:minimization} is know as the POD minimization problem, cf.~\citet{sirovich1,GubV17}. The optimization problem~\eqref{eqn:POD:minimization} has a closed-form solution, as the basis functions are given as the $\stateDimRed$ leading eigenfunctions of the nonnegative, self-adjoint, compact operator
\begin{align*}
	\calR\colon \pivotSpace\to\pivotSpace,\qquad \calR\reducedBasis = \int_{\timeInt} \langle \statePDEQ(t,\cdot), \reducedBasis\rangle_{\pivotSpace} \statePDEQ(t,\cdot) \dt,
\end{align*}
see, e.g., \citeasnoun{GubV17}. For the particular choice of the advection equation as in~\eqref{eqn:advectionEquation:periodicBC}, one can show, see \citet{holmes,BlaSU20}, that the optimal basis functions are given by Fourier modes, and thus are not tailored to the concrete solution, which only depends on the initial function $\statePDEQz$. In fact, this observation serves as an explanation to the oscillatory behavior of POD-based linear reduced models for transport dominated problems. 

\subsubsection{Shifted POD: Transformed modes}
The \POD optimization problem~\eqref{eqn:POD:minimization} leads to the best linear approximation of the given data. The idea from \citeasnoun{ReiSSM18} is to introduce multiple shift operators $\shiftOper_i$ to account for different waves inherent to the solution, and so to mitigate the issue with transport-dominated problems as in the motivational example. The shift operators are parametrized by an additional parameter $\pathVar_i$, called the path variable, that amounts to the shift amount and, in multiple dimensions, the shift direction for each mode. To simplify the presentation, we restrict ourselves to one-dimensional path variables here, i.e., $\pathVar_i\in\RR$ for $i=1,\ldots,\stateDimRed$ and refer to \citeasnoun{BlaSU20} for the vector-valued setting. One then seeks a reduced approximation of the form
\begin{equation}
	\label{eq:NMOR:transformPODParam}
	\hat{q}(t, \bfx, \param) = \sum_{i = 1}^{\stateDimRed} \coeff_i(t, \param) \left[\shiftOper_i(\pathVar_i(t,\param))\varphi_i\right](\bfx),
\end{equation}
where $\shiftOper_i$ is a shift operator acting on the basis function $\varphi_i$ and the resulting function $\left[\shiftOper_i(\pathVar_i(t,\param))\varphi_i\right]$ is the transformed mode that is then evaluated at the spatial coordinate.
If the path variables $\pathVar_1(t, \bfmu), \dots, \pathVar_{\stateDimRed}(t, \bfmu)$ are not fixed a priori, then~\eqref{eq:NMOR:transformPODParam} is a nonlinear parametrization, which constitutes the first element of nonlinear model reduction as outlined in \Cref{sec:NMOR:Ingredients}. In the language of \Cref{sec:NPM:Examples}, the online weights are therefore given by the coefficients $\hat{q}_1(t, \bfmu), \dots, \hat{q}_n(t, \bfmu)$ as well as the path variables $\pathVar_1(t, \bfmu), \dots, \pathVar_n(t, \bfmu)$. 
We emphasize that the method is not limited to $\shiftOper_i$ being a shift operator. In fact, any parameter-family of operators can be used, which is why we refer to $\shiftOper_i$ as a transformation operator. Typical requirements include that $\pathVar_i \mapsto \shiftOper_i(\pathVar_i)\varphi_i$ is continuous, the transformation operators are $\PDEspace$-invariant, and that $\shiftOper_i$ is uniformly bounded.

To generalize the \POD optimization problem~\eqref{eqn:POD:minimization} to the decomposition~\eqref{eq:NMOR:transformPODParam} with transformed modes, it is important to observe that in general there is no guarantee that $\{\shiftOper_i(\pathVar_i)\varphi_i\}_{i=1}^{\stateDimRed}$ forms a set of orthonormal basis functions. Consequently, inspired by \citet{ReiSSM18}, \citet{10.1007/978-3-319-98177-2_17} and \citet{BlaSU20} pose the transformed \POD optimization problem as
\begin{align}
	\label{eqn:shiftedPOD:minimization}
	\quad\begin{aligned}
		&\min \tfrac{1}{2} \int_\timeInt \|\statePDEQ(t,\cdot) - \sum_{i=1}^{\stateDimRed} \coeff_i(t, \param) \left[\shiftOper_i(\pathVar_i(t,\param))\varphi_i\right](\bfx) \|^2_{\pivotSpace} \dt,\\
		&\text{s.t. } \varphi_i\in \PDEspace \text{ and } \|\varphi_i\|_{\pivotSpace} = 1 \text{ for } i,j=1,\ldots,\stateDimRed,
	\end{aligned}
\end{align}
which has to be minimized over the coefficient functions $\coeff_i$, the path variables $\pathVar_i$, and the modes $\varphi_i$. Using a suitable functional analytic setting including assumptions on the transformation operators $\shiftOper_i$ and appropriate spaces and bounds for the coefficient functions $\coeff_i$ and the path variables $\pathVar_i$, one can show that the shifted \POD optimization problem~\eqref{eqn:shiftedPOD:minimization} is solvable; cf.~\citeasnoun{BlaSU20} and \citeasnoun{BlaSU22}. In comparison with the \POD optimization problem~\eqref{eqn:POD:minimization}, which can be obtained from~\eqref{eqn:shiftedPOD:minimization} by setting $\shiftOper_i$ to be the identity, the analysis of the properties of~\eqref{eqn:shiftedPOD:minimization} is more delicate; cf.~\citeasnoun{BlaSU22}. A notable exception is given, if all transformation operators and their path variables are forced to be identical, i.e., the approximation~\eqref{eq:NMOR:transformPODParam} takes the form
\begin{equation}
	\label{eq:NMOR:transformPODParam:singleShift}
	\hat{q}(t, \bfx, \param) = \sum_{i = 1}^{\stateDimRed} \coeff_i(t, \param) \left[\shiftOper(\pathVar(t,\param))\varphi_i\right](\bfx).
\end{equation}
If we additionally assume that $\shiftOper$ is an isometry, then the corresponding minimization problem can be solved by first transforming snapshot data accordingly and then applying standard POD to transformed snapshots. For further details we refer to \citeasnoun{BlaSU20}; see also \cite{ReiSSM18,10.1007/978-3-319-98177-2_17}, and \cite{CagMS19} for related transformations and their computations. The specific case where the transformation operators are given by ridge functions is discussed in \citeasnoun{GreJU22}.

\subsubsection{Shifted POD: Computational aspects for the optimization problem}
In general, there is no closed-form solution to~\eqref{eqn:shiftedPOD:minimization} available, and an approximate minimizer has to be computed numerically. Besides a spatial and temporal discretization, this also requires numerically solving the optimization problem, for instance by employing a gradient-based method as in \citeasnoun{BlaSU22}. As discussed in \citet{Sch23}, even for a moderate number $\stateDimRed$ this becomes computationally expensive, as the degrees of freedom scale with the number of temporal and spatial degrees of freedom. Instead of solving the full optimization problem, one can also use a variable projection approach as discussed in detail in \citet{Sch23}. A good initialization is often important in practice. This can for instance be obtained by initializing the path variables so that they suitably align the snapshots; see also the techniques discussed in \citet{CagMS19,Mendible2020,Rei21}. Given an initialization for the path variables, the modes can be initialized by applying the transformation operators with the path variables and infer modes as the leading singular vectors of the transformed snapshots, as for instance done in the case of a wildland fire simulation in \citet{BlaSU21b}.

\subsubsection{Shifted POD: Construction of the reduced model}
\label{subsec:shiftedPOD:ROM}
The second element of nonlinear model reduction constitutes the reduced dynamics. We start our discussion by assuming that the path variables $\pathVar_i$ are fixed, which, assuming a suitable parametrization, means that we can consider them as offline weights; cf.~\Cref{sec:NPM:Examples}. In this case, the transformed modes ansatz~\eqref{eq:NMOR:transformPODParam} can be used in a classical Galerkin projection framework similarly as in~\eqref{eq:Prelim:HatFGalerkin}. In more detail, recall the PDE problem specified in~\eqref{eqn:PDE} and 
the reduced state $\stateRedSPOD(t, \param) = [\coeff_1(t, \param), \dots, \coeff_{\stateDimRed}(t, \param)]^{\top} \in \RR^{\stateDimRed}$. We now introduce the reduced path vector as 
\begin{align*}
	\pathVarRed(t, \param) = [\hat{\pathVar}_1(t, \param), \dots, \hat{\pathVar}_\stateDimRed(t, \param)]^{\top} \in \RR^{\stateDimRed}\,,
\end{align*}
which, at the moment, we assume to be given.
The associated Galerkin reduced model is 
\begin{align}
	\label{eqn:SPOD:ROM}
	\massMatrix_1(\pathVarRed(t,\param))\dot{\stateRedSPOD}(t,\param) + \massMatrixShift(\pathVarRed(t,\param))\bfD(\stateRedSPOD(t,\param))\dot{\pathVarRed}(t,\param) = \fred_1(\stateRedSPOD(t,\param),\pathVarRed(t,\param);\param)
\end{align}
with 
\begin{align*}
	\massMatrix_1(\pathVarRed) &= \left[\langle \shiftOper_i(\hat\pathVar_i)\varphi_i, \shiftOper_j(\hat\pathVar_j)\varphi_j\rangle_{\pivotSpace}\right]_{i,j=1}^{\stateDimRed} \in \RR^{\stateDimRed\times\stateDimRed}\,,\\
	\massMatrixShift(\pathVarRed) &= \left[\langle \shiftOper_i(\hat\pathVar_i)\varphi_i, \shiftOper_j'(\hat\pathVar_j)\varphi_j\rangle_{\pivotSpace}\right]_{i,j=1}^{\stateDimRed} \in \RR^{\stateDimRed\times\stateDimRed}\,,\\
	\bfD(\stateRedSPOD) &= \diag(\coeff_1,\ldots,\coeff_{\stateDimRed})\in\RR^{\stateDimRed\times\stateDimRed}\,,\\
	\fred_1(\stateRedSPOD,\pathVarRed;\param) &= \left[\left\langle \shiftOper_i(\hat\pathVar_i)\varphi_i, f\left(\cdot, \sum\nolimits_{k=1}^{\stateDimRed} \coeff_k \shiftOper_k(\hat\pathVar_k)\varphi_k; \param\right) \right\rangle_{\pivotSpace}\right]_{i=1}^{\stateDimRed}\in\RR^{\stateDimRed},
\end{align*}
where $\shiftOper_i'$ denotes the derivative of $\shiftOper_i$. 

If the reduced path variables are unknown, i.e., if they constitute online weights, then~\eqref{eqn:SPOD:ROM} constitutes an underdetermined differential--algebraic equation; see \citeasnoun{MehU23} for a recent survey. To render the underdetermined equation~\eqref{eqn:SPOD:ROM} well-posed, we have to add additional equations, e.g., of the form
\begin{align}
	\label{eqn:phaseCondition}
	\boldsymbol{\psi}(\pathVarRed,\dot{\pathVarRed},\stateRedSPOD,\dot{\stateRedSPOD};\param) = 0,
\end{align}
which are called \emph{phase conditions} or \emph{reconstruction equations} in the literature; see for instance \citeasnoun{BeyT04, OhlR13, ClarenceWRowley_2003, RowM00}. A priori, it is unclear what the optimal choice for the phase condition~\eqref{eqn:phaseCondition} is and several instances are reported in the literature, for instance by minimizing the temporal change of the reduced  coefficients. We refer to \citeasnoun{Sch23} for an overview.  \citeasnoun{BlaSU20} propose to obtain the phase condition via a Petrov--Galerkin projection, where the test space is chosen as
\begin{align*}
	\mathrm{span}\{\coeff_1\shiftOper_1'(\pathVar_1)\varphi_1, \ldots, \coeff_{\stateDimRed}\shiftOper'_\stateDimRed(\pathVar_\stateDimRed)\varphi_\stateDimRed\}.
\end{align*}
With this particular choice, the phase condition reads
\begin{align}
	\boldsymbol{\psi}(\pathVarRed,\dot{\pathVarRed},\stateRedSPOD,\dot{\stateRedSPOD};\param) = \bfD(\stateRedSPOD)^\top \left(\massMatrixShift(\pathVarRed)^\top \dot{\stateRedSPOD} + \massMatrix_2(\pathVarRed)\bfD(\stateRedSPOD)\dot{\pathVarRed} - \fred_2(\stateRedSPOD,\pathVarRed;\param)\right)
\end{align}
with
\begin{align*}
	\massMatrix_2(\pathVarRed) &= \left[\langle \shiftOper_i'(\hat\pathVar_i)\varphi_i, \shiftOper_j'(\hat\pathVar_j)\varphi_j\rangle_{\pivotSpace}\right]_{i,j=1}^{\stateDimRed} \in \RR^{\stateDimRed\times\stateDimRed}\,,\\
	\fred_2(\stateRedSPOD,\pathVarRed;\param) &= \left[\left\langle \shiftOper_i'(\hat\pathVar_i)\varphi_i, f\left(\cdot, \sum\nolimits_{k=1}^{\stateDimRed} \coeff_k \shiftOper_k(\hat\pathVar_k)\varphi_k; \param\right) \right\rangle_{\pivotSpace}\right]_{i=1}^{\stateDimRed}\in\RR^{\stateDimRed}.
\end{align*}
The coupled reduced model for the reduced state $\stateRedSPOD$ and the reduced path variable $\pathVarRed$ is thus given by (ignoring the arguments)
\begin{align}
	\label{eqn:shiftedPOD:fullROM}
	\begin{bmatrix}
		\bfI_{\stateDimRed} & 0\\
		0 & \bfD(\stateRedSPOD)^\top
	\end{bmatrix}\begin{bmatrix}
		\massMatrix_1(\pathVarRed) & \massMatrixShift(\pathVarRed)\\
		\massMatrixShift(\pathVarRed)^\top & \massMatrix_2(\pathVarRed)
	\end{bmatrix}\begin{bmatrix}
		\bfI_{\stateDimRed} & 0\\
		0 & \bfD(\stateRedSPOD)
	\end{bmatrix}\begin{bmatrix}
		\dot{\stateRedSPOD}\\
		\dot{\pathVarRed}
	\end{bmatrix} &= \begin{bmatrix}
		\fred_1(\stateRedSPOD,\pathVarRed;\param)\\
		\bfD(\stateRedSPOD)^\top \fred_2(\stateRedSPOD,\pathVarRed;\param)
	\end{bmatrix}.
\end{align}
It can be shown, that this reduced model minimizes the residual in a certain sense, cf.~\citeasnoun{BlaSU20}, which was termed \emph{continuous optimality} in \citet{CARLBERG2017693,LeeC20}, and already used by \citeasnoun{MilM81} in the context of moving finite element methods. We emphasize that this is a particular instance of the Dirac--Frenkel variational principle, which will be discussed in more general settings in \Cref{sec:InstantaneousResMin}.

In general, there is no guarantee that the matrix on the left-hand side of~\eqref{eqn:shiftedPOD:fullROM} is non-singular, in which case a regularization might be required; see also the forthcoming \Cref{sec:RegDF}. Indeed, whenever a coefficient $\coeff_i$ of the reduced state attains the value zero, then one cannot solve for the corresponding path variable $\hat{\pathVar}_i$. To mitigate this issue it is favorable to enforce the same path variable and same transformation operator for multiple reduced coefficients. In the case of a single transformation operator and a single path variable, this amounts to the representation~\eqref{eq:NMOR:transformPODParam:singleShift}. The general case is presented in detail in \citeasnoun{Sch23}. In the case that the right-hand side of the PDE is linear, one can use semigroup theory, e.g., \citeasnoun{Paz83}, to construct an a posteriori error estimator \cite{BlaSU20}. We emphasize that the resulting error estimator is quite generic and can also be applied to certify the prediction error of a physics-informed neural network~\cite{RAISSI2019686}; see~\citeasnoun{HilU25a,HilU25b}.

\subsubsection{Shifted POD: Offline--online decomposition}
We now turn to the third element on nonlinear model reduction (cf.~\Cref{sec:NMOR:Ingredients}), i.e., which in this case amount to an efficient evaluation of the reduced model~\eqref{eqn:shiftedPOD:fullROM}. Due to the nonlinear ansatz, all matrices and the right-hand sides, formally depend on the full-model dimension, and might not be easily pre-computable. A notable exception is given, if a single transformation operator with a single path variable is used as in~\eqref{eq:NMOR:transformPODParam:singleShift}, and the transformation operator is unitary, which is the case for a one-dimensional shift operator with periodic boundary conditions. In this case, most of the inner products are independent of the full-model dimension and can be precomputed and efficiently evaluated; see \citeasnoun[Thm.~4.3.1]{Sch23}. In very specific cases, this is even true for the general ansatz~\eqref{eq:NMOR:transformPODParam}, as demonstrated by \citeasnoun{BlaSU21a} for the one-dimensional wave equation. In general, this is however not the case and the costs of the numerical evaluation of~\eqref{eqn:shiftedPOD:fullROM} scales with the dimension of the full model. 

To mitigate this issue, different strategies are proposed in \citeasnoun{BlaSU21b}, differentiating between the inner products for the matrices on the left-hand side of~\eqref{eqn:shiftedPOD:fullROM} and the evaluation of the inner products corresponding to the potentially nonlinear full-model right-hand side. As the matrices on the left-hand side of~\eqref{eqn:shiftedPOD:fullROM} only depend on $\pathVarRed$, \citeasnoun{BlaSU21b} propose to sample these matrices for typical values of $\pathVarRed$ and then assembling the matrices via interpolation of the sampled matrices. A similar strategy can also be employed if the full-model right-hand side~$f$ is linear and we refer to \citeasnoun{CagMS19} and \citeasnoun{Sch23} for details. Exploiting properties of the transformation operator and using the method of active subspaces presented in \citeasnoun{Con15} might be used to lower the dimension of the space to be sampled. If $f$ is nonlinear, then an extension of empirical interpolation as discussed in \Cref{sec:EmpiricalRegression} to account for the transformation operators can be applied. Corresponding ideas are presented for instance in \citeasnoun{BlaSU21b, doi:10.1137/20M1316998,Grundel2024}. An additional opportunity for more speedup is obtained by applying time-integrators with adaptive step-sizes, as the special ansatz~\eqref{eq:NMOR:transformPODParam} often allows for large time steps compared to linear reduced models. 

An alternative approach to the numerical integration of the reduced model~\eqref{eqn:shiftedPOD:fullROM} is given by learning the reduced coefficients and reduced paths via machine learning. In the case of transformed modes, this is pursued in \citeasnoun{PAPAPICCO2022114687} and \citeasnoun{BurKR25}.

\subsubsection{Shifted POD: Transformation operators}
The question that remains to be answered, is the design of the transformation operators $\shiftOper_i$. If a one-dimensional setting with periodic boundary conditions, a simple shift operator might be sufficient. Nevertheless, the situation changes for non-periodic boundary conditions or more involved wave propagation in two or higher-dimensional computational domains. The case of non-periodic boundary conditions is discussed in \citeasnoun[Sec.~3.3]{Sch23}; for instance by constant extrapolation, virtually extending the computational domain, or by zero padding. In general, it is desirable to construct the transformation operators based on expert knowledge on the application at hand, for instance by studying properties of the semigroup underlying the respective problem. In particular, the semigroup approach suggests the construction of transformation operators that are equivariant with respect to the full model right-hand side $f$. We emphasize that incorporating such expert knowledge can be the key to have high-fidelity approximations even far outside the training data. If no such expert knowledge is available, then one can resort to learning the transformation operators, as for instance demonstrated in \citeasnoun{KraBHR22}.

\section{Online adaptive model reduction}\label{sec:ADEIM}
This section surveys model reduction methods that adapt the reduced space over time during the online phase. The underlying motivation is that solutions to transport-dominated problems can often be approximated well in reduced spaces locally in time, even when approximations in a single, global reduced space are inefficient. 

\subsection{Model reduction with online basis updates}
We consider online adaptive model reduction methods, which represent reduced solutions as linear combinations of basis functions that evolve with time $t$ (and possibly the parameter $\bfmu$). Correspondingly, the reduced space spanned by the basis functions evolves over time. 

One common way to realize time-dependent basis functions is through a nonlinear parametrization in which the basis depends on time-varying features that are contained in the online weights $\bftheta(t, \bfmu)$, as in \eqref{eq:NonlinParam:AdaptiveInTime}. The online phase then includes an adaptation step for the weight vector $\bftheta(t, \bfmu)$ that evolves these features (and thus the basis) over time.  In the semi-discrete setting, the time-dependent reduced space $\Vcal(t)$ is represented by the time-dependent basis matrix $\bfV(t) \in \mathbb{R}^{N \times n}$, possibly also depending on the parameter $\bfmu$, which is then updated online. We remark that some of the transformation-based methods can also be interpreted as time-dependent basis methods where, for example, shifts are evolving the basis over time.

The principle underlying time-dependent basis methods that can  overcome the Kolmogorov barrier is local-in-time approximations: over short time windows, the full-model solutions can be well approximated in a low-dimensional vector space, but the space needs to be adapted over time as the solutions evolve. Unlike in linear model reduction with a single reduced space that is fixed over time, an online adaptive reduced space can approximate solutions with moving coherent structures and transport without requiring a high-dimensional online weight vector $\bftheta(t, \bfmu)$.

A consequence of the online adaptation is that the offline--online decomposition of linear model reduction is abandoned or at least weakened. 
In particular, the basis matrix $\bfV(t)$ (or equivalently the corresponding basis functions) is adapted while the reduced model is simulated online; these basis updates can incur costs that scale with the dimension $N$ of the full-model solution space, thereby limiting speedups. A central challenge is therefore to design basis-update mechanisms that limit the $N$-dependence, for example via sparse basis updates. 
At the same time, this weakening of the offline--online decomposition enables the reduced space to adapt online and capture dynamics and solution features that were not present in the training data used during the offline phase. This can be of particular importance when the parameter $\bfmu$ is high-dimensional, and the corresponding parameter domain~$\Dcal$ cannot be sampled exhaustively, or when the full model is so expensive to simulate that not many training data trajectories can be generated. 
We focus on online adaptive model reduction to overcome the Kolmogorov barrier rather than other advantages, such as avoiding a comprehensive offline training phase, for which we refer to a discussion in \citet{Peherstorfer15aDEIM} and \Cref{rm:NonTimeAdapt}.

\begin{remark}\label{rm:NonTimeAdapt}
There is a line of work on adapting reduced spaces over outer-loop iterations, such as optimization, as in, e.g., \citet{ArianFahlSachs2000_TRPOD,KunischPODAdapt,https://doi.org/10.1002/nme.4770,Amsallem2015,doi:10.1137/16M1081981}. The purpose of progressively building reduced spaces is often that the parameter domain cannot be exhaustively sampled, and so reduced models are built during outer-loop iterations. Instead, in this survey, we focus on adaptation over time to mitigate the Kolmogorov barrier.
\end{remark}

\subsection{Overview of online adaptive model reduction and time-dependent basis methods}\label{sec:ADEIM:LitOverview}
Let us set online adaptive model reduction in context to broader time-dependent basis methods, i.e., methods where the approximation space evolves with time rather than staying fixed as in linear model reduction.

\subsubsection{Localized model reduction}
Localized model reduction methods rely on  multiple low-dimensional spaces, with a selection or weighting mechanism to use them online. For example, a common approach is to pre-compute multiple reduced spaces offline and then select one online with an indicator function based on, e.g., the current reduced solution, parameter, or time. Selecting the reduced space based on online information results in a nonlinear parametrization. The key is that the decomposition into local reduced spaces is performed so that each local reduced space corresponds to a local regime that can be well approximated by a low-dimensional reduced space, even if a single global reduced space of the same dimension would lead to inaccurate approximations.  

One way to find a decomposition into local regimes is by tailoring the local reduced spaces to regions in the parameter domain $\paramSet$ \cite{doi:10.1137/090780122,Eftang01082011,Haasdonk01082011,https://doi.org/10.1002/nme.3327,Wieland2015,NonLinMORLocal}. Another way is to decompose the time interval and use different local reduced spaces for each time window   \cite{DihlmannDrohmannHaasdonk2011AdaptiveTimePartitioning,10.1007/978-3-642-20671-9_39,SHIMIZU2021114050,PARISH2021109939,COPELAND2022114259}, which can be interpreted as an approximation with a switched or hybrid system. The reduced spaces can also be localized directly based on the state space (in the semi-discrete setting) \cite{doi:10.2514/6.2012-2686,https://doi.org/10.1002/nme.4371,doi:10.1137/130924408,Amsallem2015,https://doi.org/10.1002/nme.6603,HESS2019379,Amsallem2016,10.1098/rsta.2021.0206}. One approach to do so is via clustering: A clustering method is applied to the snapshots in the offline phase so that snapshots that are similar are grouped together into clusters and a local reduced space is constructed for each cluster of snapshots. 
Localized model reduction can also be viewed as yielding nonlinear parametrizations that are induced by a union of multiple reduced spaces as in \citet{doi:10.1137/20M1380818}.  Motivated by domain decomposition and multi-scale methods, \citet{doi:10.1137/151003660, 10.1007/978-3-319-57394-6_38,doi:10.1137/17M1138480,doi:10.1137/22M148402X} localize over the spatial domain and compute local reduced spaces for each subdomain, which is specifically designed to target problems with a multi-scale structure. A central development of these works is rigorous error control, which is then leveraged for adaptive online enrichment \cite{10.1007/978-3-319-57394-6_38}. If the domain can be partitioned in such a way that almost the same localized features can be expected in each subdomain, then \citet{VanES26} propose to use a single local basis. Localization can also be used to support local modifications in the online phase, so that updates and recomputation remain local and lightweight \cite{StaticCondRBM,doi:10.1137/15M1009603,doi:10.1137/15M1054213}.

Dictionary approaches push the idea of selecting a local reduced space further by pre-computing a large set of candidate basis functions offline and then selecting a few of the candidate basis functions online to span a reduced space \cite{KaulmannHaasdonk2013OnlineGreedy,Abgrall2016,Daniel2020,Balabanov2021,Nouy2024,Herkert2024}.

\subsubsection{Dynamic low-rank approximations and related methods}
\citet{doi:10.1137/050639703} introduce dynamic low-rank approximations, which parametrize the solution via a low-rank matrix and derive reduced dynamics for how the factors of the low-rank matrix representation evolve over time. The key principle underlying the reduced dynamics of dynamic low-rank approximations is the Dirac--Frenkel variational principle \cite{dirac1930note,frenkel1934wave}, which projects the full dynamics onto the tangent space of the manifold of matrices of a fixed rank; we  will revisit the Dirac--Frenkel variational principle for general nonlinear parametrizations in \Cref{sec:InstantaneousResMin}. 

A challenge of dynamic low-rank approximations is that the reduced dynamics of the low-rank factors can become poorly conditioned due to numerically small singular values of the low-rank matrix representations. Thus, applying standard integrators can suffer from strong step-size restrictions. 
Correspondingly, tailored time integration schemes for dynamic low-rank approximations have been developed. There are projector-splitting integrators introduced by \citet{Lubich2014,doi:10.1137/15M1026791}, which are robust to the small singular values in the low-rank matrix representations. An alternative is offered by  integrators called basis-update and Galerkin (BUG) introduced by \citet{Ceruti2022,Ceruti2024}, which can be computationally more efficient. In general, efficient, stable, and robust time integration of the reduced dynamics corresponding to dynamic low-rank approximations remains an active research direction, e.g., \cite{Ceruti2022b,doi:10.1137/22M1519493,Hochbruck2023,doi:10.1137/23M1565103,lam2024randomizedlowrankrungekuttamethods,NOBILE2025117495,Bachmayr2025,CARREL2026114421}. 

Another active research direction is structure preservation. For example, for Hamiltonian wave equations, a symplectic dynamic low-rank approximation approach was developed by \citet{Musharbash2020}. For kinetic problems, in particular for Vlasov-type equations, there is a line of work on enforcing conservation laws in dynamic low-rank approximations  \cite{EINKEMMER2021110353,COUGHLIN2024113055,EINKEMMER2023112060,EINKEMMER2021110495,doi:10.1137/18M1218686,HesthavenPagliantiniRipamonti2024AdaptiveSymplecticMOR}.

Dynamic low-rank approximations have also been extended to tensor decompositions by \citet{doi:10.1137/09076578X,doi:10.1137/120885723,Arnold2014,doi:10.1137/140976546}. Additionally, there are methods to evolve bases that specifically represent only parts of the full-model dynamics such as instability directions  \cite{10.1098/rspa.2015.0779}. 

Dynamic low-rank approximations are related to the dynamic orthogonal decomposition for stochastic (partial) differential equations, which was introduced by \citet{SAPSIS20092347}. There is a line of work on the error analysis for dynamic orthogonal approximations \cite{NobileDO,MUSHARBASH2018135,7b835857089b4f798baac614569c46cc} and for time integration \cite{doi:10.1137/16M1109394,doi:10.1137/22M1534948,doi:10.1137/21M1431229}. Furthermore, there are generalizations to Petrov--Galerkin formulations for the random PDE setting \cite{10.1007/978-3-031-86169-7_45,NOBILE2025117495}.

We can provide only a very brief overview of dynamic low-rank approximation methods and their related methods. We refer to the dedicated survey by \citet{EINKEMMER2025114191} for more details.

\subsubsection{Time-adaptive basis in model reduction}
Let us now focus on time-adaptive basis methods  developed explicitly for model reduction. For a survey dedicated to model reduction of time-dependent problems, see \citet{hesthaven_pagliantini_rozza_2022}. 

One line of work on online adaptive reduced spaces derives the basis evolution by mimicking the transport dynamics. For example, \citet{PhysRevE.89.022923} build on ideas from optimal transport to evolve the basis functions spanning the reduced space. Analogously, \citet{GERBEAU2014246} evolve the reduced basis based on approximated Lax pairs. Related concepts are followed by  \citet{doi:10.1137/20M1316998} for transporting a reduced space along characteristics. A different approach is presented in \citet{BlaSU20}, which evolves the modes of the reduced basis via instantaneous residual minimization in time; see \Cref{sec:Template} for a detailed discussion.

\citet{RAPUN20103046} update a reduced model by switching back to the full model occasionally to generate new snapshots. The method introduced in  \citet{doi:10.1002/nme.4800,ETTER2020112931} aims to mimic adaptive mesh refinement in finite elements by splitting reduced-basis vectors to enrich the reduced space. The focus is on failsafe refinement in the sense that the full-model space is recovered in the limit. 

Another line of work performs data-driven online basis updates by incorporating newly available online information, such as full-model residuals and incomplete state fields, to derive low-rank updates to the reduced spaces. For example, for linear problems, dynamic data-driven reduced models introduced by  \citet{PEHERSTORFER201521} adapt reduced spaces with an incremental SVD approach and the reduced-model operators with corresponding low-rank updates, without having to rebuild the reduced model from scratch. 

For nonlinear problems, \citet{Peherstorfer15aDEIM} introduce the online adaptive discrete empirical interpolation method (ADEIM), which builds on the empirical interpolation method \cite{Barrault2004,Chaturantabut2010} and computes low-rank updates to the reduced space from sparse residual samples. The basis update is formulated as a low-rank correction inferred from a small (possibly random) sketch via a sampling operator. Relying on sparse residual samples helps to reduce online costs. The sparsity of the updates also helps to mitigate the lifting bottleneck when reducing nonlinear systems (see \Cref{sec:Prelim:LiftingBottleneckSub}).
A central component of ADEIM is how to choose and adapt the sampling points over time. \citet{P18AADEIM} couple online basis updates with adaptive sampling, which is shown to be critical for transport-dominated problems.  
Sparse sampling raises the question not only of where to sample but also when to sample: sampling the residual information corresponding to the current state can lag behind. Lookahead strategies use short-term forecasting rules to evaluate the residual and forecasted states that anticipate upcoming dynamics, which improves robustness of the online updates and reduced instabilities; see \citet{SUP23ADEIMAHEAD,https://doi.org/10.1002/fld.5240}. 

There are multiple ways to compute the low-rank updates from the sparse samples. The basis updates proposed in the original ADEIM work \cite{Peherstorfer15aDEIM} do not guarantee that the adapted basis is orthogonal. \citet{ZPW17SIMAXManifold} build on subspace tracking methods from signal processing  \cite{BalzanoNowakRecht2010} to derive an efficient rank-one update rule for ADEIM that preserves orthonormality of the adapted basis over time. 
\citet{HUANG2023112356,MOHAGHEGH2026114468} introduce a related online adaptation strategy specifically for least-squares formulations, which are key for reducing chaotic and multiscale problems.

There are also time-dependent basis methods in model reduction that do not update the basis via online sparse residual information. For example,  \citet{RAMEZANIAN2021113882,NADERI2023115813,AITZHAN2025113549} derive and integrate evolution equations for a low-rank representation directly from the governing equations. In this sense, these methods share similarities with dynamic low-rank approximations. 

There are also other online adaptive reduced basis methods inspired by dynamic low-rank approximations, such as \citet{HesthavenPagliantini2021_StructurePreservingPoisson,Pagliantini2021,HesthavenM2AN2022} that are structure-preserving. In particular, the works \cite{doi:10.1137/23M1601225,pagliantini2025adaptivehyperreductionnonsparseoperators} show how to preserve structure even when the basis is adapted via  sparse sampling-type reduction strategies such as empirical interpolation. 

Another approach to online adaptive reduced space is to precompute a time-dependent reduced space in the offline phase from training trajectories and then use it online without further adaptation or enrichment of the reduced space  \cite{doi:10.1137/16M1071493}.

\subsection{Adaptive discrete empirical interpolation method (ADEIM)}
We now discuss online adaptive model reduction using ADEIM as an example of how to realize online updates. Reduced models based on ADEIM adapt the reduced space online with low-rank updates from sparse residual information. 
The presentation loosely follows \citet{Peherstorfer15aDEIM} and \citet{P18AADEIM}. We drop the dependence of the state on $\bfmu$ for ease of exposition in this section.

\subsubsection{Time-discrete full models}
Let us consider a fully discrete full model of the form
\begin{equation}\label{eq:ADEIM:DiscFOM}
\bfq_{k-1} = \bff^{(\Delta)}(\bfq_k; \bfmu)\,,\qquad k = 1, \dots, K\,,
\end{equation}
which can be obtained from the semi-discrete full model \eqref{eqn:FOM} with, e.g., an implicit time-integration scheme. In particular, $\bfq_k$ is an approximation to $\bfq(t_k)$ for some $t_k\in[0,T]$. The form~\eqref{eq:ADEIM:DiscFOM} requires generally a nonlinear solve to obtain~$\bfq_k$ from~$\bfq_{k-1}$. 
Note that the right-hand side function $\bff$ of the semi-discrete full model~$\eqref{eqn:FOM}$ and the right-hand side function~$\bff^{(\Delta)}$ of the fully discrete full model~\eqref{eq:ADEIM:DiscFOM} are different. 

Following the steps in \Cref{sec:LinMOR:OfflineStep}, a discrete reduced model can be obtained from the discrete full model as
\begin{equation}\label{eq:ADEIM:StaticROM}
\hat{\bfq}_{k-1} = \hat{\bff}^{(\Delta)}(\hat{\bfq}_k; \bfmu)\,,\qquad k = 1, \dots, K\,,
\end{equation}
where the right-hand side function $\hat{\bff}^{(\Delta)}$ is obtained from the full-model right-hand side function $\bff^{(\Delta)}$ via empirical interpolation; see \Cref{sec:Prelim:EIM}. In particular, the map $\hat{\bff}^{(\Delta)}$ depends on the basis matrix $\bfV = [\bfv_1, \dots, \bfv_n] \in \mathbb{R}^{N \times n}$ and the interpolation points matrix $\bfP$ corresponding to the pairwise distinct interpolation points $p_1, \dots, p_n \in \{1, \dots, N\}$. 

In this section on ADEIM and in the context of online adaptive reduced modeling, we refer to the reduced model \eqref{eq:ADEIM:StaticROM} as a static reduced model because the basis matrix $\bfV$ and the interpolation points $\bfP$ are fixed over all time steps $k = 1, \dots, K$.

\subsubsection{ADEIM: Basis and points adaptation}\label{sec:ADEIM:BasisAndPointsAdaptation}
We now discuss ADEIM's update mechanism for the basis and interpolation points. 

\paragraph{Basis updates and nonlinear parametrizations}
Let now the basis matrix $\bfV_k$ depend on the time step $k = 1, \dots, K$. ADEIM adapts the basis matrix $\bfV_{k-1}$ to $\bfV_k$ as
\begin{equation}\label{eq:ADEIM:RankRUpdateToVk}
\bfV_{k} = \bfV_{k-1} + \bfalpha_{k-1}\bfbeta_{k-1}^{\top}\,,
\end{equation}
where $\bfalpha_{k-1} \in \mathbb{R}^{N \times r}$ and $\bfbeta_{k-1} \in \mathbb{R}^{n \times r}$ provide a rank $r \in \N$ update. Typically, the rank is small,  $r \ll n$, or even just a rank-one update, $r = 1$. In addition to the basis matrix, ADEIM adapts the interpolation points so that the matrix $\bfP_{k-1}$ is adapted to $\bfP_k$. 

Adapting the basis matrix $\bfV_k$ (and the interpolation points matrix $\bfP_k$) over time can be interpreted as evolving the representation of the reduced state. In the terminology of nonlinear parametrizations, this means that the online weights at time step $k$ comprise the reduced state (e.g., $\hat{\bfq}_k$) as well as the low-rank update variables (e.g., $\bfalpha_{k-1}$ and $\bfbeta_{k-1}$). More generally, the online weights may also include variables associated with the interpolation-point update, as well as additional auxiliary variables required by the update, such as sampling points and window data. The offline weights consist of quantities fixed in the offline phase, such as the initial reduced basis matrix $\bfV_0$ and initial interpolation points matrix~$\bfP_0$. From this viewpoint, ADEIM can be interpreted as realizing a fully discrete instance of the time-adaptive nonlinear parametrization discussed earlier in \eqref{eq:NonlinParam:AdaptiveInTime}: the reduced solution is still represented in a low-dimensional form, but the representation is adapted online through time-dependent basis updates.

\paragraph{Basis updates}
Let us first consider the basis adaptation. For ease of exposition, we adapt the basis matrix at every time step $k = 1, \dots, K$; however, in practice, the basis is often updated only at selected adaptation steps (e.g., every $z$-th time step), to reduce online costs. At time step $k$, the basis matrix $\bfV_{k-1}$ from the previous time step $k-1$ is given, and ADEIM adapts it to the basis matrix $\bfV_k$.

The update is constructed such that the space given by the updated basis matrix~$\bfV_k$ can approximate well the vectors in a sliding window 
\begin{equation}\label{eq:ADEIM:SlidingWindow}
\bfF_{k-1} = [\mathfrak{f}_{k-w}, \dots, \mathfrak{f}_{k-1}] 
\end{equation}
of length $w \in \N$. 
One might want to use the full-model solutions as the vectors of the sliding window $\bfF_{k-1}$, but they are unavailable. Instead, ADEIM uses surrogates $\mathfrak{f}_{k-w}, \dots, \mathfrak{f}_{k-1}$ that capture some of the information from the full-model solutions but can be computed without the expensive cost of time-stepping the full model. We will discuss how to obtain the vectors of the sliding window in \Cref{sec:ADEIM:SlidingWindow}. 
For now, let us assume the sliding window $\bfF_{k-1}$ is given. 

Let us formally define an objective $J_{\text{nl}}$ for computing the basis update as the norm of the residual when approximating the columns of $\bfF_{k-1}$ in the updated space spanned by the adapted basis matrix $\bfV_{k-1} + \bfalpha\bfbeta^{\top}$ and with the interpolation points matrix $\bfP$, i.e.,
\begin{equation}\label{eq:ADEIM:IdealJ}
J_{\text{nl}}(\bfalpha, \bfbeta, \bfP) = \|(\bfV_{k-1}  + \bfalpha\bfbeta^{\top})(\bfP^{\top}(\bfV_{k-1}  + \bfalpha\bfbeta^{\top}))^{-1}\bfP^{\top}\bfF_{k-1} - \bfF_{k-1}\|_F\,.
\end{equation}
A low-rank update and interpolation points matrix that minimize $J_{\text{nl}}$ minimize the error of approximating the vectors in the sliding window in the adapted space~\eqref{eq:ADEIM:RankRUpdateToVk}. Note that an admissible interpolation points matrix needs to satisfy additional conditions, which are not present in the objective \eqref{eq:ADEIM:IdealJ}. 
The objective $J_{\text{nl}}$ is nonlinear in the update $\bfalpha\bfbeta^{\top}$ and the interpolation points matrix $\bfP$. Solving a nonlinear optimization problem with objective $J_{\text{nl}}$ can be computationally expensive in the online phase. 

ADEIM simplifies the objective $J_{\text{nl}}$ to the objective 
\begin{equation}\label{eq:ADEIM:Update}
J(\bfalpha, \bfbeta) =  \|\left(\bfV_{k-1} + \bfalpha\bfbeta^{\top}\right)\bfC_{k-1} - \bfF_{k-1}\|_F\,,
\end{equation}
where the coefficient matrix from time step $k - 1$, \begin{equation}\label{eq:ADEIM:Ck-1}
\bfC_{k - 1} = (\bfP^{\top}_{k-1}\bfV_{k-1})^{-1}\bfP^{\top}_{k-1}\bfF_{k-1}\,,
\end{equation}
is used instead of the coefficients that depend on the adapted quantities as in \eqref{eq:ADEIM:IdealJ}. Correspondingly, the objective $J$ is independent of the adapted interpolation points matrix.
To further reduce the costs of the online update, the residual is only minimized at sparse sampling points, 
\begin{equation}\label{eq:ADEIMSparseUpdate}
\min_{\substack{\bfalpha \in \mathbb{R}^{N \times r}\\ \bfbeta \in \mathbb{R}^{n \times r}}} \|\bfS_{k-1}^{\top}\left(\left(\bfV_{k-1} + \bfalpha\bfbeta^{\top}\right)\bfC_{k-1} - \bfF_{k-1}\right)\|_F
\end{equation}
where $\bfS_{k-1} = [\bfe_{s_{1}^{(k-1)}}, \dots, \bfe_{s_{n_s}^{(k-1)}}]$ is the sampling matrix corresponding to the pairwise distinct sampling points $s_1^{(k-1)}, \dots, s_{n_s}^{(k-1)}$. Typically, $n \leq n_s \ll N$. Notice that the problem \eqref{eq:ADEIMSparseUpdate} specifies $\bfalpha_{k - 1}$ only at the $n_s$ sampling points $s_1^{(k-1)}, \dots, s_{n_s}^{(k-1)}$. 

The work \citet{Peherstorfer15aDEIM} builds on an eigenvalue decomposition to compute the update that solves \eqref{eq:ADEIMSparseUpdate}; however, again with the aim of reducing the costs of the online update, the work \cite[Algorithm~3]{P18AADEIM} shows that an approximation can be efficiently computed from a thin SVD of the $n_s \times w$ matrix $\bfS_{k-1}^{\top}\bfR_{k-1}$, which is the residual matrix $\bfR_{k-1} = \bfV_{k-1}\bfC_{k-1} - \bfF_{k-1}$ restricted to the components corresponding to the sampling points $s_1^{(k-1)}, \dots, s_{n_s}^{(k-1)}$: Let $\bfU_{1:r} \in \mathbb{R}^{n_s \times r}$ and $\bfW_{1:r} \in \mathbb{R}^{w \times r}$ contain the first $r$ left- and right-singular vectors of $\bfS_{k - 1}^{\top}\bfR_{k-1}$ corresponding to singular values $\sigma_1, \dots, \sigma_r$, respectively. Collect the first $r$ singular values $\sigma_1, \dots, \sigma_r$ on the diagonal of $\bfSigma_{1:r}$, then the ADEIM update is 
\begin{equation}\label{eq:ADEIM:BasisUpdateAlphaBeta}
\bfalpha_{k-1} = - \bfS_{k-1}\bfU_{1:r}\bfSigma_{1:r} \in \mathbb{R}^{N \times r}\,,\qquad \bfbeta_{k-1} = (\bfC_{k-1}^+)^{\top}\bfW_{1:r} \in \mathbb{R}^{n \times r}\,,
\end{equation}
which changes the basis matrix $\bfV_{k - 1}$ only at the sampling points $s_1^{(k-1)}, \dots, s_{n_s}^{(k-1)}$ when applied as in \eqref{eq:ADEIM:RankRUpdateToVk}. This is consistent with the objective \eqref{eq:ADEIMSparseUpdate} that defines the update $\bfalpha_{k - 1}$ only at the sampling points $s_1^{(k-1)}, \dots, s_{n_s}^{(k-1)}$ corresponding to $\bfS_{k-1}$. Notice, however, that subsequent orthonormalization of $\bfV_k$ generally alters all rows. The costs of computing the SVD for solving the problem scale as $\mathcal{O}(n_s w^2)$, in the common case that $w \leq n_s$.

\paragraph{Updating the interpolation points}
After the basis matrix has been updated from $\bfV_{k-1}$ to $\bfV_k$, the interpolation points matrix $\bfP_{k-1}$ is updated to $\bfP_k$.  \citet{Peherstorfer15aDEIM} propose an algorithm that computes a rank-one update to $\bfP_{k-1}$, which corresponds to an update of at most one point per basis update. In contrast, \citet{P18AADEIM} re-computes the interpolation points using the QDEIM algorithm \cite{drmac-gugercin-DEIM-2016} from scratch. The costs of computing the points from scratch scale as $\mathcal{O}(Nn^2)$ and thus linearly in the full-model dimension~$N$.

\subsubsection{ADEIM: Forward step of reduced model}\label{sec:ADEIM:ForwardStep}
At time step $k$, after the basis matrix has been updated to $\bfV_k$ and the interpolation points matrix to $\bfP_k$, the corresponding reduced model right-hand side function~$\hat{\bff}_k^{(\Delta)}$ is derived from the full-model right-hand side function $\bff^{(\Delta)}$ via empirical interpolation. 
Notice that the reduced right-hand side function $\hat{\bff}_k^{(\Delta)}$ depends on $k$ because the basis matrix and the interpolation points depend on $k$.
The reduced model is then integrated forward one step from $k-1$ to $k$ to compute the reduced solution $\hat{\bfq}_k$ at time step $k$, which means solving the generally nonlinear system of equations $
\hat{\hat{\bfq}}_{k-1} = \hat{\bff}_k^{(\Delta)}(\hat{\bfq}_k; \bfmu)$, 
where $\hat{\hat{\bfq}}_{k-1}$ is $\bfV_{k-1}\hat{\bfq}_{k-1}$ projected onto the adapted basis $\bfV_k$.

\subsubsection{ADEIM: Sliding window}\label{sec:ADEIM:SlidingWindow}
Recall that the columns of the sliding window  (see  \eqref{eq:ADEIM:SlidingWindow}) serve as surrogates for unavailable full-model solutions. ADEIM leverages that if the full-model solution $\bfq_{k}$ at time step $k$ is plugged into the fully discrete full model \eqref{eq:ADEIM:DiscFOM}, then $\bfq_{k-1}$ is obtained.  

In ADEIM, at time step $k$, the sliding window is updated from $\bfF_{k-1}$ to $\bfF_k$. After computing the reduced solution $\hat{\bfq}_k$ with the adapted basis~$\bfV_k$ (see \Cref{sec:ADEIM:ForwardStep}), it is used to derive
\begin{equation}\label{eq:ADEIM:SurrogateF}
\mathfrak{f}_{k} = \bff^{(\Delta)}(\bfV_{k}\hat{\bfq}_{k}; \bfmu)\,,
\end{equation}
which is an approximation of the full-model solution $\bfq_{k-1}$ at the previous time step $k - 1$. The window is then updated from $\bfF_{k-1}$ to $\bfF_k$ by dropping the oldest column and adding $\mathfrak{f}_{k}$. The window $\bfF_k$ is then used to drive the basis update from $\bfV_k$ to $\bfV_{k+1}$ at the next time step $k + 1$. 

The surrogate $\mathfrak{f}_{k}$ generally has components outside the reduced space spanned by the columns of $\bfV_{k}$, whereas the lifted reduced solution $\bfV_{k}\hat{\bfq}_{k}$ lies in the range of~$\bfV_{k}$ by construction. Therefore, using surrogates of the form \eqref{eq:ADEIM:SurrogateF} in the sliding window can reveal missing directions in the current basis and lead to meaningful basis updates \eqref{eq:ADEIM:BasisUpdateAlphaBeta}. In this sense, $\mathfrak{f}_{k}$ can be viewed as residual-type information: it probes the full-model right-hand side at the current reduced solution and captures components that the current basis cannot represent.

Now recall that ADEIM minimizes the residual \eqref{eq:ADEIMSparseUpdate} only at the $n_s < N$ sampling points corresponding to $\bfS_{k-1}$ and not at all $N$ components.
This means that the full-model right-hand side function needs to be evaluated only at the components corresponding to the sampling points, which we denote as 
\begin{equation}\label{eq:ADEIM:SurrogateTildeQ1}
\bfS_{k-1}^{\top}\mathfrak{f}_{k} = \bfS_{k-1}^{\top}\bff^{(\Delta)}(\bfV_{k}\hat{\bfq}_{k}; \bfmu)\,.
\end{equation}
The components that are not indexed by $\bfS_{k-1}$ can be approximated with empirical regression as 
\begin{equation}\label{eq:ADEIM:SurrogateTildeQ2}
\breve{\bfS}_{k-1}^{\top} \mathfrak{f}_{k} = \breve{\bfS}_{k-1}^{\top}\bfV_{k}(\bfS_{k-1}^{\top}\bfV_{k})^{+}\bfS_{k-1}^{\top}\bff^{(\Delta)}(\bfV_{k}\hat{\bfq}_{k}; \bfmu)\,,
\end{equation}
where $\breve{\bfS}_{k-1}$ corresponds to the sampling points matrix associated with the complement $\{1, \dots, N\} \setminus \{s_1^{(k-1)}, \dots, s_{n_s}^{(k-1)}\}$. Note we approximate the unobserved components by a least-squares regression using the pseudo inverse $(\bfS_{k-1}^{\top}\bfV_{k})^{+}$ rather than via interpolation; see \Cref{sec:EmpiricalRegression} for oversampling in empirical interpolation. With this procedure, ADEIM evaluates the right-hand side function $\bff^{(\Delta)}$ of the full model at only the $n_s$ sampling points. 

Computing $\mathfrak{f}_{k}$ requires only evaluating the time-discrete full-model right-hand side function $\bff^{(\Delta)}$ at the lifted reduced solution $\bfV_{k}\hat{\bfq}_{k}$, rather than solving the nonlinear systems corresponding to integrating the full model \eqref{eq:ADEIM:DiscFOM} for one time step.  Consequently, the cost of obtaining a vector for the sliding window is comparable to a full-model residual evaluation. Note, however, that if the full model is discretized with an explicit scheme, then evaluating $\bff^{(\Delta)}$ and advancing the full model by one step are typically of comparable cost.

The procedure described here for filling the sliding window is backward-looking: the reduced solution $\hat{\bfq}_{k}$ at time step $k$ is used to obtain the surrogate approximating the full-model solution at time step $k-1$. Lookahead strategies have been developed that use inexpensive forecasts to generate forward-looking surrogates, which avoid the one-step lag and so can improve stability and robustness \cite{SUP23ADEIMAHEAD,https://doi.org/10.1002/fld.5240}.

\subsubsection{ADEIM: Sampling}\label{sec:ADEIM:Sampling}
The ADEIM basis updates are driven by the sliding window, which is obtained from sparse evaluations of the full-model right-hand side function as given in \eqref{eq:ADEIM:SurrogateTildeQ1}--\eqref{eq:ADEIM:SurrogateTildeQ2}. 
The indices of the components at which the full-model right-hand side function is evaluated are given by the sampling points denoted as $s_1^{(k-1)}, \dots, s_{n_s}^{(k-1)}$, which are collected into the sampling points matrix $\bfS_{k-1}$. 
\citet{Peherstorfer15aDEIM} primarily focus on uniformly distributed sampling points. However, for transport-dominated problems, which can exhibit local features such as wave fronts, uniform sampling is inefficient, as shown in \cite{P18AADEIM}.
Instead, \citet{P18AADEIM} proposes an adaptive sampling strategy. 

Let us measure the distance between two spaces spanned by the basis matrices $\bar{\bfV}$ and $\bfV_k$, respectively, via
\begin{equation}\label{eq:ADEIM:MeasureDistance}
d(\bar{\bfV}, \bfV_k) =  \|\bar{\bfV} - \bfV_k\bfV_k^{\top}\bar{\bfV}\|_F^2\,,
\end{equation}
where it is assumed that $\bfV_k$ and $\bar{\bfV}$ are orthonormal for ease of exposition. 
Let us consider the situation that we want to adapt the basis matrix $\bfV_{k}$ to $\bfV_{k+1}$ so that $\bfV_{k+1}$ spans a space that is close to the space spanned by $\bar{\bfV}$ in the sense of \eqref{eq:ADEIM:MeasureDistance}.  
Notice that we now consider the adaptation from time step $k$ to $k + 1$ instead of $k - 1$ to~$k$. The reason is that the sampling points adaptation is performed in the algorithm after the basis has been updated at time step $k$, and hence is in preparation for adapting the space at the subsequent time step $k + 1$; see the algorithm given in \Cref{sec:ADEIM:Algorithm}.

When adapting from $\bfV_k$ to $\bfV_{k + 1}$,  the sliding window is $\bfF_{k}$. We assume for now that it is given as $\bfF_{k} = \bar{\bfV}\bar{\bfF}_{k}$ with $\bar{\bfF}_{k}$ full rank and $w \geq n$. This means that the column span of the sliding window $\bfF_{k}$ is the same space as the space spanned by~$\bar{\bfV}$. 
Consider the coefficient matrix $\bfC_{k}$, analogously defined to $\bfC_{k-1}$ in \eqref{eq:ADEIM:Ck-1}, and the residual matrix $\bfR_{k} = \bfV_{k}\bfC_{k} - \bfF_{k}$. 
Let now $\bfS$ be a valid sampling points matrix. \citet{P18AADEIM} shows that the ADEIM rank-$r$ update $\bfV_{k+1} = \bfV_{k} + \bfalpha_{k}\bfbeta_{k}^{\top}$ via~\eqref{eq:ADEIMSparseUpdate} with $\bfF_k$ and sampling points matrix $\bfS$ leads to the bound
\begin{equation}\label{eq:ADEIM:SamplingErrorBound}
d(\bar{\bfV}, \bfV_{k+1}) \leq \frac{1}{\sigma^2_{\text{min}}(\bfF_{k})} \left(\|\breve{\bfS}^{\top}\bfR_k\|_F^2 + \sum_{i = r + 1}^{\bar{r}}\sigma_i^2\right)\,,
\end{equation}
where $\bar{r}$ is the rank of $\bfS^\top\bfR_{k}$, $\sigma_1 \geq \dots \geq \sigma_{\bar{r}} > 0$ are the singular values of $\bfS^\top\bfR_{k}$, and $\sigma_{\text{min}}(\bfF_{k})$ is the smallest non-zero singular value of $\bfF_{k}$. Recall that $\breve{\bfS}$ is the complementary sampling points matrix to $\bfS$ in the sense that it corresponds to the sampling points $\{1, \dots, N\} \setminus \{s_1, \dots, s_{n_s}\}$. 
The bound \eqref{eq:ADEIM:SamplingErrorBound} motivates introducing the objective
\begin{equation}\label{eq:ADEIM:SamplingObjectiveB}
b(\bfS) = \|\breve{\bfS}^{\top}\bfR_k\|_F^2 + \sum_{i = r + 1}^{\bar{r}}\sigma_i^2\,.
\end{equation}

Selecting sampling points $\{s_1, \dots, s_{n_s}\} \subseteq \{1, \dots, N\}$ so that the corresponding matrix~$\bfS$ minimizes the objective $b$ incurs a combinatorial subset-selection problem in the full-model dimension $N$, which makes this direct approach intractable.
Instead, \citet{P18AADEIM} suggest to select sampling points that minimize only the norm $\|\breve{\bfS}^{\top}\bfR_k\|_F^2$ of the sparse residual $\breve{\bfS}^{\top}\bfR_k$ corresponding to the complementary sampling points $\{1, \dots, N\}\setminus \{s_1, \dots, s_{n_s}\}$, which is one term in the objective $b$.
Because the equality
\[
\|\bfR_k\|_F^2 = \|\bfS^{\top}\bfR_k\|_F^2 + \|\breve{\bfS}^{\top}\bfR_k\|_F^2
\]
holds for any sampling points matrix $\bfS$ and its complementary sampling points matrix~$\breve{\bfS}$, finding a set of sampling points so that its complement minimizes~$\|\breve{\bfS}^{\top}\bfR_k\|_F^2$ is equivalent to finding a set of sampling points that maximizes $\|\bfS^{\top}\bfR_k\|_F^2$.
This can be achieved as follows: First, compute the Euclidean norm of all $N$ rows of $\bfR_k$, then sort them in descending order, and then select the $n_s$ indices corresponding to the $n_s$ largest row norms. 
The selected indices are the new sampling points $s_1^{(k)}, \dots, s_{n_s}^{(k)}$, with sampling points matrix $\bfS_k$. 

\citet{CKMP19ADEIMQuasiOptimalPoints} show that the approach based on maximizing the residual norm leads to an objective value \eqref{eq:ADEIM:SamplingObjectiveB} that is at most a factor two worse than the objective value obtained with the optimal sampling points selection that minimizes~$b$. 
The costs of selecting the points via the residual norm are dominated by the costs of evaluating the full-model right-hand side function at all $N$ components and then sorting the norms of the $N$ rows.
There is a range of other sampling strategies for ADEIM and related online adaptive model reduction techniques \cite{doi:10.2514/1.J062869,https://doi.org/10.1002/fld.5240,MOHAGHEGH2026114468}. For example, \citet{doi:10.2514/1.J062869} focus on computational fluid-dynamics applications and proposes physics-motivated strategies. 

\subsubsection{Online algorithm for ADEIM reduced models}\label{sec:ADEIM:Algorithm}
There are several variants of ADEIM in the literature \cite{Peherstorfer15aDEIM,P18AADEIM,SUP23ADEIMAHEAD}, but they share the same core steps: 
\begin{enumerate}
    \item[(1)] low-rank basis adaptation (\Cref{sec:ADEIM:BasisAndPointsAdaptation}),
    \item[(2)] interpolation-points update (\Cref{sec:ADEIM:BasisAndPointsAdaptation}),
    \item[(3)] reduced time integration with empirical interpolation (\Cref{sec:ADEIM:ForwardStep}),
    \item[(4)] sparse full-model queries to populate the sliding window (\Cref{sec:ADEIM:SlidingWindow}), and
    \item[(5)] sampling-points adaptation (\Cref{sec:ADEIM:Sampling}).
\end{enumerate}
We summarize a generic online ADEIM procedure below. For ease of exposition, we assume that the basis is adapted at every time step; in practice, steps (5) and/or (1)--(2) can be executed only intermittently (e.g., every $z$-th step) to reduce online costs.

We follow the online ADEIM procedure given in 
\citet[Algorithm 1]{P18AADEIM}. For the first $w_{\text{init}} \in \mathbb{N}$ time steps, the full model is integrated in time to generate the snapshots $\bfQ_{1:w_{\text{init}}} \in \mathbb{R}^{N \times w_{\text{init}}}$, from which the initial basis matrix $\bfV_{w_{\text{init}}}$, interpolation points matrix $\bfP_{w_{\text{init}}}$, and reduced right-hand side function $\hat{\bff}^{(\Delta)}_{w_{\text{init}}}$ are constructed with empirical interpolation. Furthermore, the sliding window $\bfF_{w_{\text{init}}} = \bfQ_{w_{\text{init}}-w + 1:w_{\text{init}}}$ is initialized with the snapshots. 
Then, the following steps are iterated for $k = w_{\text{init}} + 1, \dots, K$: 
\begin{itemize}
\item \emph{Step 1 (basis update)}: Adapt ADEIM basis from $\bfV_{k-1}$ to $\bfV_{k}$ using the update~\eqref{eq:ADEIM:BasisUpdateAlphaBeta} based on the sliding window $\bfF_{k-1}$ and sampling points matrix~$\bfS_{k-1}$~(\Cref{sec:ADEIM:BasisAndPointsAdaptation}).
\item \emph{Step 2 (interpolation points update)}: Compute the interpolation points matrix~$\bfP_{k}$ by, e.g., applying QDEIM to the adapted basis matrix $\bfV_{k}$ (\Cref{sec:ADEIM:BasisAndPointsAdaptation}). 
\item \emph{Step 3 (forward step)}: Solve the reduced model given by the right-hand side function $\hat{\bff}^{(\Delta)}_k$ with adapted basis and interpolation points 
for $\hat{\bfq}_k$ and store the lifted reduced solution $\bfV_k\hat{\bfq}_k$ by extending $\bfQ_{1:k-1}$ to $\bfQ_{1:k}$. Note that $\hat{\bff}^{(\Delta)}_k$ now depends on the step $k$ because the underlying basis $\bfV_k$ and interpolation points matrix $\bfP_k$ depend on $k$  (\Cref{sec:ADEIM:ForwardStep}).
\item \emph{Step 4 (sliding window update) and optional Step 5 (sampling points update)}: If the sampling points should be adapted in this time step $k$, then evaluate the full-model right-hand side function $\bff^{(\Delta)}(\bfV_k\hat{\bfq}_k)$ at all $N$ components and append it to the sliding window to obtain $\bfF_k$. Drop the oldest column from $\bfF_k$. Compute the coefficient matrix $\bfC_k$ and the residual matrix $\bfR_k = \bfV_k\bfC_k - \bfF_k$ and select new sampling points $\bfS_k$ based on row-wise residual norms (\Cref{sec:ADEIM:Sampling}). 
Otherwise, if the sampling points should not be updated in this time step $k$, then evaluate $\bff^{(\Delta)}(\bfV_k\hat{\bfq}_k)$ only at the sampling points corresponding to $\bfS_{k-1}$ and update the sliding window $\bfF_k$ with \eqref{eq:ADEIM:SurrogateTildeQ1}--\eqref{eq:ADEIM:SurrogateTildeQ2}. Then set $\bfS_k = \bfS_{k-1}$.
\end{itemize}

The computational cost of an online time step with an ADEIM reduced model is typically dominated by evaluating the full-model right-hand side function used to update the sliding window. At time steps where the sampling points are adapted, the full-model right-hand side function is evaluated at all $N$ components, which results in online costs that scale at least linearly in $N$, but this is typically still cheaper than solving the nonlinear system \eqref{eq:ADEIM:DiscFOM} corresponding to the full model at this time step. In all other time steps, the full-model right-hand side function is evaluated only at $n_s \ll N$ sampling points, so the costs of these sparse full-model queries scale with $n_s$ rather than $N$.

\section{Instantaneous residual minimization methods}\label{sec:InstantaneousResMin}
We now consider generic parametrizations \eqref{eq:NMOR:Param} with a time- and parameter-dependent weight vector $\bftheta(t, \bfmu)$ that enters nonlinearly, possibly in combination with an offline weight vector as in \eqref{eq:NonLin:ParamWithOfflineWeight}. Generic means here that we do not restrict the parametrization to have a specific structural form for how the weight vector $\bftheta(t, \bfmu)$ enters. This stands in contrast to the transformation-based methods discussed in \Cref{sec:Template} and the online adaptive basis methods in \Cref{sec:ADEIM}, where the weight vector enters in a specific, prescribed way that is leveraged algorithmically for efficient computation. To derive reduced dynamics for such generic nonlinear parametrizations, we review instantaneous residual minimization methods. The central idea is to determine, at each time $t$, the time derivative $\dot{\bftheta}(t, \bfmu)$ of the weight vector by minimizing the residual norm. 

\subsection{Overview of this section}
We present the Dirac--Frenkel variational principle as a core building block for instantaneous residual minimization with generic nonlinear parametrizations in \Cref{sec:DF:DF}, and review numerical methods such as Neural Galerkin schemes that build upon this principle. The derivation and practical implementation of methods based on the Dirac--Frenkel variational principle involve several numerical considerations. First, the residual norm has to be approximated numerically. We review static and adaptive collocation approaches in the context of Neural Galerkin schemes in \Cref{sec:NeuralGalerkin}. Second, the evolution equations for the weight vector induced by the residual minimization need to be discretized in time. We briefly consider time-stepping choices and their implications in \Cref{sec:DiscreteDF}.
Third, instantaneous residual minimization may be insufficient to fully determine the weight vector, leading to singular or poorly conditioned systems. We therefore discuss regularization strategies for such cases in \Cref{sec:RegDF}.
Fourth, the residual minimization leads to a system of algebraic equations given by the first-order optimality conditions, which must be solved numerically at each time step. We outline solvers based on randomized sketching to reduce the costs, and we discuss how randomization can also improve conditioning in \Cref{sec:RegDF}. 

We further relate the Dirac--Frenkel variational principle to optimize-then-discretize schemes and discuss discretize-then-optimize schemes for instantaneous residual minimization as an alternative to the Dirac--Frenkel variational principle in \Cref{sec:DtOOtD}. In \Cref{sec:DFWithPretrained} we consider instantaneous residual minimization in the context of nonlinear parametrizations that have been trained on snapshot data, including autoencoder-based approaches and neural representations. 

\subsection{The Dirac--Frenkel variational principle}\label{sec:DF:DF}
A common way to formulate instantaneous residual minimization is via the Dirac--Frenkel variational principle. Its origin can be traced back to \citet{dirac1930note} and \citet{frenkel1934wave}. More modern presentations can be found in \citet{McLachlan01011964,Kramer1981,BROECKHOVE1988547,lubich2008quantum}. A short history of the Dirac--Frenkel variational principle is reported in \citet[Sec.~3.8]{lasser_lubich_2020}. 

\subsubsection{The Dirac--Frenkel variational principle}
A nonlinear parametrization $\hat{q}(\bftheta, \cdot)\colon \Omega \to \mathbb{R}$ with weight vector $\bftheta \in \mathbb{R}^n$ (which could depend on time $t$ and parameter $\bfmu$) induces the set $\widehat{\mathcal{M}}$ defined in \eqref{eq:IRM:HatM}, 
which contains all functions that can be represented by varying the weight vector in $\mathbb{R}^n$. Note that we do not explicitly denote the possible dependence of $\hat{q}$ on an offline weight vector for now.
Recall that $\mathcal{M}$ denotes the solution manifold defined in \eqref{eq:Prelim:SolManifold}, whereas we denote the set \eqref{eq:IRM:HatM} as $\widehat{\mathcal{M}}$. Even though we do not restrict the parametrizations to have a specific structural form, we do need to make suitable regularity assumptions, e.g.,  $\hat{q}$ is differentiable in the weight vector $\bftheta$ and  all functions $\hat{q}(\bftheta, \cdot) \in \widehat{\mathcal{M}}$ are also in the space $\mathcal{Q}$ over which the PDE problem is formulated; see \Cref{sec:FullModel}. 

Let us now interpret $\widehat{\mathcal{M}}$ as a trial manifold for solving the PDE problem \eqref{eqn:PDE} and plug the nonlinear parametrization $\hat{q}(\bftheta(t, \bfmu), \cdot)$ with time- and parameter-dependent weight vector $\bftheta(t, \bfmu)$ into the PDE \eqref{eqn:PDE} to obtain the residual function $r\colon \mathbb{R}^n \times \mathbb{R}^n \times \Omega \to \mathbb{R}$,
\begin{equation}\label{eq:DF:ResidualFun}
r(\bftheta(t, \bfmu), \dot{\bftheta}(t, \bfmu), \cdot) = \nabla_{\bftheta}\hat{q}(\bftheta(t, \bfmu), \cdot)^{\top}\dot{\bftheta}(t, \bfmu) - f(\cdot, \hat{q}(\bftheta(t, \bfmu), \cdot); \bfmu)\,.
\end{equation}
We used the chain rule 
\[
\partial_t \hat{q}(\bftheta(t, \bfmu), \cdot) = \sum_{i = 1}^n \partial_{\theta_i}\hat{q}(\bftheta(t, \bfmu), \cdot)\dot{\theta}_i(t, \bfmu) = \nabla_{\bftheta} \hat{q}(\bftheta(t, \bfmu), \cdot)^{\top}\dot{\bftheta}(t, \bfmu)
\]
to make the dependence of the time derivative $\partial_t \hat{q}(\bftheta(t, \bfmu), \cdot)$ on the time derivative of the weight vector $\dot{\bftheta}(t, \bfmu)$ more explicit. 

The Dirac--Frenkel variational principle seeks a time derivative $\dot{\bftheta}(t, \bfmu) \in \mathbb{R}^n$ of the weight vector such that the residual is orthogonal to all admissible variations of $\hat{q}(\bftheta(t, \bfmu), \cdot)$ at time $t$. If $\widehat{\calM}$ has a suitable manifold structure induced by the regularity of the parametrization, then this is equivalent to
\begin{equation}\label{eq:DF:Conditions}
    \left\langle v, r(\bftheta(t, \bfmu), \dot{\bftheta}(t, \bfmu), \cdot) \right\rangle = 0 \qquad  \forall v \in \mathcal{T}_{\hat{q}(\bftheta(t, \bfmu), \cdot)}\widehat{\mathcal{M}}\,,
\end{equation}
where $\mathcal{T}_{\hat{q}(\bftheta(t, \bfmu), \cdot)}\widehat{\mathcal{M}}$ denotes the tangent space of $\widehat{\mathcal{M}}$ at $\hat{q}(\bftheta(t, \bfmu), \cdot)$.
Notice that conditions \eqref{eq:DF:Conditions} do not necessarily determine a unique time derivative if the dimension of the space $\mathcal{T}_{\hat{q}(\bftheta(t, \bfmu), \cdot)}\widehat{\mathcal{M}}$ is less than the number $n$ of weights; a fact that we will revisit in \Cref{sec:RegDF} and that we have already encountered in \Cref{subsec:shiftedPOD:ROM}.

The choice of the inner product $\langle \cdot, \cdot \rangle$ of the orthogonality conditions \eqref{eq:DF:Conditions} is critical. Out of convenience, a common choice is the $\Ltwo$ inner product; however, the $\Ltwo$ inner product may be inappropriate for some problems in the sense that the orthogonality conditions on the residual do not control the error in the norm of interest. At the same time, numerically evaluating the orthogonality conditions formulated with other inner products can be challenging. We refer to \citet{10.1093/imanum/draf073} for an in-depth discussion. 

The orthogonality conditions \eqref{eq:DF:Conditions} are instantaneous in the sense that they act locally in time: at each time $t$, a derivative $\dot{\bftheta}(t, \bfmu)$ is chosen that sets the residual orthogonal to all elements in the tangent space $\mathcal{T}_{\hat{q}(\bftheta(t, \bfmu), \cdot)}\widehat{\mathcal{M}}$, which defines an evolution equation for the weight vector $\bftheta(t, \bfmu)$. The instantaneous conditions on the residual stand in contrast to global-in-time approaches, such as tensor-based space-time approaches \cite{NOUY20101603} and physics-informed neural networks \cite{RAISSI2019686}, which determine a global approximation by minimizing a residual-based objective over the space--time domain. A detailed comparison between instantaneous (or sequential-in-time) and global-in-time methods in the context of nonlinear parametrizations can be found in \citet{24OTDDTO}; for a discussion in the context of linear model reduction we refer to \cite{Glas2017}.

\begin{remark}
While the Dirac--Frenkel variational principle provides a natural and widely used mechanism for deriving reduced dynamics, it represents only one particular choice within a broader geometric setting.
\citet{BUCHFINK2024134299} developed a differential--geometric framework for nonlinear model reduction that adopts a more general viewpoint. In this framework, reduced models are not required to arise from the Dirac--Frenkel variational principle; instead, more general reduction maps defined on the tangent bundle of the solution manifold are admitted, which can serve as a blue print for structure-preserving model reduction using nonlinear parametrizations.
\end{remark}

\subsubsection{Dirac--Frenkel and instantaneous residual minimization}
Given the residual function \eqref{eq:DF:ResidualFun}, instantaneous residual minimization means minimizing the residual norm at the current time $t$,
\begin{equation}
    \label{eq:NG:ResNormMin}
    \dot{\bftheta}(t, \bfmu) \in \operatorname*{arg\,min}_{\bfeta \in \mathbb{R}^n} \|r(\bftheta(t, \bfmu), \bfeta, \cdot)\|^2\,,
\end{equation}
where $\|\cdot\|$ is the norm induced by the inner product $\langle \cdot, \cdot \rangle$. 
Setting the first-order optimality condition of \eqref{eq:NG:ResNormMin} to zero, i.e.,
\[
\nabla_{\bfeta} \|r(\bftheta(t, \bfmu), \bfeta, \cdot)\|^2 = 0,
\]
leads to the system of equations
\begin{equation}\label{eq:NG:NGConditions}
\langle \partial_{\theta_i} \hat{q}(\bftheta(t, \bfmu), \cdot), r(\bftheta(t, \bfmu), \bfeta, \cdot)\rangle = 0\,,\qquad i = 1, \dots, n.
\end{equation}
Now notice that if the $\widehat{\mathcal{M}}$ is equipped with a suitable manifold structure, then the tangent space at $\hat{q}(\bftheta(t, \bfmu), \cdot)$ is spanned by the partial derivatives of $\hat{q}$ in all directions of the weight vector components,
\begin{equation}\label{eq:TangentSpaceDF}
\mathcal{T}_{\hat{q}(\bftheta(t, \bfmu), \cdot)}\widehat{\mathcal{M}} = \operatorname{span}\left\{\partial_{\theta_1}\hat{q}(\bftheta(t, \bfmu), \cdot), \dots, \partial_{\theta_n}\hat{q}(\bftheta(t, \bfmu), \cdot)\right\}. 
\end{equation}
Thus, conditions \eqref{eq:NG:NGConditions} derived from the instantaneous residual norm minimization \eqref{eq:NG:ResNormMin} are equivalent to the conditions imposed by the Dirac--Frenkel variational principle \eqref{eq:DF:Conditions}.

If the parametrization is linear, i.e., of the form~\eqref{eq:NMOR:LinParam}, then~\eqref{eq:NG:NGConditions} is equivalent to,
\begin{align*}
        \langle \varphi_i,\partial_t\hat{q}(\bftheta(t,\param),\cdot) - f(\cdot,\hat{q}(\bftheta(t,\param),\cdot);\param)\rangle = 0,\qquad i=1,\ldots,\stateDimRed,
\end{align*}
which, provided that the basis functions $\varphi_i$ are orthonormal, corresponds (in coordinates) exactly to the Galerkin reduced model~\eqref{eq:Prelim:HatFGalerkin}. Consequently, the Dirac--Frenkel variational principle can be interpreted as a generalization of the Galerkin orthogonality condition.

\begin{remark}
The Dirac--Frenkel variational principle, and more generally instantaneous residual minimization, provides a unifying framework for deriving reduced dynamics on nonlinear trial manifolds. In particular, several methods discussed in the sections on transformation-based reduction (\Cref{sec:Template}) and online adaptive bases (\Cref{sec:ADEIM}) can be interpreted as instances of instantaneous residual minimization. 
Nevertheless, when the nonlinear parametrization has additional exploitable structure, the reduced dynamics can often be written in a more specific form and solved more efficiently (e.g., matrix factorizations in dynamic low-rank approximations or shift parameters in transformation-based model reduction). We therefore treat such more structured nonlinear parametrizations separately in this survey. 
\end{remark} 

\subsubsection{Literature overview of instantaneous residual minimization with nonlinear parametrizations}
There is a vast literature on the Dirac--Frenkel variational principle and its applications in science and engineering. In this section, we focus on works that develop schemes based on the Dirac-Frenkel variational principle for parametrizations that are not necessarily trained offline on snapshot data, i.e., approaches that are closer in spirit to PDE solvers. From the perspective of nonlinear model reduction, however, methods that apply the Dirac--Frenkel principle to trained nonlinear parametrizations are of particular interest; these are reviewed separately in \Cref{sec:DFWithPretrained}. 

The computational chemistry community has developed methods based on the Dirac-Frenkel variational principle for numerically solving, e.g., the Schr\"odinger equation  \cite{10.1063/1.431911,MEYER199073,BECK20001,9073ba01-c8c8-3f30-b15c-e4b52a44e9da,lubich2008quantum}. 
The Dirac--Frenkel variational principle is also used in dynamic low-rank approximations and dynamic orthogonal decompositions with matrix factorizations as parametrizations, which we survey in the context of time-dependent basis methods in \Cref{sec:ADEIM:LitOverview}. 

The moving finite-element method of \citet{MilM81} parametrizes solutions by finite-element coefficients together with time-dependent nodal positions, and their evolution equations are obtained with the Dirac--Frenkel variational principle. But Dirac--Frenkel can be applied to more generic nonlinear parametrizations. For example,  \citet{anderson2021evolution}  use Dirac--Frenkel dynamics with morphing-shape parametrizations. Neural-network parametrizations are of particular interest due to their expressivity, as studied in the context of Dirac--Frenkel schemes in \citet{LeeC20} and \citet{PhysRevE.104.045303}. In a similar vein, the Neural Galerkin schemes introduced in \citet{BPE22NG} formulate Dirac--Frenkel evolution equations for the weights of neural-network parametrizations, which have led to follow up works on high-dimensional problems \cite{wen2023coupling}, nonlinear model reduction \cite{berman2024colora,WEDER2025114249}, and fast and structure-preserving online solvers \cite{berman2023randomized,SSBP23NGE}. Related ideas have also been used for filtering in \citet{FilteringNG}. Various aspects and challenges of using the Dirac--Frenkel variational principle with neural networks have been addressed: \citet{KAST2024112986} propose to use a Laplace--Beltrami positional embedding to handle geometrically complex domains; \citet{finzi2023a} introduce a re-training step to improve numerical stability; and \citet{quantuminspired_NG} develop  a variational Monte Carlo-type sampling estimator to realizing Dirac--Frenkel parameter dynamics for high-dimensional PDE problems. 
An important variant of the Dirac--Frenkel variational principle when applied to generic nonlinear parametrization such as neural networks is its regularized version, which is systematically studied in \citet{feischl2024regularized,lubich2025regularizeddynamicalparametricapproximation} and which we revisit in \Cref{sec:RegDF}.

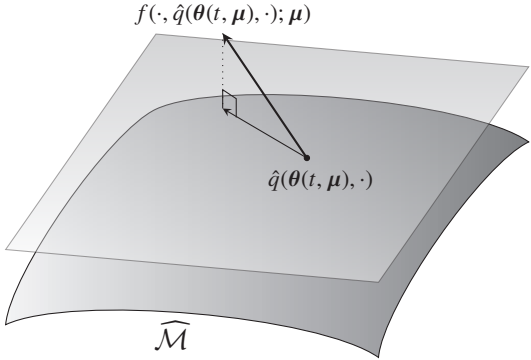
\begin{SCfigure}\centering
\resizebox{0.55\columnwidth}{!}{\begin{tikzpicture}
\tikzset{x={(175:1cm)},y={(55:.7cm)},z={(90:1cm)}}
\shade[left color=black!4,right color=black!60,looseness=.5, draw=black]  (2.5,-2.5,-1) 
to[bend left] (2.5,2.5,-1)
to[bend left] coordinate (mp) (-2.5,2.5,-1)
to[bend right] (-2.5,-2.5,-1)
to[bend right] coordinate (mm) (2.5,-2.5,-1)
-- cycle;
\shade[opacity=0.4, draw=black] (2.5,-2.5,0) -- (2.5,2.5,0) -- (-2.5,2.5,0) -- (-2.5,-2.5,0) -- cycle;
\draw[-stealth, thick] (-0.5,0.08,0) -- coordinate[pos=.3] (f) (1.0,1.0,1);
\draw[-stealth] (-0.5,0.08,0) -- coordinate[pos=.3] (f) (1.0,1.0,0);
\draw[dotted] (1.0,1.0,1) -- (1.0,1.0,0);
\draw[solid]  (1.0,1.0,0) -- (1.0,1.0,0.25) -- (0.75,0.83,0.25) -- (0.75,0.83,0);
\filldraw (-0.5,0.08,0) circle (1pt);
\node[right, scale=0.8] at (-0.1,-0.5,0) {{$\hat{q}(\bftheta(t, \bfmu), \cdot)$}};
\node[right, align=left, scale=0.8] at (1.7,-0.3,0.2) {{}};
\node[above, scale=0.8] at (1.0,1.0,1.0) {{$f(\cdot, \hat{q}(\bftheta(t, \bfmu),\cdot); \bfmu)$}};
\draw (-1.5,2.5,0) node[above, rotate=-4] {}; 
\draw (0,-4.1,0) node[right] {$\widehat{\mathcal{M}}$};
\end{tikzpicture}}
\caption{The Dirac--Frenkel variational principle evolves the weight vector $\bftheta(t, \bfmu)$ by instantaneous residual minimization: it chooses the time derivative $\dot{\bftheta}(t, \bfmu)$ as the coordinates of the orthogonal projection of the right-hand side function $f(\cdot, \hat{q}(\bftheta(t, \bfmu), \cdot); \bfmu)$ onto the tangent space of the trial manifold $\widehat{\mathcal{M}}$ at the current solution $\hat{q}(\bftheta(t, \bfmu), \cdot)$.}
\label{fig:DiracFrenkel}
\end{SCfigure}

\subsection{Neural Galerkin schemes}\label{sec:NeuralGalerkin}
We now discuss Neural Galerkin schemes that instantiate the Dirac--Frenkel variational principle using neural-network parametrizations \cite{BPE22NG,BERMAN2024389}. 
We focus on one particular numerical aspect, which is central to Neural Galerkin schemes as well as other schemes based on the Dirac--Frenkel variational principle: The orthogonality conditions given in \eqref{eq:DF:Conditions} involve inner products that must be evaluated numerically. Neural Galerkin schemes approximate these inner products, or equivalently the norm in the residual minimization formulation, via quadrature or collocation (sampling) over the spatial domain.

\subsubsection{Neural Galerkin: An instantiation of the Dirac--Frenkel variational principle}
Motivated by neural-network parametrizations, \citet{BPE22NG} introduce Neural Galerkin schemes that impose the Dirac-Frenkel dynamics on the neural-network weights $\bftheta(t, \bfmu)$. This means that the neural-network weights depend on time and are evolved following the Dirac-Frenkel variational principle. We remark that the earlier work \citet{PhysRevE.104.045303} also uses the Dirac--Frenkel variational principle in the context of neural-network parametrizations. 

From a numerical perspective, Neural Galerkin schemes directly consider instantaneous residual norm minimization \eqref{eq:NG:ResNormMin}, so that the time derivative $\dot{\bftheta}(t, \bfmu)$ is a solution of the time-dependent least-squares problem
\begin{equation}\label{eq:NG:ResExp}
\dot{\bftheta}(t, \bfmu) \in \operatorname*{arg\,min}_{\bfeta \in \mathbb{R}^n} \|\nabla_{\bftheta}\hat{q}(\bftheta(t, \bfmu), \cdot)^{\top} \bfeta - f(\cdot, \hat{q}(\bftheta(t, \bfmu), \cdot); \bfmu)\|^2\,,
\end{equation}
which is the same optimization problem as \eqref{eq:NG:ResNormMin} except that the residual function \eqref{eq:DF:ResidualFun} is written out explicitly.
The formulation \eqref{eq:NG:ResExp} highlights that this is a least-squares problem that is linear in the unknown $\bfeta$, even though a nonlinear parametrization $\hat{q}$ is used to represent the PDE solution. Furthermore, the formulation \eqref{eq:NG:ResExp} makes explicit that the time derivative $\dot{\bftheta}(t, \bfmu)$ is the coefficient vector of the best approximation (projection) of the right-hand side function onto the space \eqref{eq:TangentSpaceDF} spanned by the component functions of the gradient, as visualized in \Cref{fig:DiracFrenkel}. 

The normal equations of the least-squares problem \eqref{eq:NG:ResExp} lead to a dynamical system for $\dot{\bftheta}(t, \bfmu)$, 
\begin{equation}\label{eq:NG:NormalEquations}
\bfE(\bftheta(t, \bfmu))\dot{\bftheta}(t, \bfmu) = \bfG(\bftheta(t, \bfmu); \bfmu)\,,
\end{equation}
with matrices
\begin{align}
\bfE_{ij}(\bftheta(t, \bfmu)) = &\langle \partial_{\theta_i}\hat{q}(\bftheta(t, \bfmu), \cdot), \partial_{\theta_j}\hat{q}(\bftheta(t, \bfmu), \cdot) \rangle\,,\quad i, j = 1, \dots, n\,,\label{eq:NG:MatE}\\
\bfG_i(\bftheta(t, \bfmu); \bfmu) = &\langle \partial_{\theta_i}\hat{q}(\bftheta(t, \bfmu), \cdot), f(\cdot, \hat{q}(\bftheta(t, \bfmu), \cdot); \bfmu) \rangle\,,\quad i = 1, \dots, n\,.\label{eq:NG:MatF}
\end{align} 
The dynamical system \eqref{eq:NG:NormalEquations} can be integrated in time to obtain a weight trajectory~$\bftheta(t, \bfmu)$ so that at each point in time the orthogonality conditions \eqref{eq:NG:NGConditions} are satisfied. Alternatively, the least-squares problem \eqref{eq:NG:ResExp} can be integrated in time, which  typically leads to better conditioned intermediate numerical problems at each time step than solving the normal equations \cite{BERMAN2024389}. We remark once more that \eqref{eq:NG:NormalEquations} can be under-determined; see the forthcoming \Cref{sec:RegDF}.

\subsubsection{Estimating inner products with sampling points}\label{sec:NGEstInnerProducts}
In this section, we discuss collocation strategies for numerically estimating the norm in the least-squares problem \eqref{eq:NG:ResExp} in the context of Neural Galerkin schemes.

\paragraph{The need for numerically estimating inner products}
Consider the least-squares formulation \eqref{eq:NG:ResExp} corresponding to the Neural Galerkin schemes.
The objective of the least-squares problem is formulated via the norm $\|\cdot\|$, which is induced by the inner product $\langle \cdot, \cdot \rangle$ that underlies the residual orthogonality conditions \eqref{eq:DF:Conditions} given by the Dirac--Frenkel variational principle.
Thus, the objective has to be numerically evaluated to solve the least-squares problem for $\dot{\bftheta}(t, \bfmu)$. 

As \citet{BERMAN2024389} point out, there is an analog in numerical analysis: In finite-element methods, assembling Galerkin systems requires evaluating inner products of basis functions to compute entries of mass and stiffness matrices. For linear parametrizations, such as linear combinations of local basis functions that are centered at fixed grid points, these inner products can often be pre-computed and reused efficiently \cite{FEMTheory}. 
In contrast, for nonlinear parametrizations, this pre-computation is no longer possible because the test space changes over time with the weights of the parametrization, and so the quantities in the inner products change over time and cannot be directly re-used. Furthermore, because of the nonlinearity of the parametrizations, the superposition principle is lost, and so the residual norm cannot be obtained from separate components. 

There are situations when the Dirac--Frenkel variational principle is applied, where there is a closed form of the residual so that it can be evaluated exactly without numerical computations. For example, in some instances, nonlinear parametrizations based on Gaussian wave packets can be chosen such that they lead to closed form residual norms for variants of the Schr\"odinger equations; see, e.g., \citet{lubich2008quantum,lasser_lubich_2020}. In other cases, symbolic computations can be used as by \citet{ANDERSON2024112649}.
However, in many settings, especially when the nonlinear parametrizations are trained on snapshot data in an offline phase as in model reduction, there is no closed-form solution, so the residual norm must  be computed numerically.

\paragraph{Estimating inner products via static collocation points}
For the numerical evaluation of the norm in~\eqref{eq:NG:ResExp}, let us consider collocation points $\bfx_1, \dots, \bfx_{\spts} \in \Omega$ that have been sampled from a measure $\nu$ supported on $\Omega$. For example, if $\nu$ is the uniform measure over $\Omega$, then the points $\bfx_1, \dots, \bfx_{\spts}$ are sampled uniformly in $\Omega$. We refer to the collocation points as static because they are fixed over time, and additionally independent of the parameter $\bfmu$. 

The collocation points are used to estimate the norm in~\eqref{eq:NG:ResExp}, which leads to the optimization problem with the empirical objective
\begin{equation}\label{eq:NG:EmpObj}
    \min_{\bfeta \in \mathbb{R}^n} \frac{1}{\spts} \sum_{i = 1}^{\spts} \left|\nabla_{\bftheta} \hat{q}(\bftheta(t, \bfmu), \bfx_i)^{\top} \bfeta - f(\bfx, \hat{q}(\bftheta(t, \bfmu), \cdot); \bfmu)|_{\bfx = \bfx_i}\right|^2\,. 
\end{equation}
If the points $\bfx_1, \dots, \bfx_{\spts}$ are sampled from $\nu$, then the empirical objective in \eqref{eq:NG:EmpObj} is a Monte Carlo estimator. Estimating the objective \eqref{eq:NG:ResExp} with a Monte Carlo estimator~\eqref{eq:NG:EmpObj} is analog to situations found in machine learning, where population loss functions are estimated with empirical loss functions \cite{NIPS1991_ff4d5fbb}. 

Numerical quadrature offers an alternative to Monte Carlo estimation: The points $\bfx_1, \dots, \bfx_{\spts}$ are selected based on some quadrature rule that accounts for the measure~$\nu$. Using a quadrature rule leads to analogous estimators as \eqref{eq:NG:EmpObj} except that there can be different weights than $1/{\spts}$. Numerical quadrature for evaluating inner products in the context of the Dirac--Frenkel variational principle are considered in, e.g.,  \citet{PhysRevE.104.045303} and \citet{SSBP23NGE}.

\begin{figure}[t]
\fbox{\parbox{0.98\columnwidth}{
\begin{tikzpicture}
\node (theta0) {$\bftheta(t_0)$};
\node[below of = theta0] (x0) {$\{\bfx_i(t_0)\}_{i = 1}^{\spts}$};
\node[right of = theta0, node distance=1.0in] (theta1) {$\bftheta(t_1)$};
\draw[->] (theta0.east) -- (theta1.west);
\node[right of = x0, node distance=1.0in] (x1) {$\{\bfx_i(t_1)\}_{i = 1}^{\spts}$};
\draw[->] (x0.east) -- (x1.west);
\node[right of = theta1, node distance=1.0in] (theta2) {$\bftheta(t_2)$};
\draw[->] (theta1.east) -- (theta2.west);
\node[right of = x1, node distance=1.0in] (x2) {$\{\bfx_i(t_2)\}_{i = 1}^{\spts}$};
\draw[->] (x1.east) -- (x2.west);
\node[right of = theta2, node distance=1.0in] (theta3) {$\dots$};
\draw[->] (theta2.east) -- (theta3.west);
\node[right of = x2, node distance=1.0in] (x3) {$\dots$};
\draw[->] (x2.east) -- (x3.west);
\node[below of = x2, node distance=1.2in] (fig1) {\includegraphics[width=0.95\columnwidth]{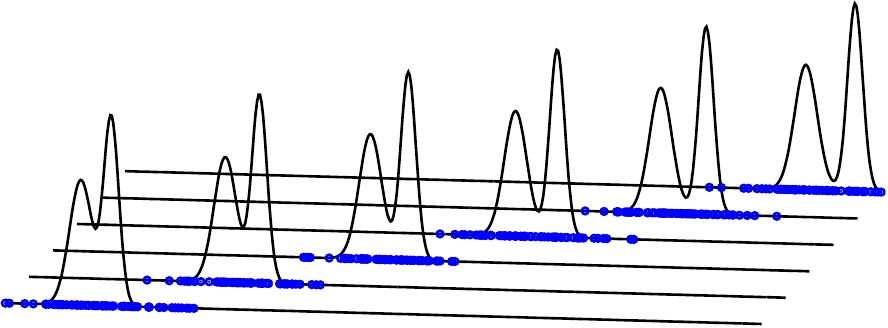}};
\node[below of = fig1, node distance = 1.0in] {spatial domain};
\node[right of = fig1, yshift=-0.75in, node distance = 2.2in] {time};
\end{tikzpicture}}}
\caption{Neural Galerkin schemes evolve the weight vector as prescribed by the Dirac–Frenkel variational principle and can simultaneously adapt sampling points corresponding to a potential informed by the residual. This adaptation concentrates the sample points in regions that dominate the residual and so yield more efficient residual-norm estimates than static sampling.}
\label{fig:AdaptiveNG}
\end{figure}
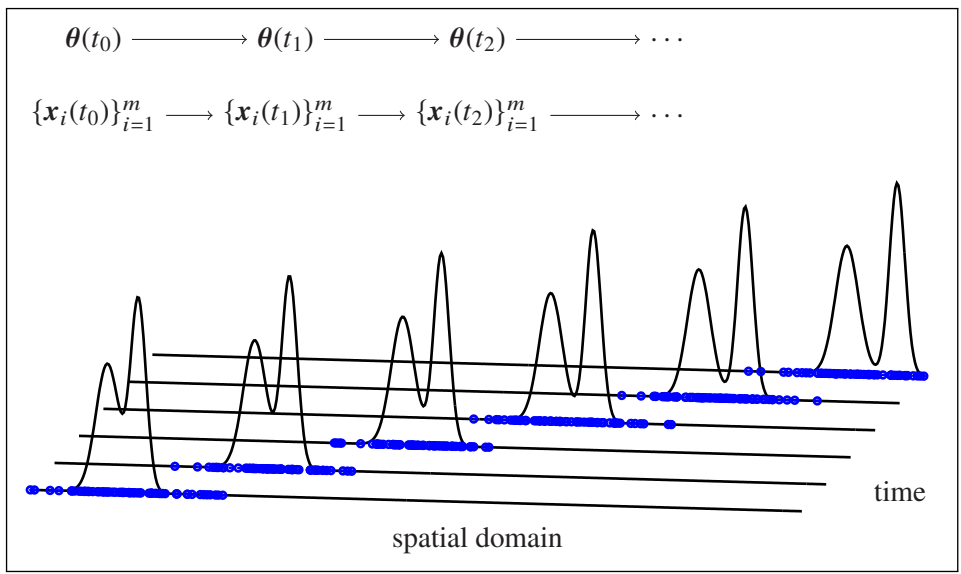

\subsubsection{Adaptive sampling points}\label{sec:NGWithAdaptiveSamplingPoints}
Neural Galerkin schemes support evolving the sampling points together with the weight vector $\bftheta(t, \bfmu)$; see \Cref{fig:AdaptiveNG}. 

\paragraph{Limitations of static collocation points for transport-dominated problems}
Returning to the motivating example in \Cref{fig:Need:NumExp} from \Cref{sec:Intro:LinNotEnough}, it is common that transport-dominated problems exhibit moving wave fronts or other local, coherent features---sharp solution gradients, bumps, wave fronts---that travel through the spatial domain. 

When using a set of static collocation points, which means that they are fixed for all times $t$, then these points must be chosen so that the moving local features are adequately resolved at all times. However, ensuring uniform accuracy at all times becomes difficult because the local feature may traverse any and all regions of the spatial domain and so a fine resolution across the whole spatial domain is needed, which often requires a large number of sampling points \cite{BPE22NG,wen2023coupling}.
In the linear transport example shown in \Cref{fig:Need:NumExp}, a uniform static sampling must include points in all parts of the domain at all times to track the narrow Gaussian bump as it moves. Even though most of these points at most of the time steps will lie in regions where the solution is close to zero, and so contribute little useful information while still requiring gradient and right-hand side evaluations to empirically estimate the objective of \eqref{eq:NG:ResExp} as in \eqref{eq:NG:EmpObj}. 
Consequentially, relying on static sampling points can require large numbers of sampling points to resolve moving local features. 

\begin{remark}
   A similar sampling challenge can be found in physics-informed machine learning methods that aim to minimize the residual norm over spatial (and time) domains. If solution features are local, then a large number of collocation points, and thus residual-norm evaluations, is needed, which can become a computational bottleneck \cite{https://doi.org/10.1111/mice.12685,pmlr-v145-rotskoff22a,WU2023115671}. 
\end{remark}

\paragraph{Estimating inner products with adaptive collocation points}
Estimating inner products with time-adaptive collocation points offers one approach to keep the number of collocation points low even if solution features are local. \citet{wen2023coupling} formulate Neural Galerkin schemes via time-dependent measures $\nu_{\bftheta(t, \bfmu), \dot{\bftheta}(t, \bfmu)}$ that are allowed to change with the weights $\bftheta(t, \bfmu)$ and the derivative $\dot{\bftheta}(t, \bfmu)$. Adaptive sampling in the context of nonlinear model reduction is also addressed in \citet{bon2025stablenonlineardynamicalapproximation}, where a sampling criterion is proposed that directly targets conditioning of the reduced dynamics and in \citet{BlaSU21b}, where interpolation points are updated over time for the efficient estimation of nonlinear terms. Also relevant is the work \citet{gruhlke2024optimalsamplingstochasticnatural} that introduces sampling for natural gradient descent, which corresponds to analog dynamics as given by Neural Galerkin schemes \cite{24OTDDTO}. We emphasize that in numerical analysis, the adaptive choice of sampling points in the context of PDE approximations is used for instance in moving finite element methods \cite{MilM81}.

Following \citet{wen2023coupling}, we formulate the orthogonality conditions via a time-dependent measure, which leads to the analogous least-squares problem of \eqref{eq:NG:ResExp} except that the $\Ltwo$ norm $\|\cdot\|_{\nu_t}$ now depends on $\nu_t$, i.e.,
\[
\min_{\dot{\bftheta}(t, \bfmu) \in \mathbb{R}^n} \|\nabla_{\bftheta}\hat{q}(\bftheta(t, \bfmu), \cdot)^{\top} \dot{\bftheta}(t, \bfmu) - f(\cdot, \hat{q}(\bftheta(t, \bfmu), \cdot); \bfmu)\|_{\nu_t}^2\,,
\]
where we denote $\nu_{\bftheta(t, \bfmu), \dot{\bftheta}(t, \bfmu)}$ as $\nu_t$ for brevity.
\citet{wen2023coupling} show that if the residual can be set to zero with respect to a static measure $\nu$, then the optimum is invariant to changes in the measure under standard conditions such as the measures are fully supported on the spatial domain $\Omega$. 
It is further argued in \citet{wen2023coupling} that the residual being close to zero is a reasonable assumption because one is typically interested in cases where the parametrization is so rich that the residual can be kept small. We note that alternatively to formulating the objective in a time-dependent norm, one can derive an importance-sampling estimator of the objective with a time-dependent biasing measure while keeping the base measure fixed; we refer to \cite{BPE22NG} where such an importance-sampling approach is discussed.

One major question is how to propagate the measures $\nu_t$ forward in time so that samples from $\nu_t$ lead to accurate estimators of the residual norm. \citet{wen2023coupling} propose to let the residual itself drive the evolution of $\nu_t$ by defining $\nu_t$ as a Gibbs measure 
\[
\nu_t \propto \exp\left(-U_{\bftheta(t, \bfmu), \dot{\bftheta}(t, \bfmu)}\right)
\]
with the potential $U_{\bftheta,\dot{\bftheta}}$ taken as the squared residual, i.e.,
\begin{equation}\label{eq:NG:PotentialResV}
U_{\bftheta(t, \bfmu), \dot{\bftheta}(t, \bfmu)}(\bfx) = |r(\bftheta(t, \bfmu), \dot{\bftheta}(t, \bfmu), \bfx)|^2\,.
\end{equation}
High probability regions of the measure $\nu_t$ correspond to spatial locations with large absolute residual. Other potentials could be used, for example by directly sampling from the solution function if it is a probability density function \cite{BPE22NG} or by building on other sampling criteria such as the stability-based criteria introduced in \citet{bon2025stablenonlineardynamicalapproximation}.

A formal way to evolve the time-dependent measure is given via the Fokker--Planck equation driven by the potential \eqref{eq:NG:PotentialResV},
\begin{equation}\label{eq:NG:Adapt:FPE}
\partial_t \nu_t = \gamma \nabla \cdot \left(\nabla \nu_t + \nu_t \nabla U_{\bftheta(t, \bfmu), \dot{\bftheta}(t, \bfmu)}\right)\,,
\end{equation}
which can be coupled to the weight dynamics,
\begin{equation}\label{eq:NG:Adapt:CoupledDynamics}
\begin{aligned}
\bfE_t(\bftheta(t, \bfmu))\dot{\bftheta}(t, \bfmu) &=  \bfG_t(\bftheta(t, \bfmu); \bfmu)\,,\\
\partial_t \nu_t &= \gamma \nabla \cdot \left(\nabla \nu_t + \nu_t \nabla U_{\bftheta(t, \bfmu), \dot{\bftheta}(t, \bfmu)}\right)\,,
\end{aligned}
\end{equation}
where the matrices $\bfE_t$ and $\bfG_t$ given in \eqref{eq:NG:MatE} and \eqref{eq:NG:MatF}, respectively, now also depend on time because they are formulated with respect to the inner product associated with $\nu_t$.  
The parameter $\gamma$ controls the time scale on which the measure $\nu_t$ adapts relative to the weights $\bftheta(t, \bfmu)$. 

Crucially, it is unnecessary to solve the Fokker-Planck equation \eqref{eq:NG:Adapt:FPE} explicitly because only samples are needed from the density corresponding to $\nu_t$. Thus, the measure $\nu_t$ can be approximated via the empirical measure  
\[
\hat{\nu}_t = \frac{1}{{\spts}} \sum_{i=1}^{\spts} \delta_{\bfx_i(t)}
\]
given by time-dependent samples $\bfx_1(t), \dots, \bfx_{\spts}(t)$.
Because the potential \eqref{eq:NG:PotentialResV} is given via the residual, one can evaluate the gradient of the residual to obtain the score $\nabla \log \nu_t$ of $\nu_t$ and therefore advance the samples $\bfx_1(t), \dots, \bfx_{\spts}(t)$ over time with particle methods rather than solving the Fokker-Planck equation \eqref{eq:NG:Adapt:FPE}. 
\citet{wen2023coupling} propose to use Langevin samplers and Stein variational gradient descent, both of which propagate samples via the score $\nabla \log \nu_t$.
In this way, the weight trajectory $\bftheta(t,\bfmu)$ and the set of samples $\{\bfx_i(t)\}_{i=1}^{\spts}$ are advanced together in time, approximating the coupled dynamics \eqref{eq:NG:Adapt:CoupledDynamics}, as illustrated in Figure~\ref{fig:AdaptiveNG}.

\subsection{Semi and fully discrete instantaneous residual minimization}\label{sec:DiscreteDF}
We now discuss semi-discrete and fully discrete variants of the least-squares problem resulting from the Dirac--Frenkel variational principle, where we follow the formulation and terminology used by Neural Galerkin schemes.

\paragraph{Semi-discrete least-squares problem}
Given ${\spts}$ sampling points $\bfx_1, \dots, \bfx_{\spts} \in \Omega$, which can be obtained statically as in \Cref{sec:NGEstInnerProducts} or adaptively as in \Cref{sec:NGWithAdaptiveSamplingPoints}, consider the so-called batch gradient
\begin{equation}
    \label{eq:NG:EmpiricalMatrixObj:BatchJacobian}
    \boldsymbol{J}(\bftheta) =  \begin{bmatrix}
        \nabla_{\bftheta}\hat{q}(\bftheta, \bfx_1)^{\top}\\
        \vdots\\
        \nabla_{\bftheta}\hat{q}(\bftheta, \bfx_{\spts})^{\top}
    \end{bmatrix} \in \mathbb{R}^{{\spts} \times n}\,,
\end{equation}
which has as rows the transpose of the gradient $\nabla_{\bftheta}\hat{q}(\bftheta, \cdot)$ of the nonlinear parametrization $\hat{q}$ evaluated at the ${\spts}$ sampling points. The name batch gradient stems from the consideration of approximating the objective as in \eqref{eq:NG:EmpObj}, where the $\bfx_1, \dots, \bfx_{\spts}$ are the current batch of samples, which can in principle change over the time steps (and optimization iterations in the time-discrete setting; see below). 
Analogously, we define the batch right-hand side function as
\begin{equation}
    \bfF(\bftheta; \bfmu) =  \begin{bmatrix}
    f(\bfx, \hat{q}(\bftheta, \cdot); \bfmu)|_{\bfx = \bfx_1}\\
    \vdots\\
    f(\bfx, \hat{q}(\bftheta, \cdot); \bfmu)|_{\bfx = \bfx_{\spts}}
    \end{bmatrix}  \in \mathbb{R}^{{\spts}}\,.\label{eq:NG:EmpiricalMatrixObj:BatchRHS}
\end{equation}

Using the batch gradient and the batch right-hand side function, we can write the semi-discrete instantaneous residual minimization problem as 
\begin{equation}\label{eq:NG:EmpiricalMatrixObj}
\dot{\bftheta}(t, \bfmu) \in \operatorname*{arg\,min}_{\bfeta \in \mathbb{R}^n} \|\bfJ(\bftheta(t, \bfmu))\bfeta - \bfF(\bftheta(t, \bfmu); \bfmu)\|_2^2\,,
\end{equation}
which is a matrix representation of the empirical objective \eqref{eq:NG:EmpObj} on the sample points $\bfx_1, \dots, \bfx_{\spts}$ up to a constant scaling factor. Note that the norm $\|\cdot\|_2$ in \eqref{eq:NG:EmpiricalMatrixObj} denotes the Euclidean vector norm, which can be readily evaluated.
It is important to notice that the least-squares problem \eqref{eq:NG:EmpiricalMatrixObj} is linear in the unknown, a property that one often wants to preserve when discretizing in time.

\paragraph{Time discretization}
Consider the discrete time steps $0 = t_0 < t_1 < \dots < t_K = T$ with constant time-step size $\timeStep > 0$. 
An explicit time discretization of the least-squares problem \eqref{eq:NG:ResExp} with Runge-Kutta methods leads to a sequence of least-squares problems that are also linear in the unknown. For example, using the explicit Euler method leads to linear least-squares problem at all time steps $k = 1, \dots, K$,
\begin{equation}\label{eq:NG:TimeDisc:Exp}
    \min_{\delta \bftheta_k(\bfmu) \in \mathbb{R}^n} \|\bfJ(\bftheta_{k-1}(\bfmu)) \delta \bftheta_k(\bfmu) - \bfF(\bftheta_{k-1}(\bfmu); \bfmu)\|_2^2\,,
\end{equation}
with the update
\begin{equation}\label{eq:NG:UpdateTheta}
    \bftheta_k(\bfmu) = \bftheta_{k-1}(\bfmu) + \timeStep\delta \bftheta_k(\bfmu)\,.
\end{equation}
Thus, if discretized with an explicit scheme, only a linear least-squares problem has to be solved at each time step. 

In contrast, implicit time discretization leads to nonlinear optimization problems at each time step. Using the implicit Euler method for demonstration purposes, we obtain
\begin{equation}\label{eq:NG:DiscTimeImp}
    \min_{\delta \bftheta_k(\bfmu) \in \mathbb{R}^n} \|\bfJ(\bftheta_{k-1}(\bfmu) + \timeStep\delta\bftheta_{k}(\bfmu))  \delta \bftheta_k(\bfmu) - \bfF(\bftheta_{k-1}(\bfmu) + \timeStep\delta\bftheta_{k}(\bfmu); \bfmu)\|_2^2,
\end{equation}
where the update $\delta \bftheta_k(\bfmu)$ is applied as in \eqref{eq:NG:UpdateTheta}. 
In the case of an implicit time discretization, the time-discrete least-squares problem becomes nonlinear in the update $\delta \bftheta_k(\bfmu)$ because both the batch gradient and the batch right-hand side function are evaluated at the next iterate $\bftheta_{k}(\bfmu)$. 
Thus, numerically solving  \eqref{eq:NG:DiscTimeImp} amounts to solving a nonlinear least-squares problem in the update $\delta \bftheta_k(\bfmu)$, and therefore inherits the typical difficulties of fitting a nonlinear parametrization such as a neural network. In particular, the fitting has to be done at each time step. 
Consequently, there is a strong interest in explicit time integration, which preserves the linearity of the least-squares problem in discrete time for many  Dirac--Frenkel-based methods.

\subsection{Regularization and randomized time integration of Dirac--Frenkel dynamics}\label{sec:RegDF}
As mentioned before, the least-squares problems arising from the orthogonality conditions in the Dirac–Frenkel variational principle can become poorly conditioned or even rank-deficient. A particular example was already presented in the context of transformation-based methods in \Cref{subsec:shiftedPOD:ROM}. In general, the ill-conditioning motivates the use of regularization. We consider explicit regularization via Tikhonov and truncated SVD. We also briefly discuss randomized time integration via right sketching, which acts as a form of regularization while additionally reducing the number of unknowns solved for at each time step by restricting the update to a low-dimensional random subspace.

\subsubsection{Poor conditioning of least-squares problem and tangent space collapse}
Let us consider the time-discrete least-squares problem given in \eqref{eq:NG:TimeDisc:Exp}, which is obtained by applying the explicit Euler method to the time-continuous problem given in \eqref{eq:NG:EmpiricalMatrixObj}.

A major challenge for schemes building on Dirac--Frenkel is numerically solving the least-squares problem. Besides computational costs, the least-squares problem can be challenging to solve numerically because of poor conditioning. The conditioning of the least-squares problem \eqref{eq:NG:TimeDisc:Exp} is dominated by the condition number of the batch gradient $\bfJ(\bftheta_{k-1}(\bfmu))$ at the weight vector $\bftheta_{k-1}(\bfmu)$.
For nonlinear parametrizations, the batch gradient can become poorly conditioned or even (numerically) rank deficient, in which case the least-squares problem is underdetermined, and the weight update is not unique.
Furthermore, rank deficiency is closely related to the accuracy in Dirac--Frenkel methods: Error bounds show that the error is controlled by the projection error incurred when projecting the right-hand side onto the column space of the batch gradient (equivalently, the component functions of the gradient in continuous time or the tangent space). When the batch gradient is rank-deficient and the corresponding space onto which the right-hand side is projected cannot represent it well, then the residual and approximation error can be large, leading to low accuracy \cite{9073ba01-c8c8-3f30-b15c-e4b52a44e9da}. 

But rank deficiency has broader implications in the context of the Dirac--Frenkel variational principle. The component functions of the gradient, which span the tangent space under suitable regularity assumptions and which correspond to the columns of the batch gradient $\bfJ$ after discretization, also span the test space against which the residual is set orthogonal with conditions \eqref{eq:NG:NGConditions}. Thus, rank deficiency also means that the Dirac--Frenkel conditions enforce orthogonality only with respect to a low-dimensional test space so that residual components outside the low-dimensional test space are uncontrolled and can grow, leading to a deterioration of accuracy despite an accurate intermediate approximation at a given time. In the context of the Dirac--Frenkel variational principle, this degeneracy phenomenon is known as the matrix singularity problem \cite{KAY1989165,10.1063/1.449204,PhysRevE.101.023313}. The work \cite{24OTDDTO} refers to it as tangent space collapse because the space spanned by the component functions of the gradient becomes low dimensional and they span the tangent space under suitable conditions.

Regularization is a natural remedy and important in practice, but a major difficulty is that one does not solve a single, isolated least-squares problem here but instead one solves a sequence of least-squares problems whose conditioning and solutions depend on the previous time steps. As a result, any bias introduced by the regularization can accumulate over time and judicious techniques and a systematic and rigorous analysis is necessary to control and interpret the bias' effect across the whole time trajectory \cite{feischl2024regularized}. Furthermore, regularization needs to be used with caution because even commonly used regularizers can lead to a selection of updates (e.g., minimal-norm updates) that preserve parameter redundancies (for instance, identical or redundant components remain identical) and so prevent recovery from a degeneracy \cite{24OTDDTO}.

\subsubsection{Regularized Dirac--Frenkel variational principle}
A systematic and rigorous study of Tikhonov regularization of Dirac--Frenkel least-squares problems is conducted by   \citet{feischl2024regularized}. To mitigate ill-conditioning in the residual minimization problem, \citet{feischl2024regularized} apply a Tikhonov regularization on the time derivative $\dot{\bftheta}(t, \bfmu)$ in the time-continuous formulation and on the update (``velocity'') $\delta \bftheta_k(\bfmu)$ in time-discrete formulation, e.g.,  
\[
    \delta\bftheta_{k}(\bfmu) = \operatorname*{arg\,min}_{\bfeta \in \mathbb{R}^n} \|\bfJ(\bftheta_{k-1}(\bfmu))\bfeta - \bfF(\bftheta_{k-1}(\bfmu); \bfmu)\|_2^2 + \alpha \|\bfeta\|_2^2\,,
\]
with the update $\bftheta_{k}(\bfmu) = \bftheta_{k-1}(\bfmu) + \timeStep \delta\bftheta_{k}(\bfmu)$. 
Here, the parameter $\alpha\geq 0$ controls the Tikhonov regularization, which acts as a smoothing filter in the sense that directions associated with small singular values of the system matrix of the least-squares problems are damped, thereby suppressing large and potentially unstable weight updates.

The analysis by \citet{feischl2024regularized} emphasizes that an accurate approximation of the solution itself is not sufficient for reliable time evolution. One additionally needs that the right-hand side (equivalently, the time derivative of the solution) is well approximated by the tangent space along the solution trajectory because the error of approximating the right-hand side in the tangent space enters the overall approximation error; see also \citet{9073ba01-c8c8-3f30-b15c-e4b52a44e9da,24OTDDTO}. This observation is not specific to nonlinear model reduction. An analogous phenomenon arises already in linear model reduction, where an accurate approximation of the state alone does not guarantee an accurate approximation of the time evolution. This insight underlies the widespread use of Petrov--Galerkin projections in linear model reduction, where distinct trial and test spaces are chosen to improve the approximation of the system dynamics \cite{OttPR22}.

Furthermore, the regularization strength has to be chosen carefully: less regularization reduces bias but can lead to updates aligned with directions corresponding to small singular values, while more regularization helps stability but ultimately can prevent the residual from being small. \citet{feischl2024regularized} introduce an adaptation approach based on the defect, which is the minimal value of the regularized least-squares objective.   
\citet{feischl2024regularized} also relates Tikhonov regularization to truncated SVD regularization, in the sense that truncated SVD entirely discards directions corresponding to small singular values rather than only smoothing them as Tikhonov regularization does. 

\citet{lubich2025regularizeddynamicalparametricapproximation} target stiff problems and therefore consider regularization for implicit instead of explicit time integration schemes with Dirac--Frenkel. 
At each time step, a regularized nonlinear least-squares problem is solved, with a regularizer on the weight update. The nonlinear problems are solved with a few Gauss--Newton iterations.

\subsubsection{Regularization via random subspace approximations}
\citet{berman2023randomized}, and further developed in \citet{DongRandom}, propose to regularize the dynamics corresponding to the Dirac--Frenkel variational principle by solving for the weight update in a random subspace, which can be interpreted as sketching from the right. 

Building on the time-discrete problem \eqref{eq:NG:TimeDisc:Exp}--\eqref{eq:NG:UpdateTheta}, the weight update becomes 
\begin{equation}\label{eq:NG:ThetaRandUpdate}
    \bftheta_{k}(\bfmu) = \bftheta_{k-1}(\bfmu) +  \frac{\timeStep}{\ell}\sum_{i = 1}^{\ell} \delta\bftheta_{k}(\bfmu; \bfGamma_{k,i})\,,
\end{equation}
where $\delta\bftheta_{k}(\bfmu; \bfGamma_{k, 1}), \dots, \delta\bftheta_{k}(\bfmu; \bfGamma_{k, \ell})$ are $\ell$ independent random updates depending on the independent random embeddings $\bfGamma_{k, 1}, \dots, \bfGamma_{k, \ell}$. Notice that the average over these $\ell$ updates is used in \eqref{eq:NG:ThetaRandUpdate}, which helps to control the variance of the random update. The updates are obtained from the least-squares problem 
\begin{equation}\label{eq:NG:RandNG}
\delta\bftheta_{k}(\bfmu; \bfGamma) = \operatorname*{arg\,min}_{\bfeta \in \operatorname{Range}(\bfGamma)} \|\bfJ(\bftheta_{k-1}(\bfmu))\bfeta - \bfF(\bftheta_{k-1}(\bfmu); \bfmu)\|_2^2\,,
\end{equation}
which restricts the solution to be in the column span of the random matrix $\bfGamma$. 
The matrix $\bfGamma$ is of size $n \times s$ with $s \leq n$. 

Typical choices for $\bfGamma$ are Gaussian embeddings and random unitary embeddings. For such embeddings, \citet{DongRandom} develop bounds on the condition number of the sketched batch gradient $\bfJ(\bftheta_{k-1}(\bfmu))\bfGamma$ with respect to the spectrum of $\bfJ(\bftheta_{k-1}(\bfmu))$ and the sketching dimension $s$, which directly controls the condition number of the least-squares problem \eqref{eq:NG:RandNG}. 
The earlier work \citet{berman2023randomized} considers sparse updates $\delta\bftheta_{k}(\bfmu; \bfGamma)$ by choosing $\bfGamma$ to correspond to a randomized sparse coordinate embedding. 

The randomized formulation \eqref{eq:NG:RandNG} sketches the least-squares problem from the right, i.e., it restricts the weight update to a random subspace. Consequently, the unknown in \eqref{eq:NG:RandNG} is an $s$-dimensional vector and thus one has to solve for $s$ unknowns only, rather than for all $n$ unknowns corresponding to all components of the weight update. In contrast, there is a line of work on sketching from the left \cite{lam2024randomizedlowrankrungekuttamethods,CARREL2026114421} for dynamic low-rank approximations as well as other more generic nonlinear parametrizations \cite{lindsey2025mne}. Sketching from the left can help to improve the conditioning but the number of unknowns of the update remain $n$ instead of $s$.


\subsection{Instantaneous residual minimization with discretize-then-optimize (DtO) versus optimize-then-discretize (OtD) approaches}\label{sec:DtOOtD}
The Dirac--Frenkel variational principle can be interpreted as being of the optimize-then-discretize (OtD) type: one first derives a continuous-time evolution equation for the weights $\bftheta(t, \bfmu)$ by enforcing orthogonality of the residual against the tangent space of the trial manifold, and only afterwards discretizes this evolution in time. 

As discussed in detail in \citet{24OTDDTO}, a conceptually different class of instantaneous (or sequential-in-time) residual minimization schemes arise from a discretize-then-optimize (DtO) strategy: In DtO schemes, time is discretized at the PDE level first, which leads to a sequence of boundary value problems in the spatial variable. Only after the time discretization, the nonlinear parametrization is introduced and a residual minimization problem is formulated for each boundary value problem corresponding to the discrete time steps. Examples of DtO residual minimization schemes for nonlinear parametrizations include the works by \citet{kvaal2023need,chen2023implicit,pmlr-v235-chen24ad}. We note for linear parametrizations, DtO corresponds to the Rothe method \cite{Rothe1930,RotheDeuflhard}. 
Because residual minimization for each boundary value problem leads to a typically non-convex optimization problem at each time step, the computational costs of DtO schemes tend to be higher than the costs of OtD schemes that lead to a linear least-squares problem at each time step---at least for explicit time integration.
At the same time, DtO schemes can avoid ill-conditioning issues that typically arise in OtD schemes such as the tangent space collapse or matrix singularity issue. We refer to \citet{24OTDDTO} for an in-depth discussion about OtD versus DtO schemes. 

We remark that there is an analog of OtD versus DtO formulations in linear model reduction. \citet{carlberg2011efficient,CARLBERG2017693} introduce and analyze the concept of least-squares Petrov--Galerkin reduced models based on linear parametrizations. In the terminology used in \citet{CARLBERG2017693}, DtO corresponds to least-squares Petrov--Galerkin and the OtD to Galerkin formulations.

\subsection{Instantaneous residual minimization with trained nonlinear parametrizations for nonlinear model reduction}\label{sec:DFWithPretrained}
We now survey model reduction methods that apply instantaneous residual minimization in the online phase to parametrizations trained offline on snapshot data. The capability to train a parametrization typically requires it to depend on offline weights, denoted in vector form as $\woff$ (see \Cref{sec:NonParam:OffOnWeights}). Once the parametrization is trained, i.e., the offline weights are fixed, the online weights $\bftheta(t, \bfmu)$ can be obtained via instantaneous residual minimization in the online phase. 

This section focuses on trained nonlinear parametrizations obtained with neural-network autoencoders, autoencoders with a priori fixed feature maps, and neural representations. We emphasize that there are also transformation-based methods and online adaptive model reduction that yield trained nonlinear parametrizations that can be used online with residual minimization; we cover these in the respective sections. 

\subsubsection{Nonlinear parametrizations given by neural-network autoencoders}\label{sec:PreTrainAutoEncoder}
The pioneering work by \citet{LeeC20} proposes to use autoencoders parametrized by deep convolutional neural networks, train them on snapshot data in the offline phase, and then use the decoder function of the autoencoder to induce a nonlinear parametrization for representing the reduced solution in the online phase. The vector $\bftheta(t, \bfmu)$ in this setting corresponds to a latent state of the reduced solution and is computed with instantaneous residual minimization. While there is earlier work on using neural-network autoencoders for model reduction \cite{7799153,8062736}, the work by \citet{LeeC20} was the first that presented a comprehensive framework that sparked a series of publications on using autoencoders for nonlinear model reduction. A general differential--geometric framework for using autoencoders in nonlinear model reduction is introduced in \citet{BUCHFINK2024134299}. Our presentation follows loosely \citet{LeeC20}. 

\paragraph{Autoencoders}
An autoencoder is the composition $\Phi^{\uparrow} \circ \Phi^{\downarrow}\colon \mathbb{R}^N \to \mathbb{R}^N$ of two maps, the encoder, $\Phi^{\downarrow}\colon \mathbb{R}^N \to \mathbb{R}^n$, and the decoder, $\Phi^{\uparrow}\colon \mathbb{R}^n \to \mathbb{R}^N$. The encoder $\Phi^{\downarrow}$ maps a high-dimensional vector of dimension $N$ onto a low-dimensional vector of dimension $n \ll N$. The decoder operates in the opposite direction and lifts a low-dimensional vector to a high-dimensional one. 

In the setting of model reduction, one is typically interested in encoder-decoder pairs so that the composition $\Phi^{\uparrow} \circ \Phi^{\downarrow}$ is close to the identity map on data corresponding to the solution manifold~\eqref{eq:Prelim:SolManifold} induced by the full model: When applying the autoencoder to a state $\bfq(t, \bfmu) \in \mathbb{R}^N$ of the semi-discrete full model  \eqref{eqn:FOM}, then the encoder maps the full-model state $\bfq(t, \bfmu)$ onto a low-dimensional approximation $\tilde{\bfq}(t, \bfmu) = \Phi^{\downarrow}(\bfq(t, \bfmu)) \in \mathbb{R}^n$. The decoder lifts the $n$-dimensional approximation $\tilde{\bfq}(t, \bfmu)$ onto an $N$-dimensional vector $\Phi^{\uparrow}(\tilde{\bfq}(t, \bfmu)) \in \mathbb{R}^N$ that approximates the full-model state $\bfq(t, \bfmu)$. In view of linear model reduction in the Galerkin setting~\eqref{eq:Prelim:HatFGalerkin}, the encoder is given by $\Phi^{\downarrow}(\bfq) = \bfV^\top \bfq$ and the decoder by $\Phi^{\uparrow}(\bfq) = \bfV\bfq$. In the Petrov--Galerkin setting, one would replace the matrix $\bfV^\top$ in the encoder with another matrix $\bfW^\top$.

\paragraph{Training autoencoders on snapshot data}
\citet{LeeC20} parametrize the encoder and decoder maps with neural networks. In the spirit of \Cref{sec:NPM:Examples}, we denote the parametrized encoder as $\PhiDECoder^{\downarrow}(\cdot; \woff^{\downarrow})\colon \mathbb{R}^N \to \mathbb{R}^n$ and the parametrized decoder as $\PhiDECoder^{\uparrow}(\cdot; \woff^{\uparrow})\colon \mathbb{R}^n \to \mathbb{R}^N$ with offline weight vectors  $\woff^{\downarrow}$ and $\woff^{\uparrow}$ of appropriate size. The architecture of the parametrization depends on the problem at hand. For example, \citet{LeeC20,BucGH21} use deep convolutional neural networks. 
In contrast, \citet{KIM2022110841} use a shallow network with a sparse mask on the output layer for reducing costs during the online phase. A comparison of different architectures is performed by \citet{GRUBER2022114764}.

The offline weights $\woff^{\downarrow}$ and $\woff^{\uparrow}$ are trained on snapshot data \eqref{eq:QtrainSnapshotMatrix} from the full model \eqref{eqn:FOM}.  A standard loss function, which is used in \citet{LeeC20}, is the mean-squared error of the reconstruction of the snapshots, 
\begin{equation}\label{eq:IRM:AutoencoderLoss}
\min_{\woff^{\downarrow}, \woff^{\uparrow}} \left\|\bfQ_{\text{train}} - \PhiDECoder^{\uparrow}(\PhiDECoder^{\downarrow}(\bfQ_{\text{train}}; \woff^{\downarrow}); \woff^{\uparrow})\right\|_F^2\,,
\end{equation}
where $\bfQ_{\text{train}}$ is the snapshot matrix and the encoder and decoder are applied column-wise to their matrix inputs.  

Additional terms can be appended to \eqref{eq:IRM:AutoencoderLoss} to penalize structure violation. For example, \citet{BucGH21} add a penalty term that penalizes violations of symplecticity of the decoder, yielding a weakly symplectic autoencoder. A physics-informed penalization term is added in \citet{Lee_Carlberg_2021}. \citet{Otto2023} additionally train a projection operator that is tailored to the learned autoencoder.

\paragraph{Online residual minimization with autoencoder-induced nonlinear parametrizations}
\citet{LeeC20} propose to use the decoder $\PhiDECoder^{\uparrow}$ to induce a nonlinear parametrization. For the ease of exposition, we now drop the dependence on the offline weights in the notation of the decoder. The induced parametrization $\PhiDECoder^{\uparrow}\colon \mathbb{R}^n \to \mathbb{R}^N$ is vector-valued and does not take the spatial coordinate as an input. Instead, we can interpret $\PhiDECoder^{\uparrow}$ as 
\begin{equation}\label{eq:DecoderAsParam}
    \PhiDECoder^{\uparrow}(\bftheta(t, \bfmu)) = \begin{bmatrix}
    \hat{q}(\bftheta(t, \bfmu), \bfx_1) & \dots & \hat{q}(\bftheta(t, \bfmu), \bfx_N)
\end{bmatrix}^{\top}\,,
\end{equation}
where $\hat{q}(\bftheta(t, \bfmu), \cdot)\colon \Omega \to \mathbb{R}$ is an implicitly given scalar-valued nonlinear parametrization of the form \eqref{eq:NMOR:Param} that is evaluated at the points $\bfx_1, \dots, \bfx_N \in \Omega$ corresponding to the spatial coordinates over which the semi-discrete full model \eqref{eqn:FOM} is defined. We stress that $\hat{q}$, which allows other spatial coordinates than $\bfx_1, \dots, \bfx_N$ as inputs, is not directly available from the decoder function $\PhiDECoder^{\uparrow}$ in general. In particular, when autoencoders are typically used in model reduction, the autoencoders operate directly on the $N$-dimensional state vectors of the semi-discrete full model and are therefore inherently tied to the underlying full-model discretization. We refer to \Cref{sec:IRM:NeuralRep} for neural representations that are independent of the underlying full-model discretization. 

Correspondingly to the vector-valued parametrization induced by the decoder~$\PhiDECoder^{\uparrow}$, one can define the vector-valued residual function $\bfr\colon \mathbb{R}^n \times \mathbb{R}^n \times \mathcal{D} \to \mathbb{R}^N$,
\[
\bfr(\bftheta(t, \bfmu), \dot{\bftheta}(t, \bfmu); \bfmu) = \bfJ(\bftheta(t, \bfmu))\dot{\bftheta}(t, \bfmu) - \bff(\PhiDECoder^{\uparrow}(\bftheta(t, \bfmu)); \bfmu)\,,
\]
where $\bff$ is the right-hand side function of the semi-discrete full model \eqref{eqn:FOM} and $\bfJ(\bftheta(t, \bfmu))$ is the Jacobian matrix of $\PhiDECoder^{\uparrow}$ evaluated at $\bftheta(t, \bfmu)$. Recall that the residual function defined in \eqref{eq:DF:ResidualFun} for scalar-valued nonlinear parametrizations depends on the gradient of $\hat{q}$ with respect to $\bftheta(t, \bfmu)$, which becomes the batch gradient \eqref{eq:NG:EmpiricalMatrixObj:BatchJacobian} corresponding to sampling points in the semi-discrete case. Here, because the parametrization is vector valued over the $N$ coordinates $\bfx_1, \dots, \bfx_N$, $\bfJ(\bftheta(t, \bfmu))$ is the Jacobian matrix of $\PhiDECoder^{\uparrow}$. 

The time derivative $\dot{\bftheta}(t, \bfmu)$ is then obtained by minimizing the norm of the vector-valued residual function,
\begin{equation}\label{eq:LeeCIRM}
\dot{\bftheta}(t, \bfmu) \in \operatorname*{arg\,min}_{\bfeta \in \mathbb{R}^n} \|\bfr(\bftheta(t, \bfmu), \bfeta; \bfmu)\|_2^2\,.
\end{equation}
\citet{LeeC20} focus on implicit time discretizations, which lead to nonlinear least-squares problems that have to be solved numerically at each time step.  Constraints can be added to the residual minimization problem \eqref{eq:LeeCIRM} (or its discrete counterpart) to preserve quantities such as mass, momentum, and energy \cite{Lee_Carlberg_2021}.

The instantaneous residual minimization approach described by \citet{LeeC20} is closely related to the Dirac--Frenkel variational principle, except that \citet{LeeC20} directly work in the $N$ coordinates corresponding to the semi-discrete full model.

\paragraph{Online costs and empirical interpolation for autoencoder-based nonlinear parametrizations}
Achieving runtime speedups in the online phase compared to the full model is challenging with autoencoder-based nonlinear parametrization as used in \citet{LeeC20}. Although the optimization variable representing~$\dot{\bftheta}(t, \bfmu)$ is low-dimensional with dimension $n\ll N$, the residual itself is defined over $\RR^{N}$. Consequently, the  minimization problem~\eqref{eq:LeeCIRM} is carried out in the full-model state space $\mathbb{R}^N$. 
This issue becomes particularly challenging when an implicit time discretization is employed, leading to a nonlinear optimization problem that must be solved at each time step. For instance, if the corresponding nonlinear optimization problem is solved with Gauss--Newton (or similar) methods, then at each optimization iteration, the weight $\bftheta(t, \bfmu)$ is lifted to a full-model state approximation of dimension $N$ with the decoder $\PhiDECoder^{\uparrow}$ and then the residual needs to be computed with the full-model right-hand side function $\bff$ at all $N$ components. Additionally, the Jacobian matrix is computed corresponding to all $N$ spatial coordinates over which the full model is defined. 

The issue of lifting back to the high-dimensional full-model state space to compute the residual is analogous to the lifting bottleneck in linear model reduction, when the full-model dynamics are nonlinear in the state; see \Cref{sec:Prelim:EIM}. 
Correspondingly, there is work that has developed autoencoder parametrizations that are compatible with empirical interpolation-like techniques. The work by \citet{KIM2022110841} proposes using a shallow, masked decoder function that can be combined with sparse sampling methods analogous to empirical interpolation. 
The work \citet{Romor2023} goes a step further and proposes two complementary mechanisms to reduce online costs. First, the time-discrete residual minimization is performed only on a selected set of a few collocation points, which builds on the over-collocation approach introduced in \citet{CHEN2021110545} to sparsely evaluate the decoder function. Second, a compressed decoder is used that evaluates to values at the set of collocation points and so avoids evaluations at all $N$ spatial coordinates corresponding to the full-model state in the online phase.

\begin{remark}
There are methods that train neural-network autoencoders to represent solution fields and additionally learn the reduced dynamics from data, rather than computing the online weights via instantaneous residual minimization \cite{XU2020113379,10.1063/5.0039986,FRESCA2022114181}. We give a brief outlook to such non-intrusive model reduction methods in \Cref{sec:Outlook}. 
\end{remark}

\subsubsection{Autoencoders with a priori feature maps and quadratic approximations}\label{sec:PreTrainQM}
We now survey nonlinear model reduction techniques based on autoencoders with decoders that are formulated via a priori chosen, fixed nonlinear feature maps and encoder functions that are linear. In particular, we focus on polynomial feature maps, and especially quadratic feature maps, which lead to so-called quadratic (trial) manifolds.

\paragraph{Autoencoders with a priori chosen, fixed feature maps}
From \Cref{sec:PreTrainAutoEncoder}, recall that an autoencoder consists of an encoder $\PhiDECoder^{\downarrow}(\cdot; \woff^{\downarrow}): \mathbb{R}^{N} \to \mathbb{R}^n$ and a decoder $\PhiDECoder^{\uparrow}(\cdot; \woff^{\uparrow}): \mathbb{R}^{n} \to \mathbb{R}^{N}$, which depend on trainable offline weights $\woff^{\downarrow}$ and $\woff^{\uparrow}$, respectively. Recall further that we can interpret the decoder $\Phi^{\uparrow}: \mathbb{R}^n \to \mathbb{R}^N$ as inducing a vector-valued nonlinear parametrization as in \eqref{eq:DecoderAsParam}. 

We now consider decoders $\PhiDECoder^{\uparrow}$, and thus equivalently vector-valued nonlinear parametrizations, that have the following form
\begin{equation}\label{eq:QMFormHatBFQ}
\PhiDECoder^{\uparrow}(\bftheta; \woff^{\uparrow}) = \bfVa\bftheta + \bfVb \featuremap(\bftheta)\,,
\end{equation}
where $\bfVa \in \mathbb{R}^{N \times n}$ and $\bfVb \in \mathbb{R}^{N \times \gdim}$ are matrices that are trainable offline, i.e., $\woff^{\uparrow} = [\operatorname{vec}(\bfVa); \operatorname{vec}(\bfVb)]$. The function $\featuremap\colon \mathbb{R}^n \to \mathbb{R}^{\gdim}$ is a nonlinear feature map function that is given a priori. 
The rationale behind ansatz~\eqref{eq:QMFormHatBFQ} is twofold: (i) enrich the linear subspace spanned by the columns of $\bfVa$ by nonlinear terms (as in a truncated Taylor expansion) and (ii) for certain applications, \citet{Cohen2023NonlinearCB} prove that coefficients of higher-order POD modes can be expressed as functions of lower-order ones.

A feature map $\featuremap$, which is common in the context of nonlinear model reduction, is the quadratic function that collects all unique degree-2 monomials
\begin{equation}\label{eq:QuadraticFeatureMap}
\featuremap(\bftheta) = \bftheta \hat{\otimes} \bftheta, 
\end{equation}
where $\hat{\otimes}$ denotes the symmetric (duplicate-free) Kronecker product, i.e., the vector containing the products $\{\theta_i\theta_j\}_{1 \leq i \leq j \leq n}$ of all $n$ components of $\bftheta$ in a fixed ordering. In case of such a quadratic feature map, the number of features is $\gdim = n(n + 1)/2$. We can interpret the two terms of \eqref{eq:QMFormHatBFQ} as follows: The first term $\bfVa\bftheta$ provides a linear approximation and the second term $\bfVb\featuremap(\bftheta)$ is a nonlinear correction. A nonlinear parametrization of the form \eqref{eq:QMFormHatBFQ} gives rise to a trial manifold, which is sometimes called quadratic manifold if the feature map is quadratic. 

Quadratic and polynomial approximations are pursued by  \citet{JainTRR2017quadratic,RutzmoserRTJ2017Generalization,Cenedese2022,BarnettF2022Quadratic,GeelenWW2023Operator} and  have led to a series of publications on quadratic approximations for nonlinear model reduction, including on methods for training quadratic approximations on snapshot data \cite{GeelenBW2023Learning,GeelenBWW2024Learning,schwerdtner2024greedyconstructionquadraticmanifolds,10.1098/rspa.2024.0670} and efficiently using them in the online phase \cite{https://doi.org/10.1002/pamm.202200049,SHARMA2023116402,GOYAL2024134158,Yldz_2024,WEDER2025114249,GlaM25}. The earlier work \citet{Gu:EECS-2012-217} also discusses nonlinear parametrizations based on polynomial feature maps in the context of model reduction. In general, we remark that quadratic systems and quadratic dynamics as well as leveraging the corresponding tensor structure are omnipresent in model reduction; see, e.g., \citet{5991229,doi:10.1137/14097255X,Schlegel_Noack_2015,doi:10.1137/16M1098280,doi:10.2514/1.J057791,QIAN2020132401}. An extension to rational manifolds is presented by \citet{KleSCPH24}.

The expressivity of nonlinear parametrizations of the form \eqref{eq:QMFormHatBFQ} is dominated by the a priori chosen feature map $\featuremap$. \citet{BuchfinkGlasHaasdonk2024} provide lower bounds on the Kolmogorov $n$-width for such nonlinear parametrizations with polynomial feature maps $\featuremap$. In particular, roughly speaking, \citet{BuchfinkGlasHaasdonk2024} show that for the linear advection example studied in \Cref{sec:KolmogorovNonLinExamples}, the lower bound of the $n$-width behaves as $\mathcal{O}(n^{-\alpha/2})$, where $\alpha$ is the degree of the polynomial. 

Relying on approximation theoretic results regarding stable manifold widths and sensing numbers, see \Cref{sec:NonLinWidth}, and popularized by \citet{BarnettF2022Quadratic,GeelenWW2023Operator,Cohen2023NonlinearCB}, it is a common choice to use a linear encoder $\PhiDECoder^{\downarrow}(\cdot; \woff^{\downarrow})\colon \mathbb{R}^N \to \mathbb{R}^n$, i.e.,   
\begin{equation}\label{eq:IRN:LinearEncoderPolyManifold}
    \PhiDECoder^{\downarrow}(\bfq) = \bfW^{\top}\bfq\,,
\end{equation}
for some matrix $\bfW\in\RR^{N\times n}$. While general conditions for the matrix $\bfW$ can be inferred from the general framework presented in \citet{BUCHFINK2024134299}, a common choice is to use $\bfW^{\top} = \bfVa^{+}$, corresponding to the Moore--Penrose pseudo-inverse of the matrix $\bfVa$. Note that if $\bfVa$ has orthonormal columns, then $\bfVa^+ = \bfVa^{\top}$. Linear encoders are also used with more generic nonlinear decoders than quadratic ones in \citet{Cohen2023NonlinearCB,doi:10.2514/6.2023-0535,BARNETT2023112420,Ballout2024}. In particular, \citet{Cohen2023NonlinearCB} investigate such linear encoder and nonlinear decoder pairs in terms of the sensing number and show that the approximation error for functions representing transported features can decay substantially faster than the Kolmogorov $n$-width with respect to the dimension $n$. A similar strategy via a linear encoder and a nonlinear decoder is taken in \citet{bensalah2025nonlinearmanifoldapproximationusing}. Furthermore, it is also observed empirically that linear encoders achieve similar reconstruction errors as more generic nonlinear encoders with the additional benefit of an improved training robustness \cite{GLAS2025550,doi:10.1137/24M1717270} .

\paragraph{Leveraging the structural form for efficient training on snapshot data}
Recall the linear encoder \eqref{eq:IRN:LinearEncoderPolyManifold} with the particular choice $\bfW^{\top} = \bfVa^{+}$ and nonlinear decoder pair \eqref{eq:QMFormHatBFQ} with the trainable matrices $\bfVa \in \mathbb{R}^{N \times n}$ and $\bfVb \in \mathbb{R}^{N \times \gdim}$, which constituted the offline weights $\woff^{\downarrow} = \operatorname{vec}(\bfVa)$ and $\woff^{\uparrow} = [\operatorname{vec}(\bfVa); \operatorname{vec}(\bfVb)]$. Recall further that the feature map $\featuremap$ is given a priori. 

A typical training strategy to train such nonlinear parametrizations is to first choose the matrix $\bfVa$ and then to fit $\bfVb$ with a linear least-squares problem,
\begin{equation}\label{eq:QM:FittingW2}
\bfVb= \operatorname*{arg\,\min}_{\bfZ \in \mathbb{R}^{N \times \gdim}} \|\bfQ_{\text{train}} - \bfVa\underbrace{\PhiDECoder^{\downarrow}(\bfQ_{\text{train}}; \woff^{\downarrow})}_{\bfVa^{+}\bfQ_{\text{train}}} - \bfZ\featuremap(\underbrace{\PhiDECoder^{\downarrow}(\bfQ_{\text{train}}; \woff^{\downarrow})}_{\bfVa^{+}\bfQ_{\text{train}}})\|_F^2 + \lambda \|\bfZ\|_F^2\,,
\end{equation}
which corresponds to the reconstruction error $\bfQ_{\text{train}} - \PhiDECoder^{\uparrow}(\PhiDECoder^{\downarrow}(\bfQ_{\text{train}}; \woff^{\downarrow}); \woff^{\uparrow})$ of  the snapshots  \eqref{eq:QtrainSnapshotMatrix} given as columns in the snapshot matrix $\bfQ_{\text{train}}$. The encoder and feature map are applied column-wise to the snapshot matrix. A regularizer is added to the least-squares problem with regularization parameter $\lambda > 0$. 

Let us now consider fitting the matrix $\bfVa$, which fully determines the encoder but also enters in the linear term in the decoder. The works \citet{BarnettF2022Quadratic,GeelenWW2023Operator} follow the interpretation that $\PhiDECoder^{\uparrow}(\bftheta; \woff^{\uparrow}) = \bfVa\bftheta + \bfVb \featuremap(\bftheta)$ as defined in \eqref{eq:QMFormHatBFQ} is a linear approximation $\bfVa\bftheta$ with a nonlinear correction term $\bfVb\featuremap(\bftheta)$. Thus, \citet{BarnettF2022Quadratic,GeelenWW2023Operator}  fit $\bfVa$ so that the corresponding linear term $\bfVa\bftheta$ leads to the lowest linear approximation error on the snapshot data by setting the columns of the matrix $\bfVa \in \mathbb{R}^{N \times n}$ to the first $n$ leading left-singular vectors of the snapshot matrix $\bfQ_{\text{train}}$.  
However, \citet{Cohen2023NonlinearCB} and \citet{schwerdtner2024greedyconstructionquadraticmanifolds} stress the perspective that the nonlinear correction term evaluates the feature map $\featuremap$ at the encoded point $\bftheta = \PhiDECoder^{\downarrow}(\bfq; \woff^{\downarrow}) = \bfVa^+\bfq$. Thus, purely fitting $\bfVa$ to minimize the approximation error of the linear approximation can lead to an encoding $\PhiDECoder^{\downarrow}(\bfq; \woff^{\downarrow})$ that discards information that is crucial for the feature map $\featuremap$ to be informative and to provide an effective correction. 
To address this, \citet{schwerdtner2024greedyconstructionquadraticmanifolds} introduce a greedy strategy to build $\bfVa$ column-by-column. Rather than fixing $\bfVa$ to have columns corresponding to the first $n$ left-singular vectors (principal components) of the snapshot matrix $\bfQ_{\text{train}}$, the greedy method introduced in \citet{schwerdtner2024greedyconstructionquadraticmanifolds} selects $n$ basis vectors from a larger pool of leading and later principal components that can have index greater than $n$. At each step of the greedy procedure, a vector that yields the largest improvement in terms of reconstruction, after fitting $\bfVb$ with the linear least-squares problem \eqref{eq:QM:FittingW2}, is chosen. The resulting~$\bfVa$ can therefore include non-leading left singular vectors (i.e., with index greater than $n$) if they are important for making the feature map effective in terms of reducing the residual norm in the objective \eqref{eq:QM:FittingW2}. 

In another line of work, \citet{GeelenBW2023Learning,GeelenBWW2024Learning} formulate optimization problems that construct $\bfVa$ and $\bfVb$ together, which can also result in matrices that achieve lower approximation errors than defining the matrix $\bfVa$ via the $n$ leading singular vectors of the snapshot matrix.

\paragraph{Leveraging the structural form for efficient online phase}
A key advantage of quadratic (or polynomial) approximations of the form \eqref{eq:QMFormHatBFQ} is that they are highly structured, which can be leveraged in the residual minimization in the online phase.

Recall the semi-discrete least-squares problem \eqref{eq:NG:EmpiricalMatrixObj} corresponding to instantaneous residual minimization with the Dirac--Frenkel variational principle. Recall further that the batch gradient becomes a Jacobian matrix because $\PhiDECoder^{\uparrow}: \mathbb{R}^n \to \mathbb{R}^N$ is a vector-valued nonlinear parametrization. For parametrizations with quadratic feature maps \eqref{eq:QuadraticFeatureMap}, the Jacobian matrix $\bfJ(\bftheta(t, \bfmu)) \in \mathbb{R}^{N \times n}$ of $\PhiDECoder^{\uparrow}$ 
admits a closed-form expression and is affine in  $\bftheta(t, \bfmu)$. More precisely, because  $\partial(\bftheta\hat{\otimes}\bftheta)/\partial\bftheta$ 
is linear in $\bftheta$, one can write
\begin{equation}\label{eq:BatchGradientStructureQM}
\bfJ(\bftheta(t, \bfmu))=\bfVa+\sum_{\ell=1}^n \theta_\ell(t, \bfmu)\,\bfJ_\ell,
\end{equation}
for matrices $\bfJ_\ell\in\mathbb{R}^{N\times n}$ that can be pre-computed in the offline phase independently of the online weight vector $\bftheta(t, \bfmu)$. 

When the full model is linear in the state and has an affine parameter dependence, then the objective of the least-squares problem for instantaneous residual minimization can be assembled in the online phase with no computations that scale with the ambient dimension $N$; this offline--online separation and $N$-independent online assembly is shown and leveraged in \citet{SHARMA2023116402} and \citet{WEDER2025114249}. When the full model has a nonlinear state dependence, the pre-computation of the Jacobian components can still be useful. In particular, \citet{WEDER2025114249} propose a quadratic approximation formulation that is not tied to the $N$ full-model spatial coordinates but can be realized on a set of selected collocation points, which allows circumventing more expensive online costs. 

Yet another way to leverage the quadratic structure is proposed by \citet{BarnettF2022Quadratic}, who plug the quadratic parametrization into an implicit-in-time least-squares Petrov--Galerkin discretization. Because the quadratic parametrization has a closed-form affine Jacobian, sparse sampling (analogous to empirical interpolation; see \Cref{sec:Prelim:EIM}) can be tailored to it and efficiently applied. As a result, online residual and gradient evaluations are restricted to a set of a few selected points, which leads to online runtime speedups compared to the full model in the presented  numerical experiments.

The work \citet{GeelenWW2023Operator} goes even a step further in leveraging structure and approximates the Jacobian \eqref{eq:BatchGradientStructureQM} with its first-order approximation, which is just the constant function that evaluates to $\bfVa$. In this case, the Jacobian matrix is independent of the online weight $\bftheta(t, \bfmu)$, which leads to faster reduced-model evaluations in the online phase; see \citet{WEDER2025114249} for a comparison in terms of online runtime. Furthermore, the work \citet{GeelenWW2023Operator} demonstrates that approximating the Jacobian matrix with its first-order approximation $\bfVa$ is a key step to derive a non-intrusive reduced modeling approach based on quadratic manifolds; see \Cref{sec:Outlook} for a more detailed discussion about intrusive versus non-intrusive model reduction.

The structure offered by quadratic (and other a priori fixed feature-map approximations) in the Jacobian \eqref{eq:BatchGradientStructureQM} is in stark contrast to the gradients corresponding to generic autoencoder parametrizations such as neural networks that have no comparable low-order structure and so must be obtained repeatedly online via automatic differentiation, possibly within each Gauss-Newton iteration at each time step. We stress though that a priori fixed feature maps lead often to nonlinear parametrizations with lower expressivity than more generic architectures given by deep neural networks with representation learning. Again, we refer to \citet{BuchfinkGlasHaasdonk2024} for lower approximation bounds in this context.

\subsubsection{Neural representations with offline--online weights}\label{sec:IRM:NeuralRep}
We now consider neural representations, which use a neural network to directly represent the nonlinear parametrization as $\hat{q}$.

\paragraph{Neural representations of solution functions in model reduction}
Neural-network parametrizations of solution fields are sometimes referred to as implicit neural representations, which originated in the computer graphics community \cite{NEURIPS2020_53c04118} but have also been used for reducing (compressing) physics solution fields, see, e.g., \citet{brunton_implcitflow}. 

A key feature of neural representations is that they can be trained on snapshot data coming from different spatial discretizations because they represent fields as functions of spatial coordinates, which is in contrast to typical neural-network autoencoders as discussed above that depend on a fixed spatial discretization. This is particularly advantageous when full-model solutions are computed on adaptive meshes or point clouds, which is leveraged in the context of nonlinear model reduction in \citet{chen_crom_2023}. 

To use neural representations in model reduction, an important design choice is how to decompose the weights of the neural network into an offline weight vector $\woff$ that can be trained offline once for all snapshot data, and an online weight vector $\bftheta(t, \bfmu)$ with few components that are obtained online with, e.g., instantaneous residual minimization. Based on an offline--online decomposition of the weights, one can represent the nonlinear parametrization as $\hat{q}(\bftheta(t, \bfmu), \cdot; \woff): \Omega \to \mathbb{R}$, which is the form given in \eqref{eq:NonLin:ParamWithOfflineWeight}, 
where $\bftheta(t, \bfmu)$ correspond to the online weights of the network that may vary with time $t$ and parameter $\bfmu$. The input to the neural network is the spatial coordinate in $\Omega$. And the offline weights $\woff$ are independent of time $t$ and the parameter $\bfmu$ and fixed during the online phase.

\paragraph{Neural representations with offline--online weights}
One approach to achieve an offline--online decomposition of the neural network weights is to let the online weights enter via modulations to the offline weights. A common weight modulation approach is via low-rank matrices introduced in \cite{hu_lora_2021}, which leads to low-rank adaptation layers. For model reduction, continuous low-rank adaptation (CoLoRA) layers have been introduced by \citet{berman2024colora}, where the modulation can depend continuously on time and a parameter $\bfmu$. 

A CoLoRA layer has the form 
\begin{equation}\label{eq:CoLoRALayers}
\mathcal{C}_{\ell}(\bfy) = \bfW_{\ell}\bfy + \alpha_{\ell}(t, \bfmu)\bfA_{\ell}\bfB_{\ell}\bfy + \bfb_{\ell}\,,
\end{equation}
where $\bfW_{\ell}$ is a weight matrix and $\bfb_{\ell}$ is a bias term. Additionally, there is a low-rank matrix $\bfA_{\ell}\bfB_{\ell}$ of rank $r$ with the rank-$r$ factors $\bfA_{\ell}$ and $\bfB_{\ell}$ and the modulation $\alpha_{\ell}(t, \bfmu) \in \mathbb{R}$.  
The layers \eqref{eq:CoLoRALayers} can be composed into CoLoRA networks as
\begin{equation}\label{eq:CoLoRANetwork}
\hat{q}(\bftheta(t, \bfmu), \bfx; \woff) = \mathcal{C}_L(\sigma(\mathcal{C}_{L-1}(\dots \sigma(\mathcal{C}_1(\bfx)))))\,,
\end{equation}
which is a nonlinear parameterization that has the form \eqref{eq:NonLin:ParamWithOfflineWeight}. 
The offline weight vector $\woff$ consists of all weights in the network \eqref{eq:CoLoRANetwork} that are independent of time~$t$ and parameter $\bfmu$, including all weight matrices $\bfW_1, \dots, \bfW_L$, low-rank matrices $\bfA_1, \dots, \bfA_{L}, \bfB_1, \dots, \bfB_L$, and bias terms $\bfb_1, \dots, \bfb_L$.
The online weight vector consists of the modulations
\[
\bftheta(t, \bfmu) = \begin{bmatrix} \alpha_1(t, \bfmu), \dots, \alpha_L(t, \bfmu)\end{bmatrix}^{\top} \in \mathbb{R}^L\,.
\]
In the particular architecture of CoLoRA networks given in \eqref{eq:CoLoRANetwork}, the number of layers $L$ determines the number of online weights $n$ as $n = L$. There are other variants of CoLoRA networks, e.g., replacing some of the $L$ CoLoRA layers with regular linear layers that do not have a modulation and so one can achieve $n < L$ online weights. Furthermore, one can introduce multiple modulations per layer by introducing a diagonal matrix to obtain
\[
\mathcal{C}_{\ell}(\bfy) = \bfW_{\ell}\bfy + \bfA_{\ell} \diag(\alpha_{\ell}^{(1)}(t, \bfmu), \dots, \alpha_{\ell}^{(r)}(t, \bfmu)) \bfB_{\ell}\bfy + \bfb_{\ell}\,,
\]
where the low-rank factors $\bfA_{\ell}$ and $\bfB_{\ell}$ have rank $r$. Such variants are discussed in \citet{berman2024colora}. In general, there a many ways how to modulate the weights \cite{10.5555/3504035.3504518}.

There are also other ways of formulating neural networks with offline and online weights, besides weight modulation. For example, the work \citet{chen_crom_2023} considers multilayer perceptrons to represent $\hat{q}$ and let a latent state enter as input to the network. The latent state is analogous to what we refer to as online weight vector $\bftheta(t, \bfmu)$, except that it enters as input. The weights of the multilayer perceptron are the offline weights in our terminology.

\paragraph{Training offline weights of neural representations on snapshot data}
A direct training of the offline weight vector $\woff$ and the online weight vector $\bftheta(t, \bfmu)$ over all times and parameters would allow the online weights to vary arbitrarily over time, which can lead to non-smooth dependence of the online weight on time.  To encourage regularity on the online weights in time $t$, \citet{berman2024colora} use a so-called modulation (or hyper-) network during training: Let $h_{\bfpsi}: [0, T] \times \Dcal \to \mathbb{R}^{n}$ be a neural network that maps the time $t$ and the parameter $\bfmu$ onto an $n$-dimensional vector $\bftheta(t, \bfmu)$. The trainable weight vector of the network $h_{\bfpsi}$ is $\bfpsi \in \mathbb{R}^{n'}$, which is typically of low dimension.

Using the network $h_{\bfpsi}$ leads to the training problem 
\begin{equation}\label{eq:CoLoRAOfflineTrainingProb}
\min_{\woff, \bfpsi} \frac{1}{\ntmu{\spts}(K+1)}\sum_{i = 1}^{\ntmu} \sum_{j = 1}^{\spts} \sum_{k = 0}^K \frac{|q(t_k, \bfx_j; \bfmu_i) - \hat{q}(h_{\bfpsi}(t_k, \bfmu_i); \bfx_j; \woff)|^2}{|q(t_k, \bfx_j; \bfmu_i)|^2}
\end{equation}
over the offline weights $\woff$ and the weights $\bfpsi$ of the network $h_{\bfpsi}$. 
Problem~\eqref{eq:CoLoRAOfflineTrainingProb} depends on ${\ntmu}$ training parameters $\bfmu_1, \dots, \bfmu_{\ntmu}$ and $K + 1$ times steps $0 = t_0 < t_1 < \dots < t_K = T$, as in the training data \eqref{eq:QtrainSnapshotMatrix}. Additionally, we use ${\spts}$ spatial points $\bfx_1, \dots, \bfx_{\spts} \in \Omega$. 
Recall that $q(t_k, \bfx_j; \bfmu_i)$ is the full-model solution \eqref{eqn:PDEapprox:FOM} at time $t_k$ and parameter $\bfmu_i$ evaluated at the spatial coordinate $\bfx_j$. In particular, when training neural representations, one can directly use the full-model solution function $q$ instead of the semi-discrete solution vector, which emphasizes the independence of the neural representation from the mesh and discretization used for the full model. 

The modulation network $h_{\bfpsi}$ can be seen as playing an analogous role to the encoder function in autoencoder-based parametrizations. The major difference however is that the modulation network receives as inputs the time $t$ and the parameter $\bfmu$ (and possibly other inputs) but not the full-model state vector.

\paragraph{Online residual minimization with trained neural representations}
Once trained, the neural representation can be combined with instantaneous residual minimization in the online phase. 

In stark contrast to previously considered nonlinear parametrizations that are formulated with respect to semi-discrete full models, which evaluate to vectors of dimension $N$ (see \Cref{sec:PreTrainAutoEncoder} and \Cref{sec:PreTrainQM}), neural representations are coordinate-to-value maps and exactly of the form of how we defined nonlinear parametrizations in \eqref{eq:NMOR:Param} (and with offline weights in \eqref{eq:NonLin:ParamWithOfflineWeight}). Consequently, when applying schemes based on the Dirac--Frenkel variational principle for residual minimization, the inner products and residual norms can be numerically estimated with a set of ${\spts}$ sampling points, where ${\spts}$ can be chosen independently of the full-model dimension $N$ and its mesh. In particular, \citet{berman2024colora} leverage Neural Galerkin schemes in the online phase. Additionally, quantities such as mass, momentum, and energy can be conserved by adding constraints and nonlinear embedding steps to the instantaneous residual minimization \cite{SSBP23NGE}.

At each time step in the online phase, a single linear least-squares problem has to be solved if an explicit time integration scheme is used \cite{berman2024colora}. The number of unknowns is $n$, with the aim of keeping $n$ small when training the parametrization on snapshot data. The online costs are therefore primarily governed by the number $\spts$ of sampling points used to approximate the residual norm in the least-squares problem and by computing the batch gradient and evaluating the full-model right-hand side function at the sampling points. In particular, the full-model right-hand side function can contain spatial derivatives, which can be obtained via automatic differentiation through the neural representation at the sampling points. 

Because the online phase operates on the $\spts$ sampling points rather than on the $N$ points corresponding to the $N$-dimensional state vector of the semi-discrete full model, it avoids the classical lifting bottleneck that is typical for reduced models formulated over the semi-discrete full model. At the same time, accuracy depends on selecting sufficiently many and sufficiently informative sampling points. In particular, while the number $\spts$ of sampling points and the number $N$ of full-model grid points are formally independent, they are intrinsically related in the sense that if the full model needs a fine grid to resolve fine-scale solution features, then also the number of sampling points $\spts$ often has to reflect this, even though adaptivity can be helpful to keep $\spts$ small, which motivates adaptive sampling strategies as discussed in the context of Neural Galerkin schemes in \Cref{sec:NGWithAdaptiveSamplingPoints}.

\begin{remark}
An alternative to online residual minimization is to reuse the trained modulation network $h_{\bfpsi}$ in the online phase to compute the online weights as $\bftheta(t, \bfmu) = h_{\bfpsi}(t, \bfmu)$, thereby obtaining the solution approximation by only evaluating a neural network. Since the modulation network $h_{\bfpsi}$ depends only on known inputs such as $t$ and $\bfmu$ (rather than on the full state as an encoder function of an autoencoder), the resulting online phase is fully non-intrusive and data-driven, but it looses the Galerkin optimality associated with Dirac--Frenkel residual minimization. Related approaches for computing weights online via modulation networks or hyper-networks for non-intrusive model reduction are considered in, e.g., \citet{belbute-peres2021hyperpinn,cho2023hypernetworkbased,berman2024colora}; see \Cref{sec:Outlook} for further discussion.
\end{remark}

\section{Conclusions and outlook}\label{sec:Outlook}
Significant progress has been made on model reduction for transport-dominated problems, as reflected by the large and rapidly growing body of literature. A large class of nonlinear model reduction methods can be viewed through the unifying lens of generic nonlinear parameterizations combined with instantaneous residual minimization. At the same time, we deliberately distinguished transformation-based and online adaptive methods from this generic perspective, as this classification highlights additional structural features that can be leveraged for improved efficiency and robustness.

Despite the breadth of the work cited in this survey, several important research directions remain underrepresented, including structure-preserving model reduction, the reduction of stochastic systems, and model reduction for control applications. An entirely different line of research, still in its early stages, recasts hyperbolic PDEs as delay equations and then applies model reduction to this alternative representation. We refer to \citet{SchUBG18} for some first developments in this direction. 

A particularly fruitful and rapidly evolving area is non-intrusive model reduction for transport-dominated problems. In contrast to the intrusive methods emphasized in this survey, non-intrusive model reduction approaches aim to learn both an efficient low-dimensional representation and the corresponding reduced dynamics primarily from snapshot data, rather than directly deriving the reduced dynamics from the governing equations. We emphasize that non-intrusive model reduction is distinct from equation discovery, which uses data to infer the governing physical laws themselves. A large literature on non-intrusive model reduction exists; see the works  \citet{Ghattas_Willcox_2021,annurev:/content/journals/10.1146/annurev-fluid-121021-025220} for an overview. 

Several non-intrusive methods explicitly target transport-dominated regimes where linear model reduction is inefficient. 
One line of work incorporates nonlinear transformations that factor out transport before applying linear non-intrusive model reduction methods \cite{SesS19,SARNA2021114168,LU2020109229,Issan2023,GOWRACHARI2025106758}, analogous to the transformation-based intrusive methods that we reviewed above.   
Another class of non-intrusive methods learns a nonlinear low-dimensional representation (e.g., via autoencoders) and then learns the reduced dynamics in the latent representation using a separate polynomial or neural-network model  \cite{Fresca2021,10.1063/5.0039986,FRESCA2022114181,Franco2023,Cenedese2022,GeelenWW2023Operator,Regazzoni2024,Tomasetto2025}. We stress that the references cited in this outlook provide only a selective glimpse into the rapidly expanding field of non-intrusive model reduction. They are intended to illustrate the   opportunities of non-intrusive methods for nonlinear model reduction and not to provide a comprehensive review.

\section*{Acknowledgments}  
B.P.~was partially supported by the National Science Foundation grant 2046521 and the Air Force Office of Scientific Research (AFOSR), USA, award FA9550-24-1-0327. B.U.~is  funded by Deutsche Forschungsgemeinschaft (DFG, German Research Foundation) – Project-ID 258734477 – SFB 1173, and by the BMBF (grant no.~05M22VSA).

\addcontentsline{toc}{section}{References}
\bibliography{journalAbbr,literature}
\label{lastpage}
\end{document}